\documentclass[12pt,leqno]{article}
\usepackage{amsmath,amssymb,amscd,latexsym}

\usepackage[all]{xy}
\newcommand{\A}{{\mathbb{A}}}
\newcommand{\C}{{\mathbb{C}}}
\newcommand{\F}{{\mathbb{F}}}
\newcommand{\Ge}{\mathbb{G}}
\newcommand{\Pa}{{\mathbb{P}}}
\newcommand{\Q}{{\mathbb{Q}}}
\newcommand{\oQ}{\overline{\Q}}
\newcommand{\Z}{{\mathbb{Z}}}
\newcommand{\oZ}{\overline{\Z}}
\newcommand{\Abb}{\mathrm{Alb}}
\newcommand{\abb}{\mathrm{ab}}
\newcommand{\Alb}{\mathrm{Alb}}
\newcommand{\cont}{\mathrm{cont}}

\newcommand{\et}{\mathrm{\acute{e}t}}
\newcommand{\ev}{\mathrm{ev}}
\newcommand{\of}{\overline{f}}
\newcommand{\fin}{\mathrm{fin}}
\newcommand{\good}{\mathrm{good}}
\newcommand{\id}{\mathrm{id}}
\renewcommand{\mod}{\;\mathrm{mod}\;}
\newcommand{\Mod}{\mathbf{Mod}}
\newcommand{\Mor}{\mathrm{Mor}}
\newcommand{\ob}{\mathrm{ob}}
\newcommand{\op}{\mathrm{op}}
\newcommand{\rank}{\mathrm{rank}}
\newcommand{\red}{\mathrm{red}}
\newcommand{\res}{\mathrm{res}}
\newcommand{\Spec}{\mathrm{Spec}\,}
\newcommand{\spec}{\mathrm{spec}\,}
\newcommand{\Aut}{\mathrm{Aut}}
\newcommand{\Ext}{\mathrm{Ext}}
\newcommand{\uExt}{\underline{\Ext}}
\newcommand{\Gal}{\mathrm{Gal}}
\newcommand{\GL}{\mathrm{GL}\,}
\newcommand{\Hom}{\mathrm{Hom}}
\newcommand{\uHom}{\underline{\Hom}}
\newcommand{\uIso}{\underline{\mathrm{Iso}}}
\newcommand{\Ker}{\mathrm{Ker}\,}
\newcommand{\Lie}{\mathrm{Lie}\,}
\newcommand{\Pic}{\mathrm{Pic}}
\newcommand{\Rep}{\mathbf{Rep}\,}
\renewcommand{\Vec}{\mathbf{Vec}}
\newcommand{\Zar}{\mathrm{Zar}}
\newcommand{\tors}{\mathrm{tors}}
\newcommand{\Ah}{{\mathcal A}}
\newcommand{\Bh}{\mathcal{B}}
\newcommand{\hAh}{\hat{\Ah}}
\newcommand{\oAh}{\overline{\Ah}}

\newcommand{\Ch}{{\mathcal C}}
\newcommand{\Fh}{{\mathcal F}}

\newcommand{\Ih}{{\mathcal I}}

\newcommand{\Lh}{{\mathcal L}}
\newcommand{\Oh}{{\mathcal O}}

\newcommand{\Wh}{\mathcal{W}}
\newcommand{\Yh}{\mathcal{Y}}
\newcommand{\Zh}{\mathcal{Z}}
\newcommand{\ea}{\mathfrak{a}}

\newcommand{\eo}{\mathfrak{o}}
\newcommand{\eB}{\mathfrak{B}}
\newcommand{\eT}{\mathfrak{T}}
\newcommand{\eX}{{\mathfrak X}}
\newcommand{\oeX}{\overline{\eX}}

\newcommand{\oK}{\overline{K}}
\newcommand{\ox}{\overline{x}}
\newcommand{\oX}{\overline{X}}
\newcommand{\oy}{\overline{y}}

\newcommand{\hA}{\hat{A}}
\newcommand{\ha}{\hat{a}}
\newcommand{\oA}{\overline{A}}

\newcommand{\ohne}{\setminus}
\newcommand{\silo}{\stackrel{\sim}{\longrightarrow}}
\newcommand{\tei}{\, | \,}
\newcommand{\ent}{\;\widehat{=}\;}
\newcommand{\hullet}{\raisebox{0.05cm}{$\scriptscriptstyle \bullet$}}
\newcommand{\verk}{\mbox{\scriptsize $\,\circ\,$}}

\newtheorem{theorem}{Theorem}
\newtheorem{lemma}[theorem]{Lemma}
\newtheorem{prop}[theorem]{Proposition}
\newtheorem{defn}[theorem]{Definition}
\newtheorem{cor}[theorem]{Corollary}

\newenvironment{rem}{\noindent {\bf Remark}}{}
\newenvironment{rems}{\noindent {\bf Remarks}}{}
\newenvironment{definition}{\noindent {\bf Definition}\it}{}
\newenvironment{proofof}{\noindent {\bf Proof of}}{\mbox{}\hfill$\Box$}
\newenvironment{proof}{\noindent {\bf Proof}}{\mbox{}\hfill$\Box$}
\begin{document}
\title{Vector bundles and $p$-adic representations I}
\author{Christopher Deninger \and Annette Werner}
\date{\ }
\maketitle
\section{Introduction}
On a compact Riemann surface every finite dimensional complex representation of the fundamental group gives rise to a flat vector bundle and hence to a holomorphic vector bundle. By a theorem of Weil, one obtains precisely the holomorphic bundles whose indecomposable components have degree zero \cite{W}.

The present paper deals with a $p$-adic analogue of this construction. 

Contrary to the classical case we start with a vector bundle and not with a representation. Namely we define a certain category of vector bundles on a $p$-adic curve to which we can associate in a functorial way finite dimensional $p$-adic representations of the algebraic fundamental group.

More precisely, let $X$ be a smooth projective curve over a finite extension $K$ of $\Q_p$. Set $\oX = X \otimes_K \oK$ and let $\pi_1 (\oX , \ox)$ be the algebraic fundamental group of $\oX$ with respect to a base point $\ox$. Our category $\eB_{X_{\C_p}}$ is a full subcategory of the category of vector bundles on $X_{\C_p} = X \otimes_K \C_p$ defined by certain conditions on the reductions of the bundles $\mod p^n$ for all $n \ge 1$. It depends functorially on $X$ and is closed under direct sums, tensor products, duals, internal homs and exterior powers.

We prove that every vector bundle $E$ in $\eB_{X_{\C_p}}$ gives rise to a continuous representation $\rho_E$ of $\pi_1 (\oX , \ox)$ on the fibre of $E$ in the point of $X_{\C_p}$ induced by $\ox$. In fact, one obtains an exact additive functor $\rho$ from $\eB_{X_{\C_p}}$ to the category of continuous representations of $\pi_1 (\oX , \ox)$ on finite dimensional $\C_p$-vector spaces. The functor $\rho$ commutes with tensor products, duals, internal homs and exterior powers and it is compatible with the natural actions of the Galois group $G_K = \Gal (\oK / K)$ on the categories.

For example, if $X$ has good reduction then a line bundle $L$ on $X_{\C_p}$ lies in $\eB_{X_{\C_p}}$ precisely if it has degree zero. In this case $\rho_L$ is a continuous $\C^*_p$-valued character of $\pi_1 (\oX , \ox)$. For $L$ in a certain open subgroup of $\Pic^0_{X / K} (\C_p)$ this character was already constructed by Tate using the $p$-divisible group of the abelian scheme $\Pic^0_{\eX / \eo_K}$, where $\eX$ is a smooth model of $X$, cf. \cite{Ta} \S\,4. However, the method of Tate using Cartier duality of $p$-divisible groups cannot be generalized to vector bundles of higher rank on $X_{\C_p}$. 

Our category contains in particular all extensions of the trivial line bundle with itself. Since the functor $\rho$ is exact it induces a homomorphism $\rho_*$ of the corresponding extension groups. The $\Ext$ group in $\eB_{X_{\C_p}}$ is identified with $H^1 (X , \Oh) \otimes_K \C_p$, and the $\Ext$ group in the representation category of $\pi_1 (\oX , \ox)$ turns out to be $H^1_{\et} (\oX , \Q_p) \otimes \C_p$. We prove that under these identifications the homomorphism $\rho_*$ coincides with the inclusion
\[
i_{HT} : H^1 (X, \Oh) \otimes_K \C_p \hookrightarrow H^1_{\et} (\oX , \Q_p) \otimes \C_p
\]
coming from the Hodge--Tate decomposition.

Incidentially, this result leads to a new explicit description of $i_{HT}$ using unipotent vector bundles of rank two, see remark 1 after theorem \ref{t29}. 

On an elliptic curve with good reduction, it follows from Atiyah's classification that an indecomposable bundle lies in our category if and only if it has degree zero.

It seems possible that quite generally, a vector bundle lies in $\eB_{X_{\C_p}}$ if and only if its indecomposable components have degree zero.

For Mumford curves, Faltings associates a vector bundle on $X$ to every $K$-rational representation of the Schottky group and proves that every semistable vector bundle of degree zero arises in this way, cf. \cite{Fa1}. In the present paper we are primarily concerned with the good reduction case. It would be interesting to compare our constructions in the case of Mumford curves with Faltings' results.

Let us now briefly describe the contents of the individual  sections. 

In section 2, we consider an abelian variety $A$ over $K$ with good reduction and define a homomorphism
\[
\alpha : \Pic^0_{A / K} (\C_p) = \hA (\C_p) \longrightarrow \Hom_{\cont} (TA , \C^*_p) \; .
\]
This is essentially Tate's construction, who restricted himself to the $\eo = \eo_{\C_p}$-valued points of the $p$-divisible group of the dual abelian scheme. It follows from Tate's work that $d\alpha$ is the Hodge--Tate map
\begin{equation}
  \label{eq:1}
i_{HT} : H^1 (A , \Oh) \otimes_K \C_p \to H^1_{\et} (\oA , \Q_p) \otimes_{\Q_p} \C_p \; .  
\end{equation}
We define a certain subgroup of the $\C^*_p$-valued characters of $TA$ and prove that $\alpha$ is a topological isomorphism onto this group.

In section 3, we show that $\alpha$ is a homomorphism of $p$-adic Lie groups. We prove directly that $\Lie \alpha$ coincides with the Hodge--Tate map $i_{HT}$ in (\ref{eq:1}) in the explicit description by Faltings and Coleman via the universal vectorial extension. 

In section 4, we consider smooth and proper varieties $X$ over $K$ whose $H^1$ has good reduction.
Combining the results of sections 2 and 3 applied to the Albanese variety of $X$ with Kummer theory, we obtain an injective homomorphism of $p$-adic Lie groups
\[
\alpha^{\tau} : \Pic^{\tau}_{X/K} (\C_p) \longrightarrow \Hom_{\cont} (\pi^{\abb}_1 (\oX) , \C^*_p) 
\]
and determine its image and $\Lie \alpha^{\tau}$. 

The preceeding constructions of $\alpha$ and $\alpha^{\tau}$ are based on Cartier duality. In order to generalize from line bundles to vector bundles we develop the following alternative description of $\alpha^{\tau}$ on curves.

Given a line bundle $L_{\C_p}$ of degree zero on the curve $X_{\C_p}$, we can extend it to a line bundle $L$ on the smooth model $\eX_{\eo}$. For every $n \ge 1$ its reduction $L_n$ modulo $p^n$ has finite order, say $N = N (n)$. Using the $N$-multiplication on the Albanese variety $A$ of $X$, one obtains a finite, $A_N$-equivariant morphism $\pi : \Yh \to \eX_{\eo}$, which is generically Galois and such that $\pi^*_n L_n$ is trivial. Here $\pi_n$ denotes the reduction of the map $\pi$ modulo $p^n$. Furthermore it can be arranged that the structural morphism $\lambda : \Yh \to \spec \eo$ satisfies $\lambda_* \Oh_{\Yh} = \Oh_{\spec \eo}$ {\it universally}. Hence we have $\Gamma (\Yh_n , \pi^*_n L_n) \cong \eo_n = \eo / p^n \eo$ and the $A_N$-action on this space is given by a character $A_N \to \eo^*_n$. Passing to limits we obtain a character of $\pi^{\abb}_1 (\oX)$ with values in $\eo^*$ which turns out to be the same as $\alpha^{\tau} (L_{\C_p})$. 

This construction suggests the definition of a category $\eB_{\eX_{\eo}}$ of vector bundles on $\eX_{\eo}$ to which one can attach $p$-adic representations in a similar way. Namely in section 5 we define for every finitely presented flat and proper model $\oeX$ of a smooth projective curve $\oX / \oK$ a category $\eT$ of ``coverings''. The objects are finitely presented proper $G$-equivariant $\eo_{\oK}$-morphisms $\pi : \Yh \to \oeX$ where $G$ is a finite group which acts $\eo_{\oK}$-linearly on $\Yh$ and trivially on $\oX$. Moreover $\Yh_{\oK}$ is supposed to be an (\'etale) $G$-torsor over $\oX$. Using a result of Raynaud we prove that every such covering can be dominated by another one, say $\Yh' \to \oeX$ whose structural morphism $\lambda' : \Yh' \to \spec \eo_{\oK}$ is flat and satisfies $\lambda_* \Oh_{\Yh'} = \Oh_{\spec \eo_{\oK}}$ universally.

In section 6 we define $\eB_{\eX_{\eo}}$ where $\eX_{\eo} = \oeX \otimes \eo$ as the category of all vector bundles $E$ on $\eX_{\eo}$ such that for all $n \ge 1$ there is a ``covering'' $\pi : \Yh \to \oeX$ in $\eT$ with $\pi^*_n E_n$ a trivial bundle. Proceeding as in the case of line bundles one obtains a canonical representation $\rho_E$ of $\pi_1 (\oX , \ox)$ on the fibre $E_{x_{\eo}}$ of $E$ in the point $x_{\eo} \in \eX_{\eo} (\eo)$ induced by $\ox \in \oX (\oK)$. For an abelian scheme $\Ah / \eo_K$ we also define a category $\eB_{\Ah_{\eo}}$ in a similar way using only $N$-multiplications as coverings. We prove among other things that the map $E \mapsto \rho_E$ is functorial and compatible with pullback and the Galois action in case the curve is defined over $K$.

We can now define the category $\eB_{X_{\C_p}}$ mentioned above: It contains all vector bundles on $X_{\C_p}$ which are isomorphic to the generic fibre of some $E$ in $\eB_{\eX_{\eo}}$ for some model $\oeX$ of $\oX$. The rest of section 6 is mainly devoted to the proof of the statements from the beginning of this introduction. 

The present work raises a number of interesting questions which are briefly discussed in section 7.

We are very grateful to Peter Schneider for pointing out to us that by Tate's work it is more natural to attach $p$-adic representations to vector bundles than vice versa. We would also like to thank H\'el\`ene Esnault, Annette Huber, Mark Kisin, Damian Roessler and Matthias Strauch for interesting discussions and suggestions. 
%\newpage
%\input{sec2}
\section{Line bundles on abelian varieties and their characters}
This section is a variation on some of Tate's ideas in \cite{Ta}, \S\,4 on the pairing between the $\eo = \eo_{\C_p}$-valued points of a $p$-divisible group and the Tate module of its dual.

Let $A$ be an abelian variety over a finite extension $K$ of $\Q_p$. We assume that $A$ has good reduction and denote by $\Ah$ the abelian scheme over $\eo_K$ whose generic fibre is $A$. Let $\hA = \Pic^0_{A/K} \quad \mbox{and} \quad \hAh = \Pic^0_{\Ah / \eo_K}$ be the dual abelian schemes. We endow the $\Ge_m$-torsor associated to the Poincar\'e bundle over $\Ah \times \hAh$ with the structure of a biextension of $\Ah$ and $\hAh$ by $\Ge_m$, so that we have $\hAh = \uExt^1 (\Ah , \Ge_m)$ in the flat topology. We will now construct a $G_K = \Gal (\oK / K)$-equivariant homomorphism:
\begin{equation}
  \label{eq:2}
  \alpha = \alpha_A : \Pic^0_{A / K} (\C_p) = \hA (\C_p) \longrightarrow \Hom_c (TA , \C^*_p) \; .
\end{equation}
Here $\Hom_c$ denotes the group of {\it continuous} homomorphisms and $TA$ is the total Tate module of $A$. Because of the decomposition $\C^*_p = p^{\Q} \times \eo^*$ we have
\[
\Hom_c (TA , \C^*_p) = \Hom_c (TA , \eo^*) \; .
\]
Let $\oK$ be the algebraic closure of $K$ in $\C_p$ and set $\oA = A \otimes \oK$ and $\oAh = \Ah \otimes \eo_{\oK}$. Then we have $\hat{\oA} = \hA \otimes \oK$ and $\hat{\oAh} = \hAh \otimes \eo_{\oK}$. Finally for $n \ge 1$ set $\Ah_n = \oAh \otimes \eo_n$ where $\eo_n = \eo_{\oK} / p^n \eo_{\oK} = \eo / p^n \eo$.

For any $N \ge 1$ the exact sequence 
\[
0 \longrightarrow \Ah_N \longrightarrow \Ah \xrightarrow{N} \Ah \longrightarrow 0 \; ,
\]
where $\Ah_N$ denotes the subscheme of $N$-torsion points, leads to the exact sequence
\[
0 \longrightarrow \uHom (\Ah_N , \Ge_m) \longrightarrow \uExt^1 (\Ah , \Ge_m) \xrightarrow{N} \uExt^1 (\Ah , \Ge_m) \longrightarrow 0 
\]
in the flat topology.

Since $\uExt^1 (\Ah , \Ge_m) = \hAh$, this induces an isomorphism
\[
\hAh_N \silo \uHom (\Ah_N , \Ge_m) \; .
\]
Hence we get a homomorphism
$\hAh_N (\eo_n) \longrightarrow \Hom (\Ah_N (\eo_n) , \eo^*_n)$. Using the reduction map $A_N (\oK) = \Ah_N (\eo_{\oK}) \longrightarrow \Ah_N (\eo_n)$ we get a map  $\hAh_N (\eo_n) \longrightarrow \Hom (A_N (\oK) , \eo^*_n)$.
Together with the projection
$TA \longrightarrow A_N (\oK)$
this gives a homomorphism
$  \hAh_N (\eo_n) \longrightarrow \Hom_c (TA , \eo^*_n)$.
Note here that $\eo_n$ carries the discrete topology. For $N \tei M$ the corresponding maps are compatible with the inclusion of $\hAh_N (\eo_n)$ in $\hAh_M (\eo_n)$.

Now, the abelian group $\hAh (\eo_n)$ is torsion since it is the union of the finite groups $\hAh (\eo_L / p^n \eo_L)$ where $L$ runs over the finite extensions of $K$ in $\oK$. Hence we get a homomorphism
$  \hAh (\eo_n) \longrightarrow \Hom_c (TA , \eo^*_n)$. Composing with the reduction map $\hA (\C_p) = \hAh (\eo) \longrightarrow \hAh (\eo_n)$
we obtain a homomorphism\\
$  \alpha_n : \hA (\C_p) \longrightarrow \Hom_{c} (TA , \eo^*_n)$. For every $\ha \in \hA (\C_p)$ and any $n \le m$ the $\eo^*_n$-valued map $\alpha_n (\ha)$ is the reduction $\mod p^n$ of the $\eo^*_n$-valued map $\alpha_m (\ha)$. Thus the $\alpha_n$'s define a homomorphism
\begin{equation}
  \label{eq:10}
  \alpha : \hA (\C_p) \longrightarrow \Hom_c (TA , \eo^*) = \Hom_c (TA , \C^*_p)
\end{equation}
as desired. It is Galois-equivariant by construction. Moreover, for any homomorphism $\varphi : B \to A$ of abelian varieties with good reduction we have a commutative diagram
\[
\xymatrix{
\hA (\C_p) \ar[r]^(.35){\alpha_A} \ar[d]_{\hat{\varphi}} & \Hom_c (TA , \C^*_p) \ar[d]^{(T \varphi)^*} \\
\hat{B} (\C_p) \ar[r]^(.35){\alpha_B} & \Hom_c (TB , \C^*_p) \; .
}
\]

We now consider the restriction $\alpha_{\tors}$ of the map $\alpha$ in (\ref{eq:10}) to the torsion part of $\hA (\C_p)$:
\[
  \alpha_{\tors} : \hA (\C_p)_{\tors} = \hA (\oK)_{\tors} \longrightarrow \Hom_c (TA , \eo^*)_{\tors} = \Hom_c (TA , \mu) \; .
\]
Here $\mu = \mu (\oK)$ is the group of roots of unity in $\eo^*$ or $\oK^*$. Note that the Kummer sequence
on $\oA_{\et}$ induces an isomorphism
\[
i_A : H^1 (\oA , \mu_N) \silo H^1 (\oA , \Ge_m)_N \; .
\]

\begin{prop}
  \label{t2}
The map $\alpha_{\tors}$ is an isomorphism. On $\hA (\oK)_N$ it coincides with the composition:
\[
\hA (\oK)_N = H^1 (\oA , \Ge_m)_N \xrightarrow{i^{-1}_A \atop \sim} H^1 (\oA , \mu_N) = \Hom_c (TA , \mu_N) \; .
\]
\end{prop}

\begin{proof}
By construction, the restriction of $\alpha$ to $\hA (\oK)_N$ is the map
\[
\hA (\oK)_N \longrightarrow \Hom (A_N (\oK) , \mu_N ) = \Hom (TA , \mu_N)
\]
coming from Cartier duality $\hat{\overline{A}}_N \simeq \uHom (\overline{A}_N , \Ge_{m_{\oK}})$ over $\oK$. The canonical identification $\Hom_c (TA , \mu_N) = H^1 (\overline{A} , \mu_N)$ can be factorized by the isomorphisms
\[
\Hom (\overline{A}_N , \mu_N) \silo \Ext^1 (\overline{A} , \mu_N) \silo H^1 (\overline{A} , \mu_N) \; ,
\]
where the first map is induced by the exact sequence $0 \to \overline{A}_N \to \overline{A} \to \overline{A} \to 0$ and the second map is the forgetful map associating to an extension the corresponding $\mu_N$-torsor. Since the diagram
\[
\xymatrix{
\Hom (\overline{A_N} , \mu_N) \ar[r]^{\sim} \ar[d]_{\wr} & \Ext^1 (\overline{A} , \mu_N) \ar[r]^{\sim} \ar[d] & H^1 (\overline{A} , \mu_N) \ar[d]^{i_A}_{\wr} \\
\Hom (\overline{A}_N , \Ge_m) \ar[r]^{\sim} & \Ext^1 (\overline{A}, \Ge_m)_N \ar[r]^{\sim} \ar@{=}[d] & H^1 (\overline{A} , \Ge_m)_N \\
\hat{\overline{A}}_N (\oK) \ar[u]_{\wr} \ar@{=}[r] & \hat{\overline{A}}_N (\oK) & 
}
\]
commutes, our claim follows. 
\end{proof}

Next, we need an elementary lemma. Consider an abelian topological group $\pi$ which fits into an exact sequence of topological groups
\[
0 \longrightarrow H \longrightarrow \pi \longrightarrow \hat{\Z}^n \longrightarrow 0
\]
where $H$ is a finite discrete group. Later $\pi$ will be the abelianized fundamental group of an algebraic variety. Applying the functor $\Hom_c (\pi , \, \underline{\ \ })$ to the exact sequence
\[
0 \longrightarrow \mu  \longrightarrow \eo^* \xrightarrow{\log} \C_p \longrightarrow 0
\]
we get the sequence
\begin{equation}
  \label{eq:12}
\small  0 \longrightarrow \Hom_c (\pi , \mu ) \longrightarrow \Hom_c (\pi , \eo^*) \xrightarrow{\log_*} \Hom_c (\pi , \C_p) \longrightarrow 0 \; .
\end{equation}

\begin{lemma}
  \label{t3}
a) The sequence (\ref{eq:12}) is exact.\\
b) In the topology of uniform convergence $\Hom_c (\pi , \eo^*)$ is a complete topological group. It contains $(\eo^*)^n$ as an open subgroup of finite index and hence acquires a natural structure as a Lie group over $\C_p$. Its Lie algebra is $\Hom_c (\pi , \C_p)$ and the logarithm map is given by $\log_*$.
\end{lemma}

\begin{proof}
a)  Since $\Hom_c (\pi , \underline{\ \ })$ is left exact, it suffices to show that $\log_*$ is surjective. As $\Hom_c (H , \C_p) = 0$ we only have to show surjectivity of $\log_*$ for $\pi = \hat{\Z}^n$, hence for $\pi = \hat{\Z}$.
We first prove that the injective evaluation map $\ev : \Hom_c (\hat{\Z} , \eo^*) \longrightarrow \eo^* \; , \; \ev (\varphi) = \varphi (1)$
is surjective. \\
Set $U_1 = \{ x \in \eo^* \tei |x-1| < 1 \}$ and $U_0 = \{ x \in \eo^* \tei |x-1| < p^{-\frac{1}{p-1}} \}$. The logarithm provides an isomorphism
\[
\log : U_0 \silo V_0 = \{ x \in \C_p \tei |x| < p^{-\frac{1}{p-1}} \} \; ,
\]
whose inverse is the exponential map. Therefore $U_0$ is a $\Z_p$-module and it follows that
$\ev : \Hom_c (\Z_p , U_0) \silo U_0$
is an isomorphism. We claim that \linebreak
$\ev : \Hom_c (\Z_p , U_1) \hookrightarrow U_1$ is an isomorphism as well.
Fix some $b$ in $U_1$. We construct a continuous map $\psi : \Z_p \to U_1$ with $\psi (1) = b$ as follows. There is some $N \ge 1$ such that $b^{p^N} \in U_0$. Hence there is a continuous homomorphism $\varphi : p^N \Z_p \longrightarrow U_1$ such that $\varphi (p^N \nu) = (b^{p^N})^{\nu}$ for all $\nu \in \Z$. Because of the decomposition $\eo^* = \mu_{(p)} \times U_1$ the group $U_1$ is divisible. Hence there is a homomorphism $\psi' : \Z_p \to U_1$ whose restriction to $p^N \Z_p$ equals $\varphi$. It follows that $\psi'$ is continuous as well. Because of $\psi' (1)^{p^N} = \psi' (p^N) = b^{p^N}$ there is a root of unity $\zeta \in \mu_{p^N}$ with $\psi' (1) = \zeta b$. Take the continuous homomorphism $\psi'' : \Z_p \to \mu_{p^{\infty}} \subset U_1$ with $\psi'' (1) = \zeta^{-1}$ and set $\psi = \psi' \cdot \psi''$.

The natural projection $\hat{\Z} \to \Z_p$ induces a commutative diagram
\[
\xymatrix{
\Hom_c (\Z_p , U_1) \ar@{^{(}->}[rr] \ar[dr]^{\sim}_{\ev} & & \Hom_c (\hat{\Z} , U_1) \ar@{^{(}->}[dl]^{\ev} \\
 & U_1 & .
}
\]
It follows that $\ev : \Hom_c (\hat{\Z} , U_1) \to U_1$ is an isomorphism as well. Using the decomposition $\eo^* = \mu_{(p)} \times U_1$ where $\mu_{(p)}$ carries the discrete topology we conclude that $\ev : \Hom_c (\hat{\Z} , \eo^*) \longrightarrow \eo^*$ is an isomorphism. Using the commutative diagram
\[
\xymatrix{
\Hom_c (\hat{\Z} , \eo^*) \ar[r]^{\log_*} \ar[d]^{\wr \, \ev} & \Hom_c (\hat{\Z} , \C_p) \ar@{=}[r] & \Hom_c (\Z_p , \C_p) \ar[d]^{\wr \, \ev} \\
\eo^* \ar@{>>}[rr]^{\log} & & \C_p
}
\]
we see that $\log_*$ is surjective for $\pi = \hat{\Z}$ and hence in general.\\
b) With the topology of uniform convergence, $\Hom_c (\pi , \eo^*)$ becomes a topological group. This topology comes from the inclusion of $\Hom_c (\pi , \eo^*)$ into the $p$-adic Banach space $C^0 (\pi , \C_p)$ of continuous functions from $\pi$ to $\C_p$ with the norm 
\[
\| f \| = \max_{\gamma \in \pi} |f (\gamma)| \; .
\]
Since $\Hom_c (\pi , \eo^*)$ is closed in $C^0 (\pi , \C_p)$ it becomes a complete metric space and hence it is a complete topological group. We now observe that the continuous evaluation map $\ev : \Hom_c (\hat{\Z} , \eo^*) \silo \eo^*$ is actually a homeomorphism. Let $x_n \to x$ be a convergent sequence in $\eo^*$ and let $\varphi_{\nu} , \varphi : \hat{\Z} \to \eo^*$ be the continuous homomorphisms with $\varphi_{\nu} (1) = x_{\nu}$ and $\varphi (1) = x$. Since $\Z^{\ge 1}$ is dense in $\hat{\Z}$ we get
\begin{eqnarray*}
  \| \varphi - \varphi_{\nu} \| & = & \max_{\gamma \in \hat{\Z}} |\varphi (\gamma) - \varphi_{\nu} (\gamma)| = \sup_{n \ge 1} |\varphi (n) - \varphi_{\nu} (n) | = \sup_{n \ge 1} |x^n - x^n_{\nu}| \\
& = & |x - x_{\nu}| \sup_{n \ge 1} \, |\sum^{n-1}_{i=0} x^i x^{n-i-1}_{\nu}| \le |x - x_{\nu}| \; .
\end{eqnarray*}
Hence $\varphi_{\nu}$ converges uniformly to $\varphi$. 

It follows that $\Hom_c (\hat{\Z}^n , \eo^*)$ and $(\eo^*)^n$ are isomorphic as topological groups. The exact sequence of topological groups
\[
0 \to \Hom_c (\hat{\Z}^n , \eo^*) \to \Hom_c (\pi , \eo^*) \to \Hom_c (H , \eo^*) = \Hom_c (N , \mu_{|H|})
\]
shows that $\Hom_c (\pi , \eo^*)$ contains $(\eo^*)^n$ as an open subgroup of finite index. Hence $\Hom_c (\pi , \eo^*)$ becomes a Lie group over $\C_p$. It is clear that the analytic structure depends only on $\pi$ and not on the choice of an exact sequence \\
$0 \to H \to \pi \to \hat{\Z}^n \to 0$ as above. The remaining assertions have to be checked for $\pi = \hat{\Z}^n$ and hence for $\pi = \hat{\Z}$ only where they are clear by the preceeding observations.
\end{proof}

\begin{rem}
  The proof shows that the topologies of uniform and pointwise convergence on $\Hom_c (\pi , \eo^*)$ coincide.
\end{rem}

By \cite{Bou}, III, \S7, no.6 there is a logarithm map on an open subgroup $U$
of the $p$-adic Lie group $\hA (\C_p)$, mapping $U \rightarrow  \Lie \hA (\C_p)$, such that the kernel consists of the torsion points in $U$. Since $\hA (\C_p) / U$ is torsion (see Theorem 4.1 in \cite{Co2}), the logarithm has a unique extension to the whole Lie group $\hA (\C_p)$.  Since the $\C_p$-vector
space $\Lie \hA (\C_p)$ is uniquely divisible, it is surjective. Therefore we have the exact sequence
\begin{equation}
  \label{eq:13}
  0 \longrightarrow \hA (\C_p)_{\tors} \longrightarrow \hA (\C_p) \xrightarrow{\log} \Lie \hA (\C_p) = H^1 (A , \Oh) \otimes \C_p \longrightarrow 0 \; .
\end{equation}
Using proposition \ref{t2} and lemma \ref{t3} we therefore get a commutative diagram with exact lines and $G_K$-equivariant maps
\begin{equation}
  \label{eq:14}
\def\objectstyle{\scriptstyle}
 \vcenter{\xymatrix@-1.2pc{
0 \ar[r] & \hA (\C_p)_{\tors} \ar[r] \ar[d]_{\alpha_{\tors}}^{\wr} & \hA (\C_p) \ar[r]^{\log} \ar[d]_{\alpha} & \Lie \hA (\C_p) \ar@{=}[r] \ar[d]^{\tilde{\alpha}} & H^1 (A , \Oh) \otimes \C_p \ar[r] & 0 \\
0 \ar[r] & \Hom_c (TA , \mu (\oK)) \ar[r] & \Hom_c (TA , \eo^*) \ar[r] & \Hom_c (TA , \C_p) \ar@{=}[r] & H^1_{\et} (\oA , \Q_p) \otimes \C_p \ar[r] & 0 
}}
\end{equation}
Here $\tilde{\alpha}$ is the map induced by $\alpha$ on the quotients. 

We will see in section 3 that $\alpha$ is a $p$-adically analytic map of $p$-adic Lie groups. It follows that $\tilde{\alpha} = \Lie \alpha$. Note here that by lemma \ref{t3} we have:
\[
\Hom_c (TA , \C_p) = \Lie \Hom_c (TA, \C^*_p) \; .
\]
Faltings in \cite{Fa2} Theorem 4, b  and Coleman in \cite{Co1} p. 379 and \cite{Co3} \S\,4 have given an elegant construction of the $G_K$-equivariant map arising in the Hodge--Tate decomposition
\begin{equation}
\label{eq:15}
\theta^*_A : H^1 (A, \Oh) \otimes \C_p \longrightarrow H^1_{\et} (\oA , \Q_p) \otimes \C_p \; .
\end{equation}
Consider the universal vectorial extension of $\Ah$ over the ring $\eo_K$
\[
0 \longrightarrow \omega_{\hAh} \longrightarrow E \longrightarrow \Ah \longrightarrow 0 \; .
\]
Here $\omega_{\hAh}$ is the vector group induced by the invariant differentials on $\hAh$, i.e. 
\[
\omega_{\hAh} (S) = H^0 (S , e^* \Omega^1_{\hAh_S / S})
\]
for any $\eo_K$-scheme $S$, where $e$ denotes the zero section . For $\nu \ge 1$ consider the map
\[
\Ah_{p^{\nu}} (\eo) \longrightarrow \omega_{\hAh} (\eo) / p^{\nu} \omega_{\hAh} (\eo)
\]
obtained by sending $a_{p^{\nu}}$ to $p^{\nu} b_{p^{\nu}}$ where $b_{p^{\nu}} \in E (\eo)$ is a lift of $a_{p^{\nu}}$. Passing to inverse limits we get a $\Z_p$-linear homomorphism
\[
\theta_A : T_p A \longrightarrow \omega_{\hAh} (\eo) \; .
\]
The $\C_p$-dual of the resulting map
\[
\theta_A : T_p A \otimes \C_p \longrightarrow \omega_{\hAh}(\eo) \otimes \C_p = \omega_{\hA}(\C_p)
\]
is the map $\theta^*_A$ in (\ref{eq:15}).

In \cite{Co1}, Coleman proved that together with a map defined by Fontaine in \cite{Fo}
\[
H^0 (A, \Omega^1) \otimes \C_p (-1) \longrightarrow H^1_{\et} (\oA , \Q_p) \otimes \C_p \; ,
\]
$\theta^*_A$ gives a Hodge--Tate decomposition
\begin{equation} \label{eq:16}
H^1_{\et} (\oA , \Q_p) \otimes \C_p \cong (H^1 (A , \Oh) \otimes \C_p) \oplus (H^0 (A, \Omega^1) \otimes \C_p (-1)) \; .
\end{equation}
In particular $\theta^*_A$ is injective.

\begin{theorem}
  \label{t4}
We have $\tilde{\alpha} = \Lie \alpha = \theta^*_A$. In particular $\alpha$ is injective.
\end{theorem}

We will give a direct proof of the formula $\Lie \alpha = \theta^*_A$ in section 4. An indirect proof is obtained as follows. On the open subgroup $\hAh (p) (\eo) \subset A (\C_p)$ Tate has shown in \cite{Ta}, \S\,4 that the diagram
\[
\xymatrix{
\hAh (p) (\eo) \ar[d]_{\Phi} \ar[r]^{\log} & \Lie \hAh (p) \ar[d]^{\Lie \Phi} \ar@{=}[r] & H^1 (A , \Oh) \otimes \C_p \\
\Hom_c (T (\Ah (p)) , U_1) \ar[r]^{\log_*} & \Hom_c (T (\Ah (p)), \C_p) \ar@{=}[r] & H^1_{\et} (\oA , \Q_p) \otimes \C_p
}
\]
commutes where $\Phi$ can be identified with the restriction of $\alpha$ and $\Lie \Phi$ is its tangent map. Combining Coleman's work in \cite{Co1} and \cite{Co3} \S\,4 with Fontaine's results, specifically \cite{Fo} Proposition 11 it follows that $\Lie \Phi = \theta^*_A$.

Applying the snake lemma to diagram (\ref{eq:14}) we see that $\alpha$ is injective because $\tilde{\alpha} = \theta^*_A$ is injective. \hspace*{\fill} $\Box$

\begin{rem}
  With the higher rank case in mind, it would be interesting to find a proof that $\alpha$ is injective which does not depend on the Hodge--Tate decomposition. 
\end{rem}

The next corollary was already observed by Tate in his context of $p$-divisible groups, \cite{Ta}, Theorem 3.

\begin{cor}
  \label{t5}
The map $\alpha$ induces an isomorphism of abelian groups
\[
\alpha : \hA (K) \silo \Hom_{c, G_K} (TA, \eo^*) \; .
\]
\end{cor}

\begin{proof}
  According to \cite{Ta}, Theorem 1 we have $H^0 (G_K , \C_p) = K$ and \\
$H^0 (G_K , \C_p (-1)) = 0$. Hence the Hodge--Tate decomposition and theorem \ref{t4} imply that $\tilde{\alpha}$ induces an isomorphism:
\[
\tilde{\alpha} : H^1 (A , \Oh) \silo H^0 (G_K , H^1_{\et} (\oA , \Q_p) \otimes \C_p) = \Hom_{c, G_K} (TA , \C_p) \; .
\]
We have $H^0 (G_K , \hA (\C_p)) = \hA (K)$. This follows for example by embedding $\hA$ into some $\Pa^N$ over $K$ and using the corresponding result for $\Pa^N$. The latter is a consequence of the decomposition $\Pa^N = \A^N \amalg \ldots \amalg \A^0$ over $K$ and the equalilty $H^0(G_K , \C_p) = K$.
The corollary follows by applying the 5-lemma to the diagram of Galois cohomology sequences obtained from (\ref{eq:14}).
\end{proof}

We next describe the image of the map $\alpha$ from (\ref{eq:10}).

\begin{defn} \label{t6}
  A continuous character $\chi : TA \to \C^*_p$ is called smooth if its stabilizer in $G_K$ is open. The group of smooth characters of $TA$ is denoted by $Ch^{\infty} (TA)$.
\end{defn}

Note that we have
\[
Ch^{\infty} (TA) = \varinjlim_{L/K} \Hom_{c, G_L} (TA , \eo^*)
\]
where $L$ runs over the finite extensions of $K$. It is also the biggest $G_K$-invariant subset $S$ of $\Hom_c (TA , \C^*_p)$ such that the $G_K$-action on $S^{\delta}$ is continuous. Here $S^{\delta}$ is $S$ endowed with the discrete topology.

Replacing $A$ by $A_L$ in Corollary \ref{t5} we find that $\alpha$ induces an isomorphism
\[
\alpha : \hA (\oK) \silo Ch^{\infty} (TA) \subset \Hom_c (TA , \eo^*) \; .
\]

Let $Ch (TA)$ be the closure of $Ch^{\infty} (TA)$ in $\Hom_c (TA , \eo^*)$ or equivalently in $C^0 (TA , \C_p)$. Then $Ch (TA)$ is also a complete topological group.

\begin{theorem} \label{t7}
  The map $\alpha$ induces an isomorphism of topological groups
\[
\alpha : \hA (\C_p) \silo Ch (TA) \; .
\]
\end{theorem}

\begin{proof}
  In Lemma \ref{t8} of the next section we will show that the map
\[
\beta : \hA (\C_p) \times TA \longrightarrow \eo^* \; , \; (\ha , \gamma) \longmapsto \alpha (\ha) (\gamma)
\]
is continuous. Since $TA$ is compact, continuity of $\alpha$ follows. We conclude that
\[
\alpha (\hA (\C_p)) = \alpha (\overline{\hA (\oK)}) \subset \overline{\alpha (\hA (\oK))} = Ch (TA) \; .
\]
It now suffices to show that $\alpha$ in (\ref{eq:10}) is a closed map. Namely, because of
\[
Ch^{\infty} (TA) \subset \alpha (\hA (\C_p)) \subset Ch (TA)
\]
it will follow that $\alpha (\hA (\C_p)) = Ch (TA)$ and $\alpha$ will be a homeomorphism onto its image.

So let $Y \subset \hA (\C_p)$ be a closed set. Let $\alpha (y_n)$ for $y_n \in Y$ be a sequence which converges to some $\chi$ in $\Hom_c (TA , \eo^*)$. Since the map $\log_*$ in (\ref{eq:12}) is continuous in the uniform topologies it follows that $\tilde{\alpha} (\log y_n) = \log_* \alpha (y_n)$ converges to $\log_* \chi$. Because of the equality 
\[
\Hom_c (TA , \C_p) = \Hom_{\Z_p} (T_p A , \C_p)
\]
the topology of uniform convergence on this space coincides with its topology as a finite dimensional $\C_p$-vector space. The map $\tilde{\alpha} = \Lie \alpha$ is a $\C_p$-linear injection by the Hodge--Tate decomposition and theorem \ref{t4}. Hence it is a closed injection and therefore the sequence $\log y_n$ converges. As $\log$ is a local homeomorphism there is a convergent sequence $\tilde{y}_n \in \hA (\C_p)$ with $\log \tilde{y}_n = \log y_n$. Writing $\tilde{y}_n = y_n + t_n$ with $t_n \in \hA (\C_p)_{\tors}$ we get $\alpha (\tilde{y}_n) = \alpha (y_n) + \alpha (t_n)$. The sequence $\alpha (y_n)$ converges by assumption and the sequence $\alpha (\tilde{y}_n)$ converges because $\alpha$ is continuous. Hence the sequence $\alpha (t_n) = \alpha_{\tors} (t_n)$ converges. The groups $\hA (\C_p)_{\tors}$ and $\Hom_c (TA , \mu )$ are the kernels of the locally topological homomorphisms $\log$ resp. $\log_*$. Hence they inherit the discrete topology from the $p$-adic topologies on $\hA (\C_p)$ resp. $\Hom_c (TA, \eo^*)$. Therefore the algebraic isomorphism $\alpha_{\tors}$ is trivially a homeomorphism and hence the sequence $t_n$ converges.  It follows that the sequence $y_n$ converges to some $y \in Y$. By continuity of $\alpha$ the sequence $\alpha (y_n)$ converges to $\alpha (y)$. Thus $\alpha (Y)$ is closed as was to be shown.
\end{proof}

The following example was prompted by a question of Damian Roessler.

{\bf Example} Fix some $\sigma$ in $G_K$. Since $\alpha$ is $G_K$-equivariant we know that if $\ha \in \hA (\C_p)$ corresponds to the character $\chi : TA \to \eo^*$ then $\sigma (\ha)$ corresponds to $^{\sigma} \chi = \sigma \verk \chi \verk \sigma^{-1}_*$. Here $\sigma_*$ is the action on $TA$ induced by $\sigma$.

How about the character $\sigma \verk \chi : TA \to \eo^*$? Using theorem \ref{t7} we will now show that it also corresponds to an element of $\hA (\C_p)$ provided that $A$ has complex multiplication over $K$. For this, we have to check that in the $CM$ case, the subgroup $Ch (TA)$ is invariant under the homeomorphism $\chi \mapsto \sigma \verk \chi$ of $\Hom_c (TA , \eo^*)$. It suffices to show that $Ch^{\infty} (TA)$ is invariant. For $\chi$ in $Ch^{\infty} (TA)$, there is a finite normal extension $N /K $ such that $\chi$ is $G_N$-invariant, i.e. $\tau^{-1} \chi \tau_* = \chi$ for all $\tau$ in $G_N$. It follows that
\[
\tau^{-1} (\sigma \chi) \tau_* = \sigma (\sigma^{-1} \tau^{-1} \sigma \chi) \tau_* = \sigma \chi (\sigma^{-1} \tau^{-1} \sigma)_* \tau_* = \sigma \chi [\sigma , \tau]_*
\]
where we define the commutator by $[\sigma , \tau] = \sigma^{-1} \tau^{-1} \sigma \tau$. By the $CM$ assumption, the image of $G_K$ in the automorphism group of $TA$ is abelian. Hence $[\sigma , \tau]_*$ acts trivially on $TA$ and we have thus shown that $\tau^{-1} (\sigma \chi) \tau_* = \sigma \chi$ for all $\tau \in G_N$. Hence $\sigma \verk \chi$ lies in $Ch^{\infty} (TA)$. This proves the claim. 

For $\ha \in \hA (\C_p)$ let $\ha_{\sigma} \in \hA (\C_p)$ be the element corresponding to $\sigma \verk \chi$ via theorem \ref{t7}. By construction, the map $(\sigma , \ha) \mapsto \ha_{\sigma}$ determines a new action of $G_K$ on $\hA (\C_p)$. It seems to be a nice exercise in $CM$-theory to give an explicit description of this action. 
%\newpage
%\input{sec4}
\section{The map $\alpha$ as a Lie group homomorphism}

The continuous homomorphism
\[
\alpha : \hA (\C_p) \longrightarrow \Hom_c (TA , \eo^*)
\]
defined in section 2 induces a pairing
\[
\beta : \hA (\C_p) \times TA \longrightarrow \eo^*
\]

\begin{lemma}
  \label{t8}
$\beta$ is continuous.
\end{lemma}

\begin{proof}
  Denote by $r_n : \hAh (\eo) \to \hAh (\eo_n)$ the reduction map. Since the kernel of $r_n$ is $p$-adically open, it contains an open neighbourhood $W \subseteq \hA (\C_p)$ of zero.

Fix $(\ha, \gamma) \in \hA (\C_p) \times TA$ mapping to $z = \beta (\ha , \gamma)$. We show that the preimage of the open neighbourhood $z (1 + p^n \eo)$ is open. Let $N \ge 1$ be big enough so that $r_n (\ha)$ is contained in $\hAh_N (\eo_n)$. If $U$ denotes the kernel of the projection $TA \to A_N (\oK)$, the neighbourhood $(\ha + W, \gamma + U )$ of $(\ha, \gamma )$ maps to $z (1 + p^n \eo)$ under $\beta$.
\end{proof}

Let us fix some $\gamma \in TA$, and denote by $\psi_{\gamma}$ the induced homomorphism
\[
\psi_{\gamma} = \beta (-,\gamma) : \hA (\C_p) \longrightarrow \eo^* \; .
\]

\begin{prop}
  \label{t9}
$\psi_{\gamma}$ is an analytic map, hence a Lie group homomorphism.
\end{prop}

\begin{proof}
  We will briefly write $\psi = \psi_{\gamma}$ in this proof. It suffices to show that $\psi$ is analytic in a neighbourhood of the zero element $e_{\C_p} \in \hA (\C_p)$.

Let $e \in \hAh_{\eo} (\eo)$ be the zero section of $\hAh_{\eo}$. Since $\hAh_{\eo}$ is smooth over $\eo$, there is a Zariski open neighbourhood $U \subseteq \hAh_{\eo}$ of $e$ of the form
\[
U = \Spec \eo [x_1 , \ldots , x_{m+r}] / (f_1 , \ldots , f_m)
\]
such that the matrix $\left( \frac{\partial f_i}{\partial x_{r+j}} (e) \right)_{i,j = 1 \ldots m}$ is invertible in $M_m(\eo)$. 

By the theorem of implicit functions (see e.g. \cite{Col}, A.3.4), $U (\eo)$ contains an open neighbourhood $V$ of $e$ in the $p$-adic topology, such that the projection map $q : U (\eo) \subseteq \eo^{m+r} \to \eo^r$ given by $(x_1 , \ldots , x_{m+r}) \longmapsto (x_1 , \ldots , x_r)$ maps $e$ to $(0, \ldots , 0)$ and induces a homeomorphism
\[
q : V \longrightarrow V_1
\]
between $V$ and an open ball $V_1 \subseteq \eo^r$ around zero. This is an analytic chart around $e_{\C_p}$. A function $f$ on $V$ is locally analytic around $e_{\C_p}$, iff it induces a function on $V_1$ which coincides on a ball $V'_1 \subseteq V_1$ around $0$ with a power series in $x_1, \ldots , x_r$ converging pointwise on $V'_1$.

Since by \cite{Bou}, III, \S\,7, Prop. 10 and 11 the logarithm map on $\hA (\C_p)$ is locally around $e_{\C_p}$ an analytic isomorphism respecting the group structures, there exists an open subgroup $H$ of $\hA (\C_p)$ such that $(p^{\nu} H)_{\nu \ge 0}$ is a basis of open neighbourhoods of $e_{\C_p}$. By shrinking $V$ if necessary, we can assume that $V \subseteq H$.

For $n \ge 1$ we denote by $r_n$ as before the reduction map
\[
r_n : \hA (\C_p) = \hAh (\eo) \longrightarrow \hAh (\eo_n) = \hAh_n (\eo_n) \; ,
\]
where $\hAh_n = \hAh \otimes_{\eo} \eo_n$.

Since the kernel of $r_n$ is an open subgroup of $\hA (\C_p)$, it contains $p^{\nu_n} H$ for a suitable $\nu_n \ge 0$. Hence $p^{\nu_n} V$ is contained in the kernel of $r_n$, which implies that for all $x \in V$ the point $r_n (x)$ lies in the scheme $\hAh_{n , p^{\nu_n}}$ of $p^{\nu_n}$-torsion points in $\hAh_n$. 

The element $\gamma$ in $TA$ induces a point in $\Ah_{n, p^{\nu_n}} (\eo_n)$, whose image under the Cartier duality morphism
\[
\Ah_{n, p^{\nu_n}} (\eo_n) \longrightarrow \Hom (\hAh_{n , p^{\nu_n}} , \Ge_{m, \eo_n})
\]
we denote by $\psi_n$. Then $\psi (x)$ for $x \in V$ is by definition the element in $\eo^* \subseteq \eo$ satisfying $\psi (x) \equiv \psi_n (r_n (x)) \mod p^n$ for all $n$.

Let $U_n = U \otimes_{\eo} \eo_n = \Spec \eo_n [x_1 , \ldots , x_{m+r}] / (\overline{f}_1 , \ldots , \overline{f}_m)$ be the reduction of the affine subscheme $U \subseteq \hAh_{\eo}$.
We write $\overline{f}$ for the reduction of a polynomial $f$ over $\eo$ modulo $p^n$.

Then $U_n \cap \hAh_{n , p^{\nu_n}} = \Spec \eo_n [x_1 , \ldots , x_{m+r}] /  \ea$ for some ideal $\ea$ containing $(\of_1 , \ldots , \of_m)$. Since $\psi_n$ is an algebraic morphism, it is given on $U_n \cap \hAh_{n , p^{\nu_n}}$ by a unit in $\eo_n [x_1 , \ldots , x_{m+r}] / \ea$, which is induced by a polynomial $\overline{g}_n \in \eo_n [x_1 , \ldots , x_{m+r}]$. Let $g_n \in \eo [x_1 , \ldots , x_{m+r}]$ be a lift of $\overline{g}_n$. 

The implicit function theorem also implies that possibly after shrinking $V$ and $V_1$, we find power series $\theta_1 , \ldots , \theta_m \in \C_p [[ x_1 , \ldots , x_r]]$ converging in all points of $V_1 \subseteq \eo^r$, such that the map $V_1 \xrightarrow{q^{-1}} V \subseteq U (\eo) \subseteq \eo^{m+r}$ is given by
\[
(x_1 , \ldots , x_r) \longmapsto (x_1 , \ldots , x_r , \theta_1 (x_1 , \ldots , x_r) , \ldots , \theta_m (x_1, \ldots, x_r)) \; .
\]

For all $i = 1 , \ldots , m$ and all $n \ge 1$ let $h_{i,n} \in \C_p [x_1 , \ldots , x_r]$ be a polynomial satisfying $\theta_i (x) - h_{i,n} (x) \in p^n \eo$ for all $x \in V_1$. We can obtain $h_{i,n}$ by truncating $\theta_i$ suitably.

Then the map $
V_1 \silo V \xrightarrow{r_n} U_n (\eo_n) \cap \hAh_{n,p^{\nu_n}} (\eo_n) \xrightarrow{\psi_n} \eo^*_n$ maps the point
$
(x_1 , \ldots , x_r)$ to $ \overline{g}_n (\ox_1 , \ldots , \ox_r , \overline{\theta_1 (x_1 ,\ldots, x_r)} , \ldots , \overline{\theta_m (x_1, \ldots, x_r)}) $.

Hence for all $x = (x_1 , \ldots , x_r) \in V_1$ we have 
\[
\psi (q^{-1} (x)) - g_n (x_1 , \ldots , x_r , h_{1,n} (x_1, \ldots, x_r), \ldots , h_{m,n} (x_1, \ldots, x_r)) \in p^n \eo \; .
\]

Thus $\psi$ is the uniform limit of polynomials with respect to the coordinate chart $V_1$. This implies our claim.
\end{proof}

\begin{cor}
  \label{t10}
The homomorphism $\alpha : \hA (\C_p) \longrightarrow \Hom_c (TA , \eo^*)$ is analytic.
\end{cor}

\begin{proof}
  Recall from section 2 that evaluation in $1 \in \hat{\Z}$ induces an isomorphism $\Hom_c (\hat{\Z} , \eo^*) \simeq \eo^*$. Hence $\Hom_c (TA , \eo^*) \simeq (\eo^*)^{2g}$ by evaluation in a $\hat{\Z}$-basis
 $\gamma_1 , \ldots , \gamma_{2g}$ of $TA$. This isomorphism induces the analytic structure on $\Hom_c (TA, \eo^*)$. Hence our claim follows from Proposition \ref{t9}.
\end{proof}

Now we can calculate the Lie algebra map induced by $\psi_{\gamma}$. 
For every commutative group scheme $G$ over the base $S$ we denote by 
$\omega_G$ the quasi-coherent sheaf $e^\ast \Omega^1_{G/S}$ on $S$, where $e$ is the unit section 
of $G$. As usual, we regard $\omega_G$ as a sheaf on the fppf site over $S$ by the rule
$\omega_G(S') = \Gamma (S', e^\ast \Omega^1_{G_{S'}/S'})$ for $S'$ over $S$. 

We always identify $\Lie \eo^*$ with $\eo$ via the invariant differential $\frac{dT}{T}$
on $\Ge_{m, \eo} = \mbox{Spec } \eo[T, T^{-1}]$. 
\begin{prop}
  \label{t11}
For every $\gamma \in T_pA$ the Lie algebra homomorphism
\[
\Lie \psi_{\gamma} : \Lie \hA (\C_p) \longrightarrow \Lie \eo^* \silo \eo
\]
is given by the invariant differential $\theta_A (\gamma)$, where $\theta_A : T_p A \to \omega_{\hAh} (\eo)$ is the map defined by Faltings and Coleman (see section 2).
\end{prop}

\begin{proof}
By definition, $\Lie \psi_{\gamma}$ is given by the invariant differential $\psi^*_{\gamma} \frac{dT}{T} \in \omega_{\hA} (\C_p)$.
 Hence we have to show that $\psi^*_{\gamma} \frac{dT}{T} = \theta_A (\gamma)$.
We use the notation from the proof of proposition 8. Fix some $n \geq 1$ and recall the morphism
\[
\psi_n : U_n \cap \hAh_{n , p^{\nu_n}} = \Spec \eo_n [x_1 , \ldots , x_{m+r}] /  \ea \longrightarrow \Ge_{m , \eo_n} \; ,
\]
induced by a polynomial $\overline{g}_n \in \eo_n [x_1 , \ldots , x_{m+r}]$.  This $\overline{g}_n$ defines a morphism
\[
\varphi^{(n)}_n : U_n \longrightarrow \Ge_{a,\eo_n} \; .
\]
Denote by $i : \Ge_{m, \eo} \rightarrow \Ge_{a, \eo}$ the obvious closed immersion, and by $i_n$ the induced morphism over $\eo_n$. The diagram 
\[
%  \label{eq:27}
  \vcenter{\xymatrix{
U_n \ar[r]^{\varphi^{(n)}_n} & \Ge_{a , \eo_n} \\
U_n \cap \hAh_{n, p^{\nu_n}} \ar[r]^-{\psi_n} \ar@{^{(}->}[u]^{\kappa_n} & \Ge_{m, \eo_n} \ar@{^{(}->}[u]_{i_n}
}}
\]
is commutative, where $\kappa_n$ is  the canonical closed immersion. 

Similarly,  the lift $g_n$ of $\overline{g}_n$ to $\eo[x_1,\ldots, x_{m+r}]$ defines a morphism $\varphi^{(n)} : U \to \Ge_{a, \eo}$ over $\eo$, which induces an analytic map $\varphi^{(n)} : U (\eo) \longrightarrow \eo$.

The space of invariant differentials $\omega_{\hAh} (\eo)$ is an $\eo$-lattice in the $\C_p$-vector space $\omega_{\hA} (\C_p)$. For any analytic map $h : V \to \C_p$ we denote by $(dh)_e$ the  element in the cotangential space $\omega_{\hA} (\C_p) \simeq M_{e_{\C_p}} / M^2_{e_{\C_p}}$ given by the class of $h- h(e)$ modulo $M^2_{e_{\C_p}}$. Here $M_{e_{\C_p}}$ is the ideal of germs of analytic functions vanishing at $e_{\C_p}$. 

If $j : U \hookrightarrow \A^{m+r}_{\eo}$ denotes the obvious closed immersion in affine space, the fact that $\left( \frac{\partial f_i}{ \partial x_{j+r}} (e) \right)_{i,j = 1 \ldots m}$ is invertible implies that
\[
\omega_{\hAh} (\eo) = \Gamma (\Spec \eo , e^* \Omega^1_{U / \eo})
\] 
is freely generated over $\eo$ by the differentials $j^* (dx_1)_e , \ldots , j^* (dx_r)_e$. 

For all $x \in V$ we have $i \verk \psi (x) \equiv i_n \verk \psi_n (r_n( x)) \mod p^n$, hence $i \verk \psi (x) - \varphi^{(n)} (x) \in p^n \eo \subseteq \eo$. Since  $\varphi^{(n)}$ is a polynomial map, we may assume by shrinking $V$, if necessary, that for {\em all} $n$ the function 
$(i \verk \psi - \varphi^{(n)}) \verk q^{-1}$ is given by a pointwise converging power series $\sum_{I = (i_1 , \ldots , i_r)} a^{(n)}_I x^{i_1}_1 \ldots x^{i_r}_r$ on the chart $V_1 \stackrel{q^{-1}}{\longrightarrow}V$, where $V_1 = (p^t \eo)^r$ for some $t \ge 0$.  For every multi-index $I$ this implies $p^{t (i_1 + \ldots + i_r)-n} a_I^{(n)} \in \eo$. 
Then 
\[d (i \verk \psi - \varphi^{(n)})_e = a^{(n)}_{(1 , 0 , \ldots , 0)} (j^* d{x_1})_e + \ldots + a^{(n)}_{(0 , \ldots , 0,1)} (j^* dx_r)_e \in p^{n-t} \omega_{\hAh} (\eo).\]

In particular for $n \geq t$ this implies $d (i \verk \psi)_e \in \omega_{\hAh} (\eo)$. Under the isomorphism \\
$\omega_{\hAh} (\eo) / p^n \omega_{\hAh} (\eo) \to \omega_{\hAh} (\eo_n)$ the element $(d \varphi^{(n)})_e \in \omega_{\hAh} (\eo)$ maps to $(d \varphi^{(n)}_n)_e \in \omega_{\hAh} (\eo_n)$.  Moreover the diagram above implies that
\[
\kappa^*_n (d \varphi^{(n)}_n)_e = (d (i_n \verk \psi_n))_e = \psi^*_n \left( \frac{dT}{T} \right)_e \; .
\]
Besides, we can assume that $\nu_n \ge n$, so that $\eo_n$ is annihilated by $p^{\nu_n}$. Then the exact sequence 
\[\omega_{\hAh_n} \stackrel{(p^{\nu_n})^*}{\longrightarrow} \omega_{\hAh_n} {\longrightarrow}  \omega_{\hAh_{n,p^{\nu_n}}} \longrightarrow 0\]
induces an isomorphism
$\kappa^*_n : \omega_{\hAh_n}  \stackrel{\sim}{\longrightarrow} \omega_{\hAh_{n,p^{\nu_n}}}$.
This implies 
\begin{eqnarray*}
  \psi^* \left( \frac{dT}{T} \right)_e = d (i \verk \psi)_e & \equiv & (d \varphi^{(n)})_e \mod p^{n-t} \omega_{\hAh} (\eo) \\
& \equiv & (\kappa^*_n)^{-1} \psi^*_n \left( \frac{dT}{T} \right)_e \mod p^{n-t} \omega_{\hAh} (\eo) \; .
\end{eqnarray*}
Now we take a closer look at the map $\theta_A$.

By definition, $\theta_A (\gamma) \equiv p^{\nu_n} b_{p^{\nu_n}} \mod p^n \omega_{\hAh} (\eo)$, where $b_{p^{\nu_n}} \in E (\eo)$ is an arbitrary lift of the $p^{\nu_n}$-torsion point $a_{p^{\nu_n}}$ in $\Ah (\eo)$ induced by $\gamma \in T_pA$ to the universal vectorial extension $E$.

Cartier duality $[\, , \,] : \Ah_{n , p^{\nu_n}} \times \hAh_{n, p^{\nu_n}} \to \Ge_{m,\eo_n}$ induces a homomorphism
\[
\tau_n : \Ah_{n,p^{\nu_n}} \longrightarrow \omega_{\hAh_{n,p^{\nu_n}}} \; ,
\]
given by $a \mapsto [a , -]^* \frac{dT}{T}$. Now we use an argument of Crew (see \cite{Cr}, section 1 and also \cite{Ch}, Lemma A.3) to show that $\theta_A (\gamma) \equiv (\kappa^*_n)^{-1} \tau_n (\overline{a_{p^{\nu_n}}}) \mod p^n \omega_{\hAh} (\eo)$, where $\overline{a_{p^{\nu_n}}} \in \Ah_n (\eo_n)$ is the reduction of $a_{p^{\nu_n}}$. 

Namely, by \cite{Ma-Me}, Chapter I, (2.6.2), the universal vectorial extension $E_n = E \otimes_{\eo} \eo_n$ of $\Ah_n$ is isomorphic to the pushout of the sequence 
\[0 \to \Ah_{n, p^{\nu_n}} \to \Ah_n \to \Ah_n \to 0\]
 by $(\kappa_n^*)^{-1} \verk \tau_n$. Hence we have a commutative diagram with exact lines
\[
\xymatrix{
0 \ar[r] & \Ah_{n , p^{\nu_n}} \ar[r] \ar[d]_{\tau_n} & \Ah_n \ar[dd]^f \ar[r]^{p^{\nu_n}} & \Ah_n \ar[r] \ar@{=}[dd] & 0 \\
 & \omega_{\hAh_{p^{\nu_n}}} & & & \\
0 \ar[r] & \omega_{\hAh_n} \ar[u]^{\kappa^*_n}_{\wr} \ar[r] & E_n \ar[r] & \Ah_n \ar[r] & 0 \; .
}
\]
Let $\overline{b_{p^{\nu_n}}}$ be the image of $b_{p^{\nu_n}}$ under the reduction map $E (\eo) \to E_n (\eo_n)$. Since multiplication by $p^{\nu_n}$ on $\Ah (\eo)$ is surjective, we can find some $\overline{c} \in \Ah_n (\eo_n)$ with $p^{\nu_n} \overline{c} = \overline{a_{p^{\nu_n}}}$. Then $f (\overline{c})$ differs from $\overline{b_{p^{\nu_n}}}$ by an element in $\omega_{\hAh_n} (\eo_n)$, which implies
  \begin{eqnarray*}
    p^{\nu_n} \overline{b_{p^{\nu_n}}} & = & f (p^{\nu_n} \overline{c}) \\
& = & (\kappa^*_n)^{-1} \tau_n (\overline{a_{p^{\nu_n}}}) \; ,
  \end{eqnarray*}
so that indeed 
\[
\theta_A (\gamma) \equiv (\kappa^*_n)^{-1} \tau_n (\overline{a_{p^{\nu_n}}}) \mod p^n \omega_{\hAh} (\eo) \; .
\]
By definition, $\tau_n (\overline{a_{p^{\nu_n}}}) = \psi^*_n \left( \frac{dT}{T} \right)_e$, so that for all $n$
{\begin{eqnarray*}
 \psi^* \left( \frac{dT}{T} \right)_e  & \equiv & (\kappa^*_n)^{-1} \psi^*_n \left( \frac{dT}{T} \right)_e \mod p^{n-t} \omega_{\hAh} (\eo) \\
& \equiv & \theta_A (\gamma) \mod p^{n-t} \omega_{\hAh} (\eo) \; ,
\end{eqnarray*}
which implies our claim.}
\end{proof}

%\newpage
%\input{sec3}
\section{Line bundles on varieties and their $p$-adic characters}

Consider a smooth and proper variety $X$ over a finite extension $K$ of $\Q_p$. As usual, varieties are supposed to be geometrically irreducible. We assume that $H^1 (X)$ has good reduction in the sense that the inertia group $I_K$ of $G_K$ acts trivially on the \'etale cohomology group $H^1 (\oX , \Q_l)$ for some prime $l \neq p$. Here $\oX = X \otimes \oK$.

The abelianization of the fundamental group $\pi_1 (\oX , \ox)$ is independent of the choice of a base point $\ox$ and will be denoted by $\pi^{\abb}_1 (\oX)$. It carries an action of the Galois group $G_K$ even if $X$ does not have a $K$-rational point.

Using the theory of section 2 we will now attach a $p$-adic character of $\pi^{\abb}_1 (\oX)$ to any line bundle $\Lh$ on $X_{\C_p}$ whose image in the N\'eron--Severi group of $X_{\C_p}$ is torsion. 

For curves and abelian varieties a direct geometric construction of these characters will be given in section 6.

 It is known, \cite{BLR} 8.2/3 and 8.4/3 that $B := \Pic^0_{X / K}$ is an abelian variety over $K$. Its dual is the Albanese variety $A = \Abb_{X / K}$ of $X$ over $K$.

Using the Kummer sequence and divisibility  of $\Pic^0_{X/K} (\oK)$ one gets an exact sequence
\begin{equation}
  \label{eq:17}
  0 \longrightarrow B (\oK)_N \longrightarrow H^1 (\oX , \mu_N) \longrightarrow NS (\oX)_N \longrightarrow 0
\end{equation}
for every $N \ge 1$. Since the N\'eron--Severi group of $\oX$ is finitely generated it follows that $T_l B = H^1 (\oX , \Z_l (1))$ for every $l$. Thus $B$ and hence also $A$ have good reduction. For sufficiently large $N$ in the sense of divisibility we have
\[
NS (\oX)_{\tors} = NS (\oX)_N \; .
\]
Applying $\Hom (\, \underline{\ \ } \, , \mu_N)$ to the exact sequence (\ref{eq:17}) and passing to projective limits therefore gives an exact sequence of $G_K$-modules:
\begin{equation}
  \label{eq:18}
  0 \longrightarrow \Hom (NS (\oX)_{\tors} , \mu) \longrightarrow \pi^{\abb}_1 (\oX) \longrightarrow TA \longrightarrow 0 \; .
\end{equation}
Here we have used the perfect Galois equivariant pairing coming from Cartier duality
\[
A (\oK)_N \times B (\oK)_N \longrightarrow \mu_N \; .
\]

For every prime number $l$ the pro-$l$ part of the sequence (\ref{eq:18}) splits continuously since $T_l A$ is a free $\Z_l$-module. Hence (\ref{eq:18}) splits continuously and applying $\Hom_c (\underline{\ \ } , \eo^*)$ we get an exact sequence of $G_K$-modules
\begin{equation}
  \label{eq:19}
  0 \longrightarrow \Hom_c (TA , \eo^*) \longrightarrow \Hom_c (\pi^{\abb}_1 (\oX) , \eo^*) \longrightarrow NS (\oX)_{\tors} \longrightarrow 0 \; .
\end{equation}
We set
\[
\Hom^0_c (\pi^{\abb}_1 (\oX) , \eo^*) = \Ker (\Hom_c (\pi^{\abb}_1 (\oX) , \eo^*) \longrightarrow NS (\oX)_{\tors}) \; .
\]
The theory from section 2 gives a continuous injective homomorphism
\begin{equation}
  \label{eq:20}
  \alpha : \Pic^0_{X/K} (\C_p) = \hA (\C_p) \hookrightarrow \Hom_c (TA , \eo^*) = \Hom^0_c (\pi^{\abb}_1 (\oX) , \eo^*) \; .
\end{equation}
Moreover $\alpha$ is a locally analytic homomorphism of $p$-adic Lie groups over $\C_p$.  Using theorem \ref{t4}, we see that $\alpha$ fits into a commutative diagram with exact lines:
\begin{equation}
  \label{eq:21}
\vcenter{ \def\objectstyle{\scriptstyle}
 \xymatrix@C=10pt{
0 \ar[r] & \Pic^0_{X/K} (\C_p)_{\tors} \ar[d]^{\wr}_{\alpha_{\tors}} \ar[r] & \Pic^0_{X/K} (\C_p) \ar@{^{(}->}[d]_{\alpha} \ar[r]^{\log} & \Lie \Pic^0_{X/K} (\C_p) \ar@{^{(}->}[d]^{\Lie \alpha} \ar@{=}[r] & H^1 (X, \Oh) \otimes \C_p \ar[r] & 0 \\
0 \ar[r] & \Hom^0_c (\pi^{\abb}_1 (\oX) , \mu ) \ar[r] & \Hom^0_c (\pi^{\abb}_1 (\oX) , \eo^*) \ar[r]^{\log_*} & \Hom_c (\pi^{\abb}_1 (\oX) , \C_p) \ar@{=}[r] & H^1 (\oX , \Q_p) \otimes \C_p \ar[r] & 0 \; .
}}
\end{equation}

Here we have set
\begin{eqnarray*}
  \Hom^0_c (\pi^{\abb}_1 (\oX) , \mu) & = & \Hom^0_c (\pi^{\abb}_1 (\oX) , \eo^*) \cap \Hom_c (\pi^{\abb}_1 (\oX) , \mu ) \\
& = & \Hom^0_c (\pi^{\abb}_1 (\oX) , \eo^*)_{\tors} \; .
\end{eqnarray*}
Furthermore, note that:
\[
\Lie \Hom^0_c (\pi^{\abb}_1 (\oX) , \eo^*) = \Hom_c (\pi^{\abb}_1 (\oX) , \C_p) \; .
\]
The map $\Lie \alpha$ coincides with the inclusion map coming from the Hodge--Tate decomposition of $H^1 (\oX , \Q_p) \otimes \C_p$. This follows from theorem \ref{t4} and the functoriality of this decomposition.

Set
\[
Ch^{\infty} (\pi^{\abb}_1 (\oX))^0 = \varinjlim_{L/K} \Hom^0_{c , G_L} (\pi^{\abb}_1 (\oX) , \eo^*)
\]
and let $Ch (\pi^{\abb}_1 (\oX))^0$ be its closure in $\Hom^0_c (\pi^{\abb}_1 (\oX) , \eo^*)$. We make similar definitions with the $^0$'s omitted.

It follows from theorem \ref{t7} that $\alpha$ induces a topological isomorphism of complete topological groups:
\begin{equation}
  \label{eq:22}
  \alpha : \Pic^0_{X/K} (\C_p) \silo Ch (\pi^{\abb}_1 (\oX))^0 \; .
\end{equation}
We will now extend the domain of definition of $\alpha$ to $\Pic^{\tau}_{X/K} (\C_p)$. This is the group of line bundles on $X_{\C_p}$ whose image in $NS (X_{\C_p}) = NS (\oX)$ is torsion. For the latter equality, note that $NS$ is the group of connected components of $\Pic$ and that $\Pic$ and $\Pic^0$ commute with base change. We thus have an exact sequence 
\begin{equation}
  \label{eq:23}
  0 \longrightarrow \Pic^0_{X/K} (\C_p) \longrightarrow \Pic^{\tau}_{X/K} (\C_p) \longrightarrow NS (\oX)_{\tors} \longrightarrow 0 \; .
\end{equation}

\begin{theorem}
  \label{t12}
There is a $G_K$-equivariant map $\alpha^{\tau}$ which makes the following diagrams with exact lines commute:
\begin{equation}
  \label{eq:24}
  \vcenter{\xymatrix@C=10pt{
0 \ar[r] & \Pic^0_{X/K} (\C_p) \ar[d]^{\alpha} \ar[r] & \Pic^{\tau}_{X/K} (\C_p) \ar[d]^{\alpha^{\tau}} \ar[r] & NS (\oX)_{\tors} \ar@{=}[d] \ar[r] & 0 \\
0 \ar[r] & \Hom^0_c (\pi^{\abb}_1 (\oX) , \eo^*) \ar[r] & \Hom_c (\pi^{\abb}_1 (\oX) , \eo^*) \ar[r] & NS (\oX)_{\tors} \ar[r] & 0 
}}
\end{equation}
and
\begin{equation}
  \label{eq:25}
 \vcenter{\def\objectstyle{\scriptstyle}
 \xymatrix@C=10pt{
0 \ar[r] & \Pic^{\tau}_{X/K} (\C_p)_{\tors} \ar[d]^{\wr\,\alpha^{\tau}_{\tors}} \ar[r] & \Pic^{\tau}_{X/K} (\C_p) \ar[d]^{\alpha^{\tau}} \ar[r]^{\log} & \Lie \Pic^{\tau}_{X/K} (\C_p) \ar@{^{(}->}[d]^{\Lie \alpha^{\tau} = \Lie \alpha} \ar@{=}[r] & H^1 (X, \Oh) \otimes \C_p \ar[r] & 0 \\
0 \ar[r] & \Hom_c (\pi^{\abb} _1 (\oX), \mu ) \ar[r] & \Hom_c (\pi^{\abb}_1 (\oX) , \eo^*) \ar[r] & \Hom_c (\pi^{\abb}_1 (\oX) , \C_p) \ar@{=}[r] & H^1 (\oX , \Q_p) \otimes \C_p \ar[r] & 0 \; .
}}
\end{equation}
The map $\alpha^{\tau}$ is an injective and locally analytic homomorphism of $p$-adic Lie groups. Its restriction $\alpha^{\tau}_{\tors}$ to torsion subgroups is the inverse of the Kummer isomorphism:
\[
\scriptstyle i_X : \Hom_c (\pi^{\abb}_1 (\oX) , \mu) = \varinjlim_{N} H^1 (\oX , \mu_N) \silo \varinjlim_N H^1 (\oX , \Ge_m)_N = \Pic^{\tau}_{X/K} (\C_p)_{\tors} \; .
\]
The map $\alpha^{\tau}$ induces a topological isomorphism of complete topological groups
\[
\alpha^{\tau} : \Pic^{\tau}_{X/K} (\C_p) \silo Ch (\pi^{\abb}_1 (\oX)) \; .
\]
\end{theorem}

\begin{proof}
  We first note that $\Pic_{X/K} (\oK)_N = \Pic_{X/K} (\C_p)_N$ and $\Pic^{\tau}_{X/K} (\oK)_N = \Pic^{\tau}_{X/K} (\C_p)_N$ because this holds for $\Pic^0$ and because $NS (X_{\C_p}) = NS (\oX)$.

As $\Pic^0_{X/K} (\C_p)$ is divisible, the sequence (\ref{eq:23}) gives a short exact sequence
\[
0 \longrightarrow \Pic^0_{X/K} (\C_p)_{\tors} \longrightarrow \Pic^{\tau}_{X/K} (\C_p)_{\tors} \longrightarrow NS (\oX)_{\tors} \longrightarrow 0 \; .
\]
We claim that the following diagram with exact lines commutes:
\begin{equation}
  \label{eq:26}
\def\objectstyle{\scriptstyle}
\def\labelstyle{\scriptstyle}
\vcenter{\xymatrix@C=10pt{
0 \ar[r] & \Pic^0_{X/K} (\C_p)_{\tors} \ar[d]^{\alpha_{\tors}} \ar[r] & \Pic^{\tau}_{X/K} (\C_p)_{\tors} \ar[d]^{i^{-1}_X} \ar[r] & NS (\oX)_{\tors} \ar@{=}[d] \ar[r] & 0 \\
0 \ar[r] & \Hom^0_c (\pi^{\abb}_1 (\oX) , \mu) \ar[r] & \Hom_c (\pi^{\abb}_1 (\oX) , \mu) \ar[r]  & NS (\oX)_{\tors} \ar[r] & 0 \; .
}}
\end{equation}
If we make the identifications explicit which define the maps of the left square, we see that on $N$-torsion it is the outer rectangle of the following diagram
\[
\def\objectstyle{\scriptstyle}
\def\labelstyle{\scriptstyle}
\xymatrix@C=5pt{
\Pic^0_{X / K} (\C_p)_N \ar@{=}[r] \ar[d]_{\alpha_{\tors}} \ar@{}[dr] |{\fbox{\tiny 1}}& \hA_N \ar[d]_{\alpha_{\tors}} \ar@{=}[r] \ar@{}[dr] |{\fbox{\tiny 2}} & B_N \ar@{^{(}->}[rr] & & H^1 (\oX , \Ge_m)_N \ar[d]^{\wr i^{-1}_X} \\
\Hom^0_c (\pi^{\abb}_1 (\oX) , \mu_N) \ar@{=}[r] & \Hom_c (TA , \mu_N) \ar@{=}[r] & \Hom (A_N , \mu_N) \ar@{=}[r] & B_N \subset H^1 (\oX , \Ge_m)_N \ar[r]^-{\overset{i^{-1}_X}{\sim}} & H^1 (\oX , \mu_N) \; .
}
\]
Now, \fbox{\tiny 1} is commutative by definition, \fbox{\tiny 2} commutes because as noted in the proof of proposition \ref{t2} the restriction of $\alpha$ to $\hA_N$ is the map $\hA_N \to \Hom_c (TA , \mu_N) = \Hom (A_N , \mu_N)$ coming from Cartier duality. Hence the outer rectangle commutes as well.

The right square in diagram (\ref{eq:26}) is commutative since the second map in the exact sequence (\ref{eq:17}) is induced by $i_X$.

We now define $\alpha^{\tau}$ on \[
\Pic^{\tau}_{X/K} (\C_p) = \Pic^0_{X/K} (\C_p) + \Pic^{\tau}_{X/K} (\C_p)_{\tors}\] 
by setting it equal to $\alpha$ on $\Pic^0_{X/K} (\C_p)$ and to $i^{-1}_X$ on $\Pic^{\tau}_{X/K} (\C_p)_{\tors}$. This is well defined since by the commutativity of (\ref{eq:26}) the maps $\alpha$ and $i^{-1}_X$ agree on 
\[
\Pic^0_{X/K} (\C_p) \cap \Pic^{\tau}_{X/K} (\C_p)_{\tors} = \Pic^0_{X/K} (\C_p)_{\tors} \; . 
\]
The remaining assertions follow without difficulty. Note that $\Pic^{\tau}_{X/K} (\C_p)_{\tors}$ carries the discrete topology as a subspace of $\Pic^{\tau}_{X/K} (\C_p)$ since it is the kernel of the locally topological log-map.
\end{proof}

%\newpage
%\input{sec5}
\section{Categories of trivializing ``coverings''}
\label{sec:5}

In the next section we will single out classes of vector bundles on curves and abelian varieties to which $p$-adic representations can be attached. Our constructions require certain auxiliary categories of trivializing ``coverings''. These will be introduced and studied in the present section.

So let $\oX$ be a smooth connected projective curve over $\oK = \oQ_p$ and fix a point $\ox \in \oX (\oK)$. Let $\oeX$ be a finitely presented, flat and proper scheme over $\spec \eo_{\oK} = \spec \oZ_p$ such that $\oX = \oeX \otimes_{\eo_{\oK}} \oK$. Note that $\eo_{\oK}$ is not Noetherian.

We now introduce the category $\eT = \eT_{\oeX}$. Objects of $\eT$ are finitely presented, proper $G$-equivariant morphisms over $\spec \eo_{\oK}$
\[
\pi : \Yh \longrightarrow \oeX
\]
where $G$ is a finite (abstract) group which acts $\eo_{\oK}$-linearly from the left on $\Yh$ and trivially on $\oeX$. 
Moreover, we assume that $\Yh_{\oK} \rightarrow \oX$ is an (\'etale) $G_{\oX}$-torsor. 

It follows that $\Yh_{\oK}$ is a finite \'etale covering of $\oX$ and that every connected component $Y^0$ of $\Yh_{\oK}$ is a connected torsor over $\oX$ (hence a Galois covering) with group $G_{Y^0}$, the stabilizer of $Y^0$ in $G$. Hence for every $Y^0$ the group $G_{Y^0}$ acts faithfully and transitively on $\Mor_{\oX} (\ox , Y^0)$, the set of $\oK$-valued points of $Y^0$ lying over $\ox$.

Note that by our assumptions the structural morphism
\[
\lambda : \Yh \longrightarrow \spec \eo_{\oK}
\]
is finitely presented and proper with smooth generic fibre.

A morphism from the $G_1$-equivariant morphism $\pi_1 : \Yh_1 \to \oeX$ to the $G_2$-equivariant morphism $\pi_2 : \Yh_2 \to \oeX$ in $\eT$ is given by a commutative diagram:
\[
\xymatrix{
\Yh_1 \ar[rr]^{\varphi} \ar[dr]_{\pi_1} & & \Yh_2 \ar[dl]^{\pi_2}\\
 & \oeX &
}
\]
and a homomorphism of groups
\[
\gamma : G_1 \longrightarrow G_2
\]
such that $\varphi$ is $G_1$-equivariant if $G_1$ acts on $\Yh_2$ via $\gamma$.

{\bf Convention} If such a morphism $(\varphi , \gamma)$ exists, we say that $\pi_1 : \Yh_1 \longrightarrow \oeX$ dominates $\pi_2 : \Yh_2 \to \oeX$.

We now define two full subcategories
\[
\eT^{\good}_{ss} \subset \eT^{\good} \subset \eT
\]
of $\eT$ by imposing further conditions on the objects of $\eT$: For $\eT^{\good}$, we require that the structural morphism $\lambda : \Yh \to \spec \eo_{\oK}$ is flat and satisfies $\lambda_* \Oh_{\Yh} = \Oh_{\spec \eo_{\oK}}$ universally. In particular this implies that $\Yh_{\oK}$ is connected and hence that $\pi_{\oK} : \Yh_{\oK} \to \oX$ is a Galois covering with group $G$. 

For $\eT^{\good}_{ss}$ we require in addition that $\Pic^0_{\Yh / \eo_{\oK}}$ is representable by a semiabelian scheme finitely presented over $\eo_{\oK}$ (with generic fibre an abelian variety).

The following category is also of interest: $\eT_{\fin}$ is the full subcategory of $\eT$ of objects $\pi : \Yh \to \oeX$ where $\pi$ is a finite morphism. Set $\eT^{\good}_{\fin} = \eT^{\good} \cap \eT_{\fin}$. 

Finally, for an abelian scheme $\Ah / \eo_K$ we let $\eT_{\oAh}$ be the category whose objects are the $N$-multiplication morphism $\pi_N = N_{\oAh} : \oAh \to \oAh$ for some $N \ge 1$ where $\oAh = \Ah \otimes \eo_{\oK}$. Morphisms are commutative diagrams
\[
\xymatrix{
\oAh \ar[rr]^{N/M} \ar[dr]_{\pi_N} & & \oAh \ar[dl]^{\pi_M}\\
& \oAh & 
}
\]
where $M$ divides $N$. Note that the generic fibre $(\pi_N)_{\oK} : \oA \to \oA$ of \linebreak
$\pi_N : \oAh \to \oAh$ is Galois over $\oA = \Ah \otimes \oK$ with Galois group $G = A_N (\oK)$ acting by translation. By properness we have $G = \Ah_N (\eo_{\oK})$, and hence $G$ acts on $\pi_N : \oAh \to \oAh$ over $\oAh$. Note that the morphism $N / M$ becomes $G = A_N (\oK)$-equivariant if $A_N (\oK)$ acts on $\pi_M : \oAh \to \oAh$ via the homomorphism
\[
N / M : A_N (\oK) \longrightarrow A_M (\oK)\; .
\]
We need the following facts about these categories

\begin{theorem}
  \label{t13}
{\bf a} Finite direct products exist in $\eT = \eT_{\oeX}$ and $\eT_{\fin}$. \\
{\bf b} Any finite number of objects in $\eT = \eT_{\oeX}$ are dominated by a common object in $\eT^{\good}$ and even in $\eT^{\good}_{ss}$.\\
{\bf c} Any finite number of objects in $\eT_{\fin}$ is dominated by an object of $\eT^{\good}_{\fin}$.\\
{\bf d} Any finite number of objects in $\eT_{\oAh}$ is dominated by an object of $\eT_{\oAh}$.\\
{\bf e} Let $\oX_1$ and $\oX_2$ be curves as above with finitely presented flat and proper models $\oeX_1$ and $\oeX_2$ over $\eo_{\oK}$. For every $\eo_{\oK}$-morphism $f : \oeX_1 \to \oeX_2$, pullback along $f$ induces a functor
\[
f^{-1} : \eT_{\oeX_2} \longrightarrow \eT_{\oeX_1}
\]
which commutes with finite direct products and maps $\eT_{\oeX_2 , \fin}$ into $\eT_{\oeX_1 , \fin}$.\\
{\bf f} For every $\eo_{\oK}$-morphism $f : \oeX \to \oAh$ pullback along $f$ induces a functor
\[
f^{-1} : \eT_{\oAh} \longrightarrow \eT_{\oeX} \; .
\]
\end{theorem}

\begin{rem}
  For the constructions of the present paper, the assertion about $\eT^{\good}_{ss}$ in {\bf b} is not required. However, we expect that it will become important for a deeper study of the categories of vector bundles introduced in the next section.
\end{rem}

\begin{proof}
  {\bf a} Let $\pi_1 : \Yh_1 \to \oeX$ and $\pi_2 : \Yh_2 \to \oeX$ be $G_1$- resp. $G_2$-equivariant morphisms that are objects in $\eT$. The fibre product
\[
\pi : \Yh_1 \times_{\oeX} \Yh_2 \longrightarrow \oeX
\]
is a finitely presented, proper and $G = G_1 \times G_2$-equivariant morphism. Its generic fibre is a $G$-torsor over $\oX$. 

We now prove parts {\bf b} and {\bf c} of theorem \ref{t13}. Because of {\bf a} if suffices to show that any object $\pi : \Yh \to \oeX$ in $\eT$ (resp. $\eT_{\fin}$) is dominated by an object of $\eT^{\good}$ and even of $\eT^{\good}_{ss}$ (resp. by an object of $\eT^{\good}_{\fin}$). Using \cite{EGAIV} \S\,8, in particular (8.8.3) and (8.10.5), together with \cite{EGAIV}, (17.7.8) we get the following situation:

There is a finite extension $K / \Q_p$ and a smooth projective geometrically connected curve $X / K$ with a flat and proper model $\eX / \eo_K$ such that $\oeX = \eX \otimes_{\eo_K} \eo_{\oK}$ and $\oX = X \otimes_K \oK$. There is a $K$-rational point $x \in X (K) = \eX (\eo_K)$ which induces the given point $\ox \in \oX (\oK)$. Moreover there is a proper (resp. finite) morphism:
\[
\pi_0 : \Yh_0 \to \eX
\]
of schemes over $\eo_K$ which is equivariant with respect to an $\eo_K$-linear $G$-action on $\Yh_0$, such that $\Yh_{0 K}$ is a $G$-torsor over $X$. All connected components $Y^0_0$ of $\Yh_{0 K}$ are geometrically irreducible, have a $K$-rational point over $x$ and are such that
\[
\pi_0 \, |_{Y^0_0} : Y^0_0 \longrightarrow X
\]
is Galois with group $G_{Y^0_0}$. Finally there is a commutative diagram where $\beta$ is $G$-equivariant
\[
\xymatrix{
\Yh_0 \otimes \eo_{\oK} \ar[r]^{\overset{\beta}{\sim}} \ar[d]_{\pi_0 \otimes \id} & \Yh \ar[d]^{\pi} \\
\eX \otimes \eo_{\oK}  \ar@{=}[r] & \oeX \; .
}
\]
The connected components of $\Yh_{\oK}$ are the spaces $Y^0 = \beta (Y^0_0 \otimes_K \oK)$ and we have $G_{Y^0} = G_{Y^0_0}$. For later purposes note that the following arguments only require that $\eo_K$ is a discrete valuation ring. So instead of a finite extension $K / \Q_p$ we could work over a finite extension of the maximal unramified extension of $\Q_p$. 

Now we choose a connected component $Z_0$ of $\Yh_{0K}$ and let $\Zh_0$ be its closure in $\Yh_0$ with the reduced induced scheme structure. Then the action of $H = G_{Z_0}$ on $\Yh_0$ induces an action on $\Zh_0$ such that the closed immersion
\[
i : \Zh_0 \hookrightarrow \Yh_0
\]
is  $H$-equivariant. It follows that the composition
\[
\pi_{\Zh_0} : \Zh_0 \overset{i}{\hookrightarrow} \Yh_0 \xrightarrow{\pi_0} \eX
\]
is $H$-equivariant, proper (resp. finite) and finitely presented. The scheme $\Zh_0$ is integral and its generic fibre $Z_0 = \Zh_0 \otimes K$ is Galois over $X$ with group $H$. As $\Zh_0$ is of finite type over $\spec \eo_K$ it is an excellent scheme. Hence the normalization $\Wh_0$ of $\Zh_0$ in its function field is finite over $\Zh_0$. By functoriality of the normalization, the group $H$ acts on $\Wh_0$. Hence the composition
\[
\pi_{\Wh_0} : \Wh_0 \longrightarrow \Zh_0 \xrightarrow{\pi_{\Zh_0}} \eX
\]
is $H$-equivariant, proper (resp. finite) and finitely presented. Moreover $\Wh_0$ is integral and its generic fibre $\Wh_{0 K} = \Zh_{0K} = Z_0$ is Galois over $X$ with group $H$. By construction $\Wh_{0K}$ contains a $K$-rational point which by properness extends to a section of $\Wh_0$ over $\spec \eo_K$. Apart from being proper, the structural morphism $\Wh_0 \to \spec \eo_K$ is also flat by \cite{Ha} III, 9.7. Hence, $\Wh_0$ maps surjectively to $\spec \eo_K$. Hence condition (5.6.6.2) in \cite{EGAIV} (5.6.6) is verified and it follows that we have $\dim \Wh_0 = 2$. From \cite{EGAIV} (13.1.1) it follows that all irreducible components of the special fibre of $\lambda_0 : \Wh_0 \to \spec \eo_K$ are at least one-dimensional. Because of $\dim \Wh_0 = 2$, all these components must be one-dimensional. Since $\lambda_0$ is proper the sheaf $\lambda_{0*} \Oh_{\Wh_0}$ on $\spec \eo_K$ is coherent and hence given by the finitely generated $\eo_K$-module
\[
\Gamma (\spec \eo_K , \lambda_{0*} \Oh_{\Wh_0}) = \Gamma (\Wh_0 , \Oh) \; .
\]
Since $\Wh_0$ is integral, this module is torsion free, hence free so that $\lambda_{0*} \Oh_{\Wh_0} = \Oh^r_{\spec \eo_K}$ for some $r \ge 0$. Alternatively this follows from the flatness of $\lambda_0$. Since the generic fibre $\Wh_{0K} = Z_0$ is connected it follows from flat base change that $r = 1$. Hence we have $\lambda_{0*} \Oh_{\Wh_0} = \Oh_{\spec \eo_K}$. We now quote the following special case of a result of Raynaud.

{\bf Theorem (Raynaud)} {\it  Let $S$ be the spectrum of a discrete valuation ring, $f : V \to S$ a proper flat morphism all of whose fibres have one-dimensional irreducible components. Assume that $V$ is normal and that $f_* \Oh_V = \Oh_S$. If $f$ has a section, then $f$ is cohomologically flat in degree zero i.e. the formation of $f_* \Oh_V$ commutes with arbitrary base change. Moreover $\Pic^0_{V/S}$ exists as a separated group scheme.}

This follows from \cite{R} Th\'eor\`eme (8.2.1): Condition $(N)^*$ of \cite{R} (6.1.4) is satisfied because $V$ is normal (see (6.1.6) in \cite{R}) and because $f_* \Oh_V = \Oh_S$ by assumption. Moreover the existence of a section implies that at least one irreducible component of the special fibre $f^{-1} (s)$ has total multiplicity one in $f^{-1} (s)$. 

We now define an object $\pi^{\good} : \Yh^{\good} \to \oeX$ of $\eT^{\good}$ (resp. $\eT^{\good}_{\fin}$) by setting $\Yh^{\good} = \Wh_0 \otimes_{\eo_K} \eo_{\oK}$ and $\pi^{\good} = \pi_{\Wh_0} \otimes \eo_{\oK}$. Then $\pi^{\good}$ is $H$-equivariant and $\Yh^{\good}_{\oK}$ is Galois over $\oX$ with group $H$. The structural morphism \\
$\lambda^{\good} = \lambda_0 \otimes \eo_{\oK} : \Yh^{\good} \to \spec \eo_{\oK}$ is flat because $\lambda_0$ is and we have $\lambda^{\good}_* \Oh_{\Yh^{\good}} = \Oh_{\spec \eo_{\oK}}$ universally since we have $\lambda_{0*} \Oh_{\Wh_0} = \Oh_{\spec \eo_K}$ and since $\lambda_0$ is cohomologically flat in degree zero by Raynaud's theorem. By construction we have a commutative diagram
\[
\xymatrix{
\Wh_0 \ar[rr]^{\kappa} \ar[dr]_{\pi_{\Wh_0}} & & \Yh_0 \ar[dl]^{\pi_0} \\
 & \eX & 
}
\]
where $\kappa$ is a finite $H$-equivariant morphism. Extending the base to $\eo_{\oK}$ gives the desired $H$-equivariant morphism $\kappa \otimes \eo_{\oK} : \Yh^{\good} \to \Yh$ over $\oeX$ which together with the inclusion $H \hookrightarrow G$ defines a morphism from $\pi^{\good}$ to $\pi$ in $\eT$ (resp. $\eT_{\fin}$).

It remains to show that the objects of $\eT$ can even be dominated by objects of $\eT^{\good}_{ss}$.  Let $g = \dim_K H^1 (Z_0 , \Oh)$ be the genus of $Z_0$ and let $N (g)$ be the number depending only on $g$ mentioned in \cite{AW} Theorem 2.8. We may choose $K$ so big that for some prime number $l > N (g) , l \neq p$ we have
\[
\dim_{\F_l} H^1 (Z_0 , \Oh^*)_l = 2g \; .
\]
Namely, since $K$ was chosen big enough so that $Z_0$ is geometrically connected, passing to the finite extension $K'$ obtained by adjoining the $l$-torsion points of $\Pic^0_{Z_0 / K}$ to $K$ leaves $Z_0 \otimes_K K'$ geometrically connected. Let $L$ be the maximal unramified extension of $K'$ in $\oK$. We pass to the base change 
 $\Yh_{0,\eo_L}$ and proceed as before by normalizing the closure of the connected 
component $Z_{0 L} $ in $\Yh_{0 ,\eo_L}$. Abusing notation, we call the
resulting scheme also $\Wh_0$. The same arguments as before show
that $\Wh_0 \rightarrow \Yh_{0, \eo_L}$ is  a finite $H$-equivariant morphism,
which is a Galois covering on the generic fibre.

As in \cite{Ar} we call a Noetherian, normal, connected and excellent scheme of dimension $2$ a surface. We have seen that $\Wh_0$ is a surface in this sense. A surface has finitely many closed points where the local ring is not regular. These finitely many singularities are therefore contained in the special fibre of $\Wh_0$ over $\spec \eo_L$. They are permuted by the action of $H$. Define a sequence of surfaces
\[
\Wh_0 \xleftarrow{f_1} \Wh_1 \xleftarrow{f_2} \Wh_2 \xleftarrow{} \ldots
\]
 as follows: Let $S_i \subset \Wh_i$ be the singular locus with the reduced structure. Then $\Wh_{i+1}$ is the normalization of the blow-up of $S_i$ in $\Wh_i$. Since blowing-up and normalizing are functorial it follows that each $\Wh_i$ inherits an $H$-action such that the morphism $f_i$ become $H$-equivariant. According to a theorem of Lipman \cite{Li}, $\Wh_n$ is regular for sufficiently large $n$. Fix some such $n$ and set $\tilde{\Wh_0} = \Wh_n$. Since the singularities all lie in the special fibres of the $\Wh_i$'s we have $\tilde{\Wh}_{0L} = \Wh_{0L} = Z_{0L}$. The morphism
\[
\pi_{\tilde{\Wh_0}} : \tilde{\Wh}_0 \longrightarrow \Wh_0 \xrightarrow{\pi_{\Wh_0}} \eX_{\eo_L}
\]
is proper but in general it will no longer be finite even if $\pi$ and hence $\pi_{\Wh_0}$ is finite. Applying the above arguments to $\tilde{\Wh}_0$ instead of $\Wh_0$, we see that the structural morphism $\tilde{\lambda}_0 : \tilde{\Wh}_0 \to \spec \eo_L$ satisfies $\tilde{\lambda}_{0*} \Oh_{\tilde{\Wh}_0} = \Oh_{\spec \eo_L}$ universally. In particular the closed fibre of $\tilde{\Wh}_0$ is geometrically connected. All of its irreducible components are one-dimensional. Moreover, since one of them has multiplicity one due to the existence of a section, the greatest common divisor of all their multiplicities is one. Since the closed and the generic fibre of $\tilde{\Wh}_0$ have the same $h^0$, they have the same $h^1$ since the Euler characteristic of the structure sheaf is constant in a proper flat family over a Noetherian connected base \cite{Mu} II, \S\,5. Hence the closed fibre of $\tilde{\Wh}_0$ has genus $g$ as well and we may finally apply \cite{AW} Theorem 2.8 to conclude that the $\Pic^0$ of the closed fibre of $\tilde{\Wh}_0$ has trivial unipotent radical. By \cite{BLR} Theorem 4 in 9.5, $\Pic^0_{\tilde{\Wh}_0 / \eo_L}$ exists as a separated scheme which is isomorphic to the identity component of the N\'eron model of the abelian variety $\Pic^0_{\tilde{\Wh}_{0 L} / L} = \Pic^0_{Z_{0L} / L}$. In conclusion, we have seen that $\Pic^0_{\tilde{\Wh}_0 / \eo_L}$ is a semiabelian scheme over $\eo_L$. 

Performing a base extension to $\eo_{\oK}$ as before we obtain the desired object of $\eT^{\good}_{ss}$ dominating $\pi : \Yh \to \oeX$. 

Now consider assertion {\bf d} of theorem \ref{t13}. If $\pi_{N_i} : \oAh \to \oAh$ are objects of $\eT_{\Ah}$ where $N_1 , \ldots , N_r \ge 1$ are integers, set $N = N_1 \cdots N_r$. Then $\pi_N : \oAh \to \oAh$ maps via $N / N_i$ to $\pi_{N_i} : \oAh \to \oAh$. \\
{\bf e} Let $\pi_2 : \Yh_2 \to \oeX_2$ be an object of $\eT_{\oeX_2}$. Its pullback \\
$\pi_1 : \Yh_1 := \oeX_1 \times_{\oeX_2} \Yh_2 \to \oeX_1$ along $f$ is again finitely presented and proper. If $G$ is the group for $\pi_2$ then we have a natural $\eo_{\oK}$-linear $G$-action on $\Yh_1$ with respect to which $\pi_1$ is $G$-equivariant. It is clear that $\pi_{1 \oK}: \Yh_{1 \oK} \rightarrow \oX_1$ is a $G$-torsor. 

Thus pullback along $f$ induces a functor from $\eT_{\oeX_2}$ into $\eT_{\oeX_1}$ and it is clear that if $\pi_2$ is finite, $\pi_1$ is finite as well. That finite direct products are preserved is equally clear by the explicit construction of the direct product in $\eT$ given in the proof of part {\bf a}.

{\bf f} Same argument as in {\bf e} with $\pi_2$ replaced by an $N$-multiplication morphism $\pi_N : \oAh \to \oAh$. 
\end{proof}
%\newpage
%\input{sec6}

\section{Vector bundles that give rise to $p$-adic representations}
\label{sec:6}

In this section we define and study certain categories of vector bundles to which one can attach $p$-adic representations of the geometric fundamental group of the base variety. This base will be either a curve or an abelian variety.

A vector bundle $E \to X$ of rank $r$ is a scheme over $X$ which Zariski-locally on $X$ is isomorphic to $r$-dimensional affine space $\A^r$, together with linear transition functions. Its sheaf of section, which we denote by $\Oh (E)$ is a locally free sheaf of rank $r$ on $X$. 

Consider first the case of a smooth projective curve $\oX / \oK$ as in the preceeding section i.e. with a finitely presented flat and proper model $\oeX$ over $\eo_{\oK}$ and a rational point $\ox \in \oX (\oK)$. Let $\Vec_{\eX_{\eo}}$ be the category of vector bundles $E$ on $\eX_{\eo} = \oeX \otimes \eo$. For every $n \ge 1$ we denote the restriction of $E$ to $\eX_n = \eX_{\eo} \otimes \eo_n = \oeX \otimes \eo_n$ by $E_n = E \otimes \eo_n$, where $\eo_n = \eo / p^n \eo = \eo_{\oK} / p^n \eo_{\oK}$. 

\begin{definition}
  The category $\eB = \eB_{\eX_{\eo}}$ (resp. $\eB^{\fin} = \eB^{\fin}_{\eX_{\eo}}$) is the full subcategory of $\Vec_{\eX_{\eo}}$ whose objects $E$ satisfy the following condition:\\
For every $n \ge 1$ there is an object $\pi : \Yh \to \oeX$ of $\eT_{\oeX}$ (resp. $\eT_{\fin}$) such that $\pi^*_n E_n$ is a trivial bundle on $\Yh_n$.\\
Here $\pi_n , \Yh_n$ are the reductions $\mod p^n$ of $\pi$ and $\Yh$.
\end{definition}

It is clear that $\eB^{\fin}$ is a full subcategory of $\eB$.

\begin{definition}
For an abelian scheme $\Ah / \eo_K$ let $\Vec_{\Ah_{\eo}}$ be the category of vector bundles on $\Ah_{\eo}$. The full subcategory $\eB = \eB_{\Ah_{\eo}}$ of $\Vec_{\Ah_{\eo}}$ is defined similarly to $\eB_{\eX_{\eo}}$ using the category $\eT_{\oAh}$.
\end{definition}

Note that if $E$ is a bundle in $\eB_{\eX_{\eo}}$ then all bundles $E'$ in $\Vec_{\eX_{\eo}}$ which are isomorphic to $E$ are also objects of $\eB_{\eX_{\eo}}$, similarly for $\eB_{\Ah_{\eo}}$. 

For a topological group $\Sigma$ let $\Rep_{\Sigma} (\eo)$ be the category of continuous representations of $\Sigma$ on free $\eo = \eo_{\C_p}$-modules of finite rank.

We are now going to construct functors $\rho$ from $\eB_{\eX_{\eo}}$ into $\Rep_{\pi_1 (\oX, \ox)} (\eo)$ and from $\eB_{\Ah_{\eo}}$ into $\Rep_{\pi_1 (\oA , 0)} (\eo)$.

Via the inclusion $\oX (\oK) = \oeX (\eo_{\oK}) \subset \eX_{\eo} (\eo)$ we may view the point $\ox$ as a section $x_{\eo} : \spec \eo \to \eX_{\eo}$. For a bundle $E$ on $\eX_{\eo}$ we write $E_{x_{\eo}} = x^*_{\eo} E$ viewed as a free $\eo$-module of rank $r = \rank E$. Via reduction $\eX_{\eo} (\eo) \to \eX_{\eo} (\eo_n) = \eX_n (\eo_n)$ the section $x_{\eo}$ gives rise to a morphism
\[
x_n : \spec \eo_n \longrightarrow \spec \eo \xrightarrow{x_{\eo}} \eX_{\eo}
\]
and we set $E_{x_n} = x^*_n E = E_{x_{\eo}} \otimes \eo_n$ viewed as a free $\eo_n$-module of rank $r$. Note that we have
\[
E_{x_{\eo}} = \varprojlim_n E_{x_n}
\]
as topological $\eo$-modules, the topology on $E_{x_n}$ being the discrete one.

We now recall some facts about the algebraic fundamental group of $\oX$. Let $F$ be the functor from the category of finite \'etale coverings $Y$ of $\oX$ to the category of sets, defined by
\[
F (Y) = \Mor_{\oX} (\ox, Y) \; .
\]
This functor is strictly pro-represented by a projective system \\
$\tilde{\oX} = (Y_i , \oy_i , \phi_{ij})_{i \in I}$ of pointed Galois coverings of $\oX$. Thus $I$ is a directed set, $\oy_i \in Y_i (\oK)$ are points over $\ox$. Moreover, for $i \ge j$ the map $\phi_{ij} : Y_i \to Y_j$ is an epimorphism over $\oX$ such that $\phi_{ij} (\oy_i) = \oy_j$ and the natural map
\[
\varinjlim_{i} \Mor_{\oX} (Y_i , Y) \longrightarrow F (Y)
\]
induced by evaluation on the $\oy_i$'s is a bijection for every $Y$. Then we have:
\[
\pi_1 (\oX , \ox) = (\varprojlim_i \Aut_{\oX} (Y_i))^{\op} = \varprojlim_i \Aut_{\oX} (Y_i)^{\op}
\]
where $\Aut_{\oX} (Y_i)$ is the group of $\oX$-automorphisms of $Y_i$ acting on the left. The transition maps
\[
\psi_{ij} : \Aut_{\oX} (Y_i) \longrightarrow \Aut_{\oX} (Y_j) \quad \mbox{for} \; i \ge j
\]
are defined by the condition that $\psi_{ij} (\sigma) (\oy_j) = \phi_{ij} \verk \sigma \verk \oy_i$. Note here that $Y_i$ being Galois over $\oX$, the group $\Aut_{\oX} (Y_i)$ acts simply transitively on $F (Y_i)$.

{\bf Construction} For a vector bundle $E$ in $\eB_{\eX_{\eo}}$ resp. $\eB_{\Ah_{\eo}}$ we define a continuous representation
\[
\rho_E : \pi_1 (\oX, \ox) \longrightarrow \Aut_{\eo} E_{x_{\eo}} \quad \mbox{resp.} \quad \rho_E : \pi_1 (\oA , 0) \longrightarrow  \Aut_{\eo} E_{x_{\eo}}
\]
as the projective limit of certain continuous representations
\[
\rho_{E,n} : \pi_1 (\oX , \ox) \longrightarrow \Aut_{\eo_n} E_{x_n}  \quad \mbox{resp.} \quad \rho_{E,n} : \pi_1 (\oA , 0) \longrightarrow \Aut_{\eo_n} E_{x_n} \; .
\]
We begin with the curve case. For any given $n \ge 1$ there exists an object $\pi : \Yh \to \oeX$ in $\eT_{\oeX}$ such that $\pi^*_n E_n$ is trivial. Because of theorem \ref{t13} part {\bf b} we may choose $\pi : \Yh \to \oeX$ to lie in $\eT^{\good}_{\oeX}$. In particular $\Yh_{\oK}$ is then Galois over $\oX$. Fix a point $\oy$ in $\Yh_{\oK} (\oK)$ lying over $\ox \in \oX (\oK)$. It determines a map of projective systems $\tilde{\oX} \to \Yh_{\oK}$ represented by a morphism $a_i : Y_i \to \Yh_{\oK}$ over $\oX$ with $\oy_i \mapsto \oy$.

Consider the induced epimorphism:
\[
\psi_i : \Aut_{\oX} Y_i \twoheadrightarrow \Aut_{\oX} \Yh_{\oK} =: G 
\]
defined by $\psi_i (\sigma) (\oy) = a_i \verk \sigma \verk \oy_i$. The composition
\[
\varphi_{\oy} : \pi_1 (\oX , \ox) \twoheadrightarrow \Aut^{\op}_{\oX} Y_i \overset{\psi_i}{\twoheadrightarrow} G^{\op}
\]
is independent of the choice of $i \in I$ but depends on $\oy$. By definition of $\eT^{\good}_{\oeX}$, the action of the group $G$ on $\Yh_{\oK}$ extends to an action on $\Yh$ over $\oeX$ and hence reduces to an action on $\Yh_n$ over $\eX_n$. Therefore $G$ acts $\eo_n$-linearly from the left on $\Gamma (\Yh_n , \pi^*_n E_n)$. For a given $g$ in $G$ we denote by $g^*$ the induced automorphism of $\Gamma (\Yh_n , \pi^*_n E_n)$. Composing with $\varphi_{\oy}$ we get a homomorphism
\[
\pi_1 (\oX , \ox) \xrightarrow{\varphi_{\oy}} G^{\op} \longrightarrow \Aut_{\eo_n} \Gamma (\Yh_n , \pi^*_n E_n) \; .
\]
As $\Yh \to \oeX$ lies in $\eT^{\good}_{\oeX}$, the structural morphism $\lambda : \Yh \to \spec \eo_{\oK}$ satisfies $\lambda_* \Oh_{\Yh} = \Oh_{\spec_{\eo_{\oK}}}$ universally. In particular, we have $\lambda_* \Oh_{\Yh_n} = \Oh_{\spec \eo_n}$  and therefore equality $\Gamma (\Yh_n , \Oh_{\Yh_n}) = \eo_n$. The bundle $\pi^*_n E_n$ being trivial, it follows that the pullback map under $y_n : \spec \eo_n \to \Yh_n$
\[
y^*_n : \Gamma (\Yh_n , \pi^*_n E_n) \silo \Gamma (\spec \eo_n , y^*_n \pi^*_n E_n) = \Gamma (\spec \eo_n , x^*_n E_n) = E_{x_n}
\]
is actually an isomorphism. The representation $\rho_{E,n}$ is defined to be the composition:
\[
\rho_{E,n} : \pi_1 (\oX , \ox) \xrightarrow{\varphi_{\oy}} G^{\op} \longrightarrow \Aut_{\eo_n} \Gamma (\Yh_n , \pi^*_n E_n) \xrightarrow{\overset{\mathrm{via}\,y^*_n}{\sim}} \Aut_{\eo_n} E_{x_n} \; .
\]
For vector bundles $E$ in $\eB_{\Ah_{\eo}}$ we proceed similarly. Recall that we have $\pi_1 (\oA , 0) = TA$ canonically. Let $x_{\eo} , x_n$ denote the zero sections of $\Ah_{\eo}$ resp. $\Ah_n = \Ah_{\eo} \otimes \eo_n = \oAh \otimes \eo_n$ and define $E_{x_{\eo}}$ and $E_{x_n}$ as before. For a given $n \ge 1$, by definition of $\eB_{\Ah_{\eo}}$, there exists some $N = N (n) \ge 1$ such that $\pi^*_{N,n} E_n$ is a trivial bundle. Here $\pi_{N,n}$ is the reduction $\mod p^n$ of the $N$-multiplication map $\pi_N = N_{\oAh} : \oAh \to \oAh$ on $\oAh$. The structural morphism $\lambda : \Ah \to \spec \eo_K$ satisfies $\lambda_* \Oh_{\Ah} = \Oh_{\spec \eo_K}$ universally  (use \cite{EGAIII} 7.8.6). Hence we have $\Gamma (\Ah_n , \Oh) = \eo_n$ and therefore the pullback map
\[
x^*_n : \Gamma (\Ah_n , \pi^*_{N,n} E_n) \silo \Gamma (\spec \eo_n , x^*_n E_n) = E_{x_n}
\]
is an isomorphism of free $\eo_n$-modules. Note that $\pi_{N,n} \verk x_n = x_n$. On \\
$\Gamma (\Ah_n , \pi^*_{N,n} E_n)$ the group $A_N (\oK)$ acts by translation. Define the representation $\rho_{E,n}$ as the composition:
\[
\rho_{E,n} : \pi_1 (\oA , 0) = TA \longrightarrow A_N (\oK) \longrightarrow \Aut_{\eo_n} \Gamma (\Ah_n , \pi^*_{N,n} E_n) \overset{\mathrm{via}\,x^*_n}{\silo} \Aut_{\eo_n} E_{x_n} \; .
\]

\begin{theorem}
  \label{t16}
For $E$ in $\eB_{\eX_{\eo}}$ resp. $\eB_{\Ah_{\eo}}$ the representations $\rho_{E_n}$ are independent of all choices and form a projective system when composed with the natural projection maps $\Aut_{\eo_{n+1}} E_{x_{n+1}} \to \Aut_{\eo_n} E_{x_n}$.
\end{theorem}

\begin{proof}
  We give the proof only for $E$ in $\eB_{\eX_{\eo}}$. The other case is similar. We first check that $\rho_{E,n}$ is independent of the point $\oy$. So let $\oy'$ be another point in $\Yh_{\oK} (\oK)$ lying over $\ox$. We then have two maps of projective systems $\tilde{\oX} \to \Yh_{\oK}$ corresponding to $\oy$ and $\oy'$. We may choose a common index $i \in I$ such that these maps are represented by morphisms over $\oX$:
\[
a_i : Y_i \longrightarrow \Yh_{\oK} \quad \mbox{and} \quad a'_i : Y_i \longrightarrow \Yh_{\oK}
\]
determined by $a_i (\oy_i) = \oy$ resp. $a'_i (\oy_i) = \oy'$. As $\Yh_{\oK}$ is Galois over $\oX$, there is a unique $\tau \in G = \Aut_{\oX} \Yh_{\oK}$ with $\tau \oy = \oy'$. The epimorphisms induced by $a_i$ and $a'_i$:
\[
\psi_i , \psi'_i : \Aut_{\oX} Y_i \longrightarrow G
\]
are related by conjugation with $\tau:$
\begin{equation}
  \label{eq:28}
  \psi'_i (\sigma) = \tau \psi_i (\sigma) \tau^{-1} \quad \mbox{for all} \; \sigma \in \Aut_{\oX} Y_i \; .
\end{equation}
Namely, $\psi_i (\sigma)$ and $\psi'_i (\sigma)$ are determined by the relations
\[
\psi_i (\sigma) (\oy) = a_i (\sigma \oy_i) \quad \mbox{and} \quad \psi'_i (\sigma) (\oy') = a'_i (\sigma \oy_i) \; .
\]
Moreover, the relation
\[
(\tau \verk a_i) (\oy_i) = \oy' = a'_i (\oy_i)
\]
implies that $\tau \verk a_i = a'_i$. Now, from the equalities
\[
(\psi'_i (\sigma) \tau) (\oy) = \psi'_i (\sigma) (\oy') = a'_i (\sigma \oy_i) = \tau (a_i (\sigma \oy_i)) = (\tau \psi_i (\sigma)) (\oy)
\]
relation (\ref{eq:28}) follows.

Let $y_n , y'_n$ be the images of $\oy$ and $\oy'$ under the reduction map
\[
\Yh_{\oK} (\oK) = \Yh (\eo_{\oK}) \longrightarrow \Yh (\eo_n) = \Yh_n (\eo_n) \; .
\]
We have $\tau y_n = y'_n$ i.e. the following diagram commutes
\[
\xymatrix{
\Yh_n \ar[dd]^{\wr}_{\tau} & \\
 & \spec \eo_n \ar[ul]_{y_n} \ar[dl]^{y'_n} \\
\Yh_n & 
}
\]
This diagram induces a commutative diagram of isomorphisms
\begin{equation}
  \label{eq:29}
  \vcenter{\xymatrix{
\Gamma (\Yh_n , \pi^*_n E_n) \ar[dr]^{\overset{y^*_n}{\sim}} & \\
 & \Gamma (\spec \eo_n , x^*_n E_n) = E_{x_n} \; .\\
\Gamma (\Yh_n , \pi^*_n E_n) \ar[uu]^{\tau^*}_{\wr} \ar[ur]^{\sim}_{y^{'*}_n} & 
}}
\end{equation}
An element $\sigma$ of $\Aut_{\oX}Y_i$ acts on $E_{x_n}$ by the automorphism $y^*_n \verk \psi_i (\sigma)^* \verk (y^*_n)^{-1}$ resp. $y^{'*}_n \verk \psi'_i (\sigma)^* \verk (y^{'*}_n)^{-1}$ if we base our construction on the point $\oy$ resp. $\oy'$.

Using (\ref{eq:28}) and (\ref{eq:29}) we see that these two automorphisms of $E_{x_n}$ are actually the same. It follows that $\rho_{E,n}$ does not depend on the choice of the auxiliary point $\oy$.

Now we check that $\rho_{E,n}$ is independent of the trivializing cover $\pi : \Yh \to \oeX$ in $\eT^{\good}_{\oeX}$. So let $\pi' : \Yh' \to \oeX$ be another object of $\eT^{\good}_{\oeX}$ such that $\pi^{'*}_n E_n$ is a trivial bundle. By theorem \ref{t13}, part {\bf b} we may assume that $\pi$ dominates $\pi'$. Thus there is a commutative diagram
\[
\xymatrix{
\Yh \ar[rr]^{\varphi} \ar[dr]_{\pi} & & \Yh' \ar[dl]^{\pi'}\\
 & \oeX & 
}
\]
and a homomorphism of the corresponding groups $\gamma : G \to G'$ such that $\varphi$ is $G$-equivariant if $G$ acts on $\Yh'$ via $\varphi$. 

Fix a point $\oy \in \Yh_{\oK} (\oK)$ over $\ox$ and set $\oy' = \varphi_{\oK} (\oy)$. This determines a commutative diagram of morphisms over $\oX$
\[
\xymatrix{
 & \Yh_{\oK} \ar[dd]^{\varphi_{\oK}} \\
Y_i \ar[ur]^{a_i} \ar[dr]_{a'_i} & \\
 & \Yh'_{\oK}
}
\]
where $a_i (\oy_i) = \oy$ and $a'_i (\oy_i) = \oy'$. The corresponding epimorphisms
\[
\psi_i : \Aut_{\oX} Y_i \longrightarrow G = \Aut_{\oX} \Yh_{\oK} \quad \mbox{and} \quad \psi'_i : \Aut_{\oX} Y_i \longrightarrow G' = \Aut_{\oX} \Yh'_{\oK}
\]
are related by the formulas
\begin{equation}
  \label{eq:30}
  \psi'_i (\sigma) \verk \varphi_{\oK} = \varphi_{\oK} \verk \psi_i (\sigma) \quad \mbox{for all} \; \sigma \; \mbox{in} \; \Aut_{\oX} Y_i
\end{equation}
and
\begin{equation}
  \label{eq:31}
  \psi'_i = \gamma \verk \psi_i \; .
\end{equation}
Namely, we have
\begin{eqnarray*}
  (\psi'_i (\sigma) \verk \varphi_{\oK}) (\oy) & = & \psi'_i (\sigma) (\oy') = a'_i (\sigma \oy_i) = \varphi_{\oK} (a_i (\sigma \oy_i)) \\
& = & (\varphi_{\oK} \verk \psi_i (\sigma)) (\oy)
\end{eqnarray*}
which gives (\ref{eq:30}) and
\begin{eqnarray*}
  \gamma (\psi_i (\sigma)) (\oy') & = & \gamma (\psi_i (\sigma)) (\varphi_{\oK} (\oy)) = (\gamma (\psi_i (\sigma)) \verk \varphi_{\oK}) (\oy) \\
& = & (\varphi_{\oK} \verk \psi_i (\sigma)) (\oy) \overset{(\ref{eq:30})}{=} (\psi'_i (\sigma) \verk \varphi_{\oK}) (\oy) = \psi'_i (\sigma) (\oy')
\end{eqnarray*}
which gives (\ref{eq:31}).

For every $\sigma$ in $\Aut_{\oX} Y_i$ we have the equality
\begin{equation}
  \label{eq:32}
  \psi'_i (\sigma) \verk \varphi = \varphi \verk \psi_i (\sigma) \quad \mbox{as morphisms} \; \Yh \to \Yh' \; .
\end{equation}
This follows from the $G$-equivariance of $\varphi$
\[
\gamma (\psi_i (\sigma)) \verk \varphi = \varphi \verk \psi_i (\sigma)
\]
together with relation (\ref{eq:31}).

The commutative diagram
\[
\xymatrix{\Yh_n \ar[dd]_{\varphi_n} & \\
 & \spec \eo_n \ar[ul]_{y_n} \ar[dl]^{y'_n} \\
\Yh'_n & 
}
\]
induces the commutative diagram of isomorphisms
\begin{equation}
  \label{eq:33}
  \vcenter{
\xymatrix{
\Gamma (\Yh_n , \pi^*_n E_n) \ar[dr]^{\overset{y^*_n}{\sim}} & \\
 & \Gamma (\spec \eo_n , x^*_n E_n) = E_{x_n} \; . \\
\Gamma (\Yh'_n , \pi^{'*}_n E_n) \ar[uu]^{\varphi^*_n}_{\wr} \ar[ur]^{y^{'*}_n}_{\sim}
}}
\end{equation}
An element $\sigma$ in $\Aut_{\oX} Y_i$ acts on $E_{x_n}$ by the automorphism $y^*_n \verk \psi_i (\sigma)^* \verk (y^*_n)^{-1}$ resp. $y^{'*}_n \verk \psi'_i (\sigma)^* \verk (y^{'*}_n)^{-1}$ if we base our construction on $(\Yh , \oy)$ resp. $(\Yh' , \oy')$. Combining the reduction $\mod p^n$ of (\ref{eq:32}) with (\ref{eq:33}), it follows that these automorphisms are equal. Hence $\rho_{E,n}$ is independent of all choices.

Now we check that the $\rho_{E,n}$ form a projective system. For a given $n \ge 1$ choose $\pi : \Yh \to \oeX$ with group $G$ in $\eT^{\good}$ such that $\pi^*_{n+1} E_{n+1}$ is a trivial bundle. Then $\pi^*_n E_n$ is trivial as well. For $\oy$ in $\Yh_{\oK} (\oK)$ over $\ox$ consider the commutative diagram
\[
\xymatrix{
\Yh_n \ar[d]_a & \spec \eo_n \ar[l]_{y_n} \ar[d]^b \ar[r]^{x_n} & \eX_n \ar[d]\\
\Yh_{n+1} & \spec \eo_{n+1} \ar[l]_{y_{n+1}} \ar[r]^{x_{n+1}} & \eX_{n+1} \; .
}
\]
It gives rise to the commutative diagram
\begin{equation}
  \label{eq:34}
  \vcenter{\xymatrix{
\Gamma (\Yh_n , \pi^*_n E_n) \ar[r]^{\overset{y^*_n}{\sim}} & \Gamma (\spec \eo_ n , x^*_n E_n) \ar@{=}[r] & E_{x_n} \\
\Gamma (\Yh_{n+1} , \pi^*_{n+1} E_{n+1}) \ar[u]^{a^*} \ar[r]^{\overset{y^*_{n+1}}{\sim}} & \Gamma (\spec \eo_{n+1} , x^*_{n+1} E_{n+1}) \ar[u]_{b^*} \ar@{=}[r] & E_{x_{n+1}} 
}}
\end{equation}
Here $b^*$ is identified with the natural reduction map from the $\eo_{n+1}$-module $E_{x_{n+1}}$ to the $\eo_n$-module $E_{x_n}$. Namely, we have:
\[
\Gamma (\spec \eo_n , x^*_n E_n) = \Gamma (\spec \eo_{n+1} , x^*_{n+1} E_{n+1}) \otimes_{\eo_{n+1}} \eo_n \quad \mbox{and} \quad b^* \ent \id \otimes_{\eo_{n+1}} \eo_n 
\]
since $x^*_n E_n$ and $x^*_{n+1} E_{n+1}$ are trivial bundles. Now it follows from (\ref{eq:34}) that 
\[
\Gamma (\Yh_n , \pi^*_n E_n) = \Gamma (\Yh_{n+1} , \pi^*_{n+1} E_{n+1}) \otimes_{\eo_n} \eo_{n+1} \quad \mbox{with} \quad a^* \ent \id \otimes_{\eo_{n+1}} \eo_n \; .
\]
In particular there are natural reduction maps
\[
\Aut_{\eo_{n+1}} \Gamma (\spec \eo_{n+1} , x^*_{n+1} E_{n+1}) \longrightarrow \Aut_{\eo_n} \Gamma (\spec \eo_n , x^*_n E_n)
\]
and
\[
\Aut_{\eo_{n+1}} \Gamma (\Yh_{n+1} , \pi^*_{n+1} E_{n+1}) \longrightarrow \Aut_{\eo_n} \Gamma (\Yh_n, \pi^*_n E_n) \; .
\]
It now suffices to show that the following diagram is commutative:
\[
\xymatrix{
 \ar@{}[ddr] |{\fbox{\tiny 1}} & \Aut_{\eo_{n+1}} \Gamma (\Yh_{n+1} , \pi^*_{n+1} E_{n+1}) \ar[dd]^{\red} \ar[r]^{\overset{\mathrm{via}\,y^*_{n+1}}{\sim}} & \Aut_{\eo_{n+1}} \Gamma (\spec \eo_{n+1} , x^*_{n+1} E_{n+1})  \ar@{}[ddl] |{\fbox{\tiny 2}} \ar[dd]^{\red} \\
G^{\op} \ar[ur] \ar[dr] \\
 & \Aut_{\eo_n} \Gamma (\Yh_n , \pi^*_n E_n) \ar[r]^{\overset{\mathrm{via}\,y^*_n}{\sim}} & \Aut_{\eo_n} \Gamma (\spec \eo_n , x^*_n E_n) \; .
}
\]
Commutativity of \fbox{\tiny 1} is clear. Diagram \fbox{\tiny 2} commutes because diagram (\ref{eq:34}) does. Thus theorem \ref{t16} is proved.
\end{proof}

With the discrete topology on $\Aut_{\eo_n} E_{x_n}$ the map $\rho_{E,n}$ is continuous since it factors over a finite quotient of $\pi_1 (\oX , \ox)$ resp. $\pi_1 (\oA , 0)$. It follows that
\[
\rho_E := \varprojlim_n \rho_{E,n} : \pi_1 (\oX , \ox) \longrightarrow \Aut_{\eo} E_{x_{\eo}}
\]
resp. 
\[
\rho_E := \varprojlim_n \rho_{E,n} : \pi_1 (\oA , 0) \to \Aut_{\eo} E_{x_{\eo}}
\] 
is a continuous homomorphism if $\Aut_{\eo} E_{\eo} \cong \GL_r (\eo)$ for $r = \rank E$ carries the topology induced from the one on $\eo$.

We now extend the map $E \mapsto \rho_E$ to a functor
\begin{equation}
  \label{eq:35}
  \rho : \eB_{\eX_{\eo}} \longrightarrow \Rep_{\pi_1 (\oX , \ox)} (\eo) \; \mbox{resp.} \; \rho : \eB_{\Ah_{\eo}} \longrightarrow \Rep_{\pi_1 (\oA , 0)} (\eo) \; .
\end{equation}
Again, we give the details only in the curve case.

Let $f : E \to E'$ be a morphism of vector bundles in $\eB_{\eX_{\eo}}$. Then we have an induced $\eo$-linear map
\[
\xymatrix{
E_{x_{\eo}} \ar[r]^-{\rho_f := f_{x_{\eo}}} \ar@{=}[d]& E'_{x_{\eo}} \ar@{=}[d]\\
\Gamma (\spec \eo , x^*_{\eo} E) \ar[r]^{x^*_{\eo} f} & \Gamma (\spec \eo , x^*_{\eo} E') \; .
}
\]
We have to check that $\rho_f$ is $\pi_1 (\oX , \ox)$-equivariant if $\pi_1 (\oX , \ox)$ acts on $E_{x_{\eo}}$ via $\rho_E$ and on $E'_{x_{\eo}}$ via $\rho_{E'}$. For this it suffices to show that for every $n \ge 1$ the natural $\eo_n$-linear map
\[
\rho_{ f,n} = f_{x_n} : E_{x_n} \longrightarrow E'_{x_n}
\]
is $\pi_1 (\oX , \ox)$-equivariant with respect to the actions $\rho_{E,n}$ and $\rho_{E' , n}$. Choose $\pi : \Yh \to \oeX$ in $\eT^{\good}$ such that both $\pi^*_n E_n$ and $\pi^*_n E'_n$ are trivial. Let $\oy$ and $G$ be as before. Consider the diagram defining the $\pi_1 (\oX , \ox)$-actions on $E_n$ and $E'_n$
\[
\xymatrix{
 & \Aut_{\eo_n} \Gamma (\Yh_n , \pi^*_n E_n) \ar[r]^-{\overset{\mathrm{via}\,y^*_n}{\sim}} & \Aut_{\eo_n} E_n \\
\pi_1 (\oX , \ox) \xrightarrow{\varphi_{\oy}} G^{\op} \ar[ur] \ar[dr] & \\
 & \Aut_{\eo_n} \Gamma (\Yh_n , \pi^*_n E'_n) \ar[r]^-{\overset{\mathrm{via}\,y^*_n}{\sim}} & \Aut_{\eo_n} E'_n \; .
}
\]
In the commutative diagram
\[
\xymatrix{
\Gamma (\Yh_n , \pi^*_n E_n) \ar[d]^{\Gamma (\Yh_n , \pi^*_n f_n)} \ar[r]^-{\overset{y^*_n}{\sim}} & \Gamma (\spec \eo_n , x^*_n E_n) \ar@{=}[r] \ar[d]^{\Gamma (\spec \eo_n , x^*_n f_n)} & E_{x_n} \ar[d]^{\rho_{f,n}} \\
\Gamma (\Yh_n , \pi^*_n E'_n) \ar[r]^-{y^*_n} & \Gamma (\spec \eo_n , x^*_n E'_n) \ar@{=}[r] & E'_{x_n}
}
\]
the map $\Gamma (\Yh_n , \pi^*_n f_n)$ is $G^{\op}$-equivariant. It follows that the map $\rho_{f,n}$ is $G^{\op}$- and hence $\pi_1 (\oX , \ox)$-equivariant as well.

It is clear that the $\rho_E$'s and $\rho_f$'s define a functor (\ref{eq:35}).

\begin{rem}
  The same constructions using $\eT^{\fin}_{\oeX}$ instead of $\eT_{\oeX}$ define a functor $\rho^{\fin}$ on $\eB^{\fin}_{\eX_{\eo}}$. It is clear that $\rho^{\fin}$ is the restriction of $\rho$ to $\eB_{\eX_{\eo}}$ via the full embedding $\eB^{\fin}_{\eX_{\eo}} \subset \eB_{\eX_{\eo}}$. 
\end{rem}

We will call a sequence of vector bundles in $\eB_{\eX_{\eo}}$ exact if it is exact in $\Vec_{\eX_{\eo}}$. Same for $\eB_{\Ah_{\eo}}$.

\begin{theorem}
  \label{t17}
The categories $\eB_{\eX_{\eo}}$ and $\eB_{\Ah_{\eo}}$ are additive and closed under tensor products, duals, internal homs and exterior powers of vector bundles. The functors
\[
\rho : \eB_{\eX_{\eo}} \longrightarrow \Rep_{\pi_1 (\oX , \ox)} (\eo) \quad \mbox{resp.} \quad \rho : \eB_{\Ah_{\eo}} \longrightarrow \Rep_{\pi_1 (\oA , 0) } (\eo)
\]
are additive and commute with the above operations on vector bundles. They are exact in the sense that they transform exact sequences into exact sequences.
\end{theorem}

\begin{proof}
  We first have to show that $\eB_{\eX_{\eo}}$ is closed under $\oplus , \otimes , ^* , \uHom$ and $\Lambda^i$. Assume that $E$ and $E'$ are in $\eB_{\eX_{\eo}}$ and fix $n \ge 1$. By theorem \ref{t13}, part {\bf b} there is an object $\pi : \Yh \to \oeX$ in $\eT^{\good}$ such that both $\pi^*_n E_n$ and $\pi^*_n E'_n$ are trivial bundles on $\Yh_n$. Hence the following bundles are trivial on $\Yh_n$ as well: \\
$\pi^*_n (E \oplus E')_n = \pi^*_n E_n \oplus \pi^*_n E'_n$ and $\pi^*_n (E \otimes_{\eX} E')_n = \pi^*_n E_n \otimes_{\Yh_n} \pi^*_n E'_n$, \\
$\pi^*_n (E^*)_n = \pi^*_n (E^*_n)$ and $\pi^*_n \uHom_{\eX_n} (E , E')_n = \uHom_{\Yh_n} (\pi^*_n E_n , \pi^*_n E'_n)$,\\
and $\pi^*_n (\Lambda^i E)_n = \Lambda^i \pi^*_n E_n$. 

Hence $E \oplus E' , E \otimes E' , E^* , \uHom (E,E')$ and $\Lambda^i E$ are in $\eB_{\eX_{\eo}}$ as well. As $y^*_n$ is an isomorphism on global sections the natural isomorphisms
\[
(E \otimes E')_{x_n} = E_{x_n} \otimes_{\eo_n} E'_{x_n} \quad \mbox{and} \quad \uHom (E , E')_{x_n} = \Hom_{\eo_n} (E_{x_n} , E'_{x_n}) \quad \mbox{etc.} 
\]
show that we have
\[
\Gamma (\Yh_n , \pi^*_n (E \otimes E')_n) = \Gamma (\Yh_n , \pi^*_n E_n) \otimes_{\eo_n} \Gamma (\Yh_n , \pi^*_n E'_n) \; \mbox{etc.} 
\]
For $\oplus$ this is clear. It follows that we have
\[
\rho_{E \oplus E',n} = \rho_{E,n} \oplus \rho_{E',n} \quad \mbox{and} \quad \rho_{E \otimes E',n} = \rho_{E,n} \otimes \rho_{E',n} \quad \mbox{etc.} 
\]
Exactness of $\rho$ is clear since the fibre functor $E \mapsto E_{x_{\eo}}$ from $\Vec_{\eX_{\eo}}$ into the category of $\eo$-modules is exact. Same arguments in the case of abelian varieties.
\end{proof}

We now discuss the effect of Galois conjugation on our constructions in the following case: We have $\oX = X \otimes_K \oK$ and $\oeX = \eX \otimes_{\eo_K} \eo_{\oK}$ for a smooth projective geometrically connected $X / K$ with a flat and proper model $\eX / \eo_K$. Moreover there is a $K$-rational point $x \in X (K)$ and $\ox \in \oX (\oK)$ is the induced point. For a scheme $\Yh$ over $\spec \eo_{\oK}$ and an element $\sigma$ of the absolute Galois group $G_K$ of $K$ one sets $^{\sigma} \Yh = \Yh \otimes_{\eo_{\oK} , \sigma} \eo_{\oK}$. This scheme is isomorphic over $\spec \eo_{\oK}$ to the spec $\eo_{\oK}$-scheme $(\Yh \to \spec \eo_{\oK} \xrightarrow{\spec (\sigma^{-1})} \spec \eo_{\oK})$. We have $^{\tau} (\,\!^{\sigma} \Yh) = \; ^{\tau \sigma} \Yh$ and there is a commutative diagram
\[
\xymatrix{
\Yh \ar[r]^{\overset{\sigma}{\sim}} \ar[d] & ^{\sigma} \Yh \ar[d]\\
\spec \eo_{\oK} \ar[r]^{\overset{\spec \sigma^{-1}}{\sim}} & \spec \eo_{\oK} \; .
}
\]
For $\oeX = \eX \otimes_{\eo_K} \eo_{\oK}$ the map $\id \times_{\spec \eo_K} \spec (\sigma^{-1})$ gives an isomorphism $^{\sigma} \oeX \to \oeX$ over $\spec \eo_{\oK}$ which will be used to identify $^{\sigma} \oeX$ with $\oeX$.

Similar remarks apply to schemes and morphisms over $\spec \eo , \spec \eo_n , \spec \oK , \linebreak
\spec \C_p$ and to vector bundles on such schemes.

In particular every $\sigma$ induces a functor from $\Vec_{\Yh_{\eo}}$ to $\Vec_{\sigma_{\Yh_{\eo}}}$ by sending a vector bundle $E$ over $\Yh_{\eo}$ to $^{\sigma}E$ over $^{\sigma} \Yh_{\eo}$ and a morphism $f$ to $^{\sigma} f$. This functor is an isomorphism of categories. For $\Yh = \oeX = \eX \otimes_{\eo_K} \eo_{\oK}$, one obtains an automorphism of $\Vec_{\eX_{\eo}}$ and by construction a homomorphism $G_K \to \Aut (\Vec_{\eX_{\eo}})$. We say that $G_K$ acts from the left on the category $\Vec_{\eX_{\eo}}$.  

There is also an action of $G_K$ on the category $\Rep_{\pi_1 (\oX , \ox)} (\eo)$. The representation $\rho : \pi_1 (\oX , \ox) \to \Aut_{\eo} \Gamma$ is sent to the representation
\[
^{\sigma}\rho = d_{\sigma} \verk \rho \verk a_{\sigma}^{-1} : \pi_1 (\oX , \ox) \longrightarrow \Aut_{\eo} \; ^{\sigma}\Gamma \; .
\]
Here $^{\sigma} \Gamma$ is $\Gamma$ with the twisted $\eo$-module structure $\lambda \cdot \gamma = \sigma^{-1} (\lambda) \gamma$ for $\lambda \in \eo$ and $\gamma \in \Gamma$ i.e. $^{\sigma} \Gamma = \Gamma \otimes_{\eo , \sigma} \eo$. We write the identity map $\Gamma \xrightarrow{\id} \;^{\sigma}\Gamma$ as $\sigma : \Gamma \to \;^{\sigma} \Gamma$ since it is $\sigma$-linear and define the conjugation isomorphism
\[
d_{\sigma} : \Aut_{\eo} \Gamma \longrightarrow \Aut_{\eo} \; ^{\sigma} \Gamma
\]
by setting $d_{\sigma} (f) = \sigma \verk f \verk \sigma^{-1}$. To define $a_{\sigma}$ consider the canonical exact sequence
\[
1 \longrightarrow \pi_1 (\oX , \ox) \longrightarrow \pi_1 (X , \ox) \longrightarrow G_K \longrightarrow 1 \; .
\]
The $K$-rational point $x$ provides a splitting
\[
x_* : G_K = \pi_1 (\spec K , \spec \oK) \longrightarrow \pi_1 (X , \ox)
\]
and hence an action of $G_K$ on $\pi_1 (\oX , \ox)$. Namely $\sigma$ acts by the automorphism
\[
a_{\sigma} : \pi_1 (\oX , \ox) \longrightarrow \pi_1 (\oX , \ox)
\]
defined by conjugation with $x_* (\sigma)$ i.e. $a_{\sigma} (\gamma) = x_* (\sigma) \gamma x_* (\sigma)^{-1}$. In case $X = A$ is an abelian variety and $x = 0$, the action of $a_{\sigma}$ on $\pi_1 (\oA , 0) = TA$ is given by applying $\sigma$ to the components in $TA = \varprojlim{} A_{p^n} (\oK)$. 

Thus we have defined $^{\sigma} \rho$. Now, consider a morphism $h : \rho_1 \to \rho_2$ in $\Rep_{\pi_1 (\oX , \ox)} (\eo)$, i.e. an $\eo$-linear map $h : \Gamma_1 \to \Gamma_2$ which is equivariant with respect to the $\pi_1 (\oX , \ox)$-actions via $\rho_1$ on $\Gamma_1$ and $\rho_2$ on $\Gamma_2$. We define $^{\sigma} h : \, ^{\sigma} \rho_1 \to \, ^{\sigma} \rho_2$ to be $h$ itself on the underlying abelian groups $^{\sigma} \Gamma_1 = \Gamma_1$ and $^{\sigma} \Gamma_2 = \Gamma_2$. In this way $\sigma$ induces an automorphism of $\Rep_{\pi_1 (\oX, \ox)} (\eo)$ and we get an action of $G$ on this category. 

Coming back to $c_{\sigma}$, assume that $\Gamma = \eo^n$. Then $^{\sigma} \Gamma = \eo^n$ as well, but with a different $\eo$-module structure. One has a commutative diagram
\[
\xymatrix{
\Aut_{\eo} (\eo^n) \ar@{=}[d] \ar[r]^{d_{\sigma}} & \Aut_{\eo} (\,\!^{\sigma} \eo^n)\ar@{=}[d]\\
\GL_n (\eo) \ar[r]^{\sigma} & \GL_n (\eo) 
}
\]
where the lower map sends a matrix $A = (a_{ij})$ to the matrix $\sigma (A) = (\sigma (a_{ij}))$. Namely, we have
\[
d_{\sigma} (f) (e_i) = (\sigma \verk f \verk \sigma^{-1}) (e_i) = (\sigma \verk f) (e_i) = \sigma \sum_j a_{ij} e_j = \sum_j \sigma (a_{ij}) e_j\; .
\]

\begin{theorem}
  \label{t18}
The above actions of $G_K$ on $\Vec_{\eX_{\eo}}$ and $\Vec_{\Ah_{\eo}}$ induce actions on $\eB_{\eX_{\eo}}$ and $\eB_{\Ah_{\eo}}$. The functors
\[
\rho : \eB_{\eX_{\eo}} \longrightarrow \Rep_{\pi_1 (\oX , \ox)} (\eo) \quad \mbox{and} \quad \rho : \eB_{\Ah_{\eo}} \longrightarrow \Rep_{\pi_1 (\oA , 0)} (\eo)
\]
commute with the $G_K$-actions on the categories.
\end{theorem}

\begin{proof}
  For a bundle $E$ over $\eX_{\eo}$ in $\eB_{\eX_{\eo}}$ and an integer $n \ge 1$ let $\pi : \Yh \to \oeX$ with group $G$ be an object of $\eT_{\oeX}$ such that $\pi^*_n E_n$ is a trivial bundle. Then $^{\sigma} \pi : \, ^{\sigma} \Yh \to \,^{\sigma} \oeX = \oeX$ is again an object of $\eT_{\oeX}$. The group is $G$ but with the conjugate action
\[
G \longrightarrow \Aut_{\oeX} \Yh \longrightarrow \Aut_{\oeX} \Yh^{\sigma} \; , \; g \longmapsto \sigma \verk g \verk \sigma^{-1}
\]
where $\sigma : \Yh \to \,^{\sigma} \Yh$ is the above isomorphism of schemes. Since $^{\sigma} \pi^*_n (\,\!^{\sigma} E_n) = \,^{\sigma} (\pi^*_n E_n)$ is again a trivial bundle, $^{\sigma} E$ lies in $\eB_{\eX_{\eo}}$ as well. Hence $G_K$ acts on $\eB_{\eX_{\eo}}$ and similarly on $\eB_{\Ah_{\eo}}$. We now check that the functor $\rho$ respects the Galois actions. Choose $\pi : \Yh \to \oeX$ in $\eT^{\good}$ and $\oy$ in $\Yh_{\oK} (\oK)$. The commutative diagram
\[
\xymatrix{
\Yh \ar[r]^{\sigma} \ar[d]_{\pi_n} & \,^{\sigma} \Yh \ar[d]^{^{\sigma}\pi_n} \\
\eX_n \ar[r]^{\sigma} & \eX_n \ar@{=}[r] & \,^{\sigma} \eX_n
}
\]
together with the relation $E_n = \sigma^* (\,\!^{\sigma} E_n)$, shows that $\sigma$ induces an isomorphism
\[
\sigma^* : \Gamma (\,\!^{\sigma} \Yh_n , \,^{\sigma} \pi^*_n (\,\!^{\sigma} E_n)) \silo \Gamma (\Yh_n , \pi^*_n E_n) \; .
\]
Consider the point $\oy^{\sigma} = \sigma \verk \oy \verk \sigma^{-1}$ in $^{\sigma} \Yh (\oK)$. It lies over $\ox^{\sigma} = \sigma \verk \ox \verk \sigma^{-1} = \ox$. For the latter equality note that $\ox$ is induced from the $K$-rational point $x \in \eX (K)$.

By definition of $a_{\sigma}$ the left square in the following diagram commutes, where $b_{\sigma} (g) = \sigma \verk g \verk \sigma^{-1}$ (composition of morphisms) and $c_{\sigma} (\beta) = (\sigma^{-1})^* \verk \beta \verk \sigma^*$
\[
\xymatrix{
\pi_1 (\oX , \ox) \ar[r]^{\varphi_{\oy}} \ar[d]_{a_{\sigma}}^{\wr} & \Aut^{\op}_{\oX} \Yh_{\oK} \ar[d]^{\wr \, b_{\sigma}} \ar@{=}[r] & G^{\op} \ar[r] & \Aut_{\eo_n} \Gamma (\Yh_n , \pi^*_n E_n) \ar[d]^{\wr \, c_{\sigma}} \ar[r]^-{\overset{\mathrm{via}\,y^*_n}{\sim}} & \Aut_{\eo_n} E_{x_n} \ar[d]^{\wr \, d_{\sigma}} \\
\pi_1 (\oX , \ox) \ar[r]^{\varphi_{\oy^{\sigma}}} & \Aut^{\op}_{\oX} \, ^{\sigma}\Yh_{\oK} \ar@{=}[r] & G^{\op} \ar[r] & \Aut_{\eo_n} \Gamma (\,\!^{\sigma} \Yh_n , \,^{\sigma } \pi^*_n (\,\!^{\sigma} E_n)) \ar[r]^-{\overset{\mathrm{via}\,(y^{\sigma}_n)^*}{\sim}} & \Aut_{\eo_n} \, ^{\sigma} E_{x_n} \; .
}
\]
Let $y^{\sigma}_n$ be the $\mod p^n$ reduction of the point $\oy^{\sigma}$. It is clear that the middle square commutes. As for the right hand square note that the commutative diagram
\[
\xymatrix{
\Yh_n \ar[d]_{\sigma} & \spec \eo_n\ar[l]_{y_n} \ar[d]^{\sigma = \spec (\sigma^{-1})} \\
^{\sigma} \Yh_n & \spec \eo_n \ar[l]^{y^{\sigma}_n}
}
\]
induces a commutative diagram:
\[
\xymatrix{
\Gamma (\Yh_n , \pi^*_n E_n) \ar[r]^{\overset{y^*_n}{\sim}} & \Gamma (\spec \eo_n , x^*_n E_n) \ar@{=}[r] & E_{x_n} & \\
\Gamma (\,\!^{\sigma} \Yh_n , \, ^{\sigma} \pi^*_n (\,\!^{\sigma} E_n)) \ar[u]^{\sigma^* \, \wr} \ar[r]^{\overset{(y^{\sigma}_n)^*}{\sim}} & \Gamma (\spec \eo_n , x^*_n (\,\!^{\sigma} E_n)) \ar[u]_{\wr \, \sigma^* = \sigma^{-1}} \ar@{=}[r] & \ar[u]_{\sigma^{-1}} (\,\!^{\sigma} E)_{x_n} \ar@{=}[r] & ^{\sigma} (E_{x_n}) \; .
}
\]
Here, the vertical maps are $\sigma^{-1}$-linear as maps of $\eo_n$-modules. It follows that $c_{\sigma}$ and $d_{\sigma}$ map $\eo_n$-linear automorphisms to $\eo_n$-linear automorphisms. We have shown that
\[
\rho_{^{\sigma}\!E} \verk a_{\sigma} = d_{\sigma} \verk \rho_E
\]
i.e. that $\rho_{^{\sigma}\!E} = \,^{\sigma}\!\rho_E$. It is clear that for a morphism of vector bundles $f : E_1 \to E_2$ in $\eB_{\eX_{\eo}}$ we have $^{\sigma}\!\rho_f = \rho_{^{\sigma}\!f}$. Similar arguments in the case of abelian schemes.
\end{proof}

We now discuss further functorialities which will become useful later on.

\begin{theorem}
  \label{t19}
Let $\oX , \oX'$ be smooth and proper curves over $\oK$ with finitely presented flat and proper models $\oeX , \oeX'$ over $\eo_{\oK}$. Let $\ox \in \oX (\oK)$ be a point. Consider also abelian varieties $A , A'$ over $K$ with good reduction. Let $\Ah , \Ah'$ be the corresponding abelian schemes over $\eo_K$. Let
\[
f : \oeX \longrightarrow \oeX' \quad \mbox{resp.} \quad g : \oAh \longrightarrow \oAh' \quad \mbox{resp.} \quad h : \oeX \longrightarrow \oAh
\]
be morphisms over $\eo_{\oK}$ such that $g (0) = 0$ and $h (\ox) = 0$. Then the pullback functors $f^* , g^*$ and $h^*$ from $\Vec_{\eX'_{\eo}}$ to $\Vec_{\eX_{\eo}}$ etc. restrict to functors \\
$f^* : \eB_{\eX'_{\eo}} \to \eB_{\eX_{\eo}}$ etc. These functors are additive, exact and respect tensor products, internal homs, duals and exterior powers. The induced diagrams of functors are commutative up to canonical isomorphisms. E.g.:
\begin{equation}
  \label{eq:36}
  \vcenter{\xymatrix{
\eB_{\eX'_{\eo}} \ar[r]^{f^*} \ar[d]_{\rho} & \eB_{\eX_{\eo}} \ar[d]^{\rho} \\
\Rep_{\pi_1 (\oX' , f(\ox))} (\eo) \ar[r]^F & \Rep_{\pi_1 (\oX , \ox)} (\eo)
}}
\end{equation}
where $F$ is the functor induced by the map $f_* : \pi_1 (\oX , \ox) \to \pi_1 (\oX' , f (\ox))$. In particular, for a bundle $E'$ in $\eB_{\eX'_{\eo}}$ we have
\begin{equation}
  \label{eq:37}
  \rho_{f^* E'} = \rho_{E'} \verk f_* \; ,
\end{equation}
and similarly in the other case.
\end{theorem}

\begin{proof}
  For $E$ in $\eB_{\eX'_{\eo}}$ and an integer $n \ge 1$ there exists an object \\
$\pi' : \Yh' \to \oeX'$ of $\eT^{\good}_{\oeX'}$ such that $\pi^{'*}_n E'_n $ is a trivial bundle. According to theorem \ref{t13} {\bf e} the morphism $\tilde{\pi} : f^{-1} \Yh' = \Yh' \times_{\oeX'} \oeX \to \oeX$ is an object of $\eT_{\oeX}$. Using the commutative diagram
  \begin{equation}
    \label{eq:38}
    \vcenter{ \xymatrix{
(f^{-1} \Yh')_n \ar@{=}[r] & \Yh'_n \times_{\eX'_n} \eX_n \ar[r] \ar[d]_{\tilde{\pi}_n} & \Yh'_n \ar[d]^{\pi'_n}\\
 & \eX_n \ar[r]^{f_n} & \eX'_n
}}
  \end{equation}
it follows that $\tilde{\pi}^*_n (f^* E')_n = \tilde{\pi}^*_n (f^*_n E'_n)$ is trivial. Hence $f^* E'$ lies in $\eB_{\eX_{\eo}}$. The same argument, but using theorem \ref{t13}, {\bf f} shows that $h^*$ maps $\eB_{\Ah_{\eo}}$ into $\eB_{\eX_{\eo}}$. Finally, for $E'$ in $\eB_{\Ah'_{\eo}}$ and $n \ge 1$ choose some $N \ge 1$ such that $\pi^{'*}_{N,n} E'_n$ is trivial, where $\pi'_N : \oAh' \to \oAh'$ is the $N$-multiplication morphism. The assumption $g (0) = 0$ implies that $g$ is a homomorphism of abelian schemes. Hence the diagram
\begin{equation}
  \label{eq:39}
  \vcenter{\xymatrix{
\Ah_n \ar[r]^{g_n} \ar[d]_{\pi_{N,n}} & \Ah'_n \ar[d]^{\pi'_{N,n}} \\
\Ah_n \ar[r]^{g_n} & \Ah'_n
}}
\end{equation}
commutes. Thus $\pi^*_{N,n} (g^* E')_n = \pi^*_{N,n} g^*_n E'_n$ is trivial and thus $g^* E'$ lies in $\eB_{\Ah_{\eo}}$. It is clear that the functors $f^* , g^*$ and $h^*$ are exact, additive etc. The composition of $\rho \verk f^*$ and $F \verk \rho$ in diagram (\ref{eq:36}) with the forgetful functor
\[
\Rep_{\pi_1 (\oX , \ox)} (\eo) \longrightarrow \Mod_{\eo}
\]
agrees in both cases with the restriction of the fibre functor
\[
\Vec_{\eX'_{\eo}} \longrightarrow \Mod_{\eo} \; , \; E' \longmapsto E'_{f (x_{\eo})} \; , \; \varphi' \longmapsto \varphi'_{f (x_{\eo})}
\]
to the full subcategory $\eB_{\eX'_{\eo}}$ of $\Vec_{\eX'_{\eo}}$. 

For the commutativity of (\ref{eq:36}) it therefore remains to show this: The two $\pi_1 (\oX , \ox)$-module structures defined by $\rho \verk f^*$ resp. $F \verk \rho$ on $E'_{f (x_{\eo})}$ agree with each other. In other words, we have to show equation (\ref{eq:37}). For this, choose $\pi : \Yh \to \oeX$ in $\eT^{\good}_{\oeX}$ dominating $\tilde{\pi} : f^{-1} \Yh' \to \oeX$. Let $G$ be the group for $\pi$ and $G'$ the group for $\pi'$ and hence for $\tilde{\pi}$. Define the morphism $\psi : \Yh \to \Yh'$ as the composition $\psi = \beta \verk \alpha$ of the maps $\alpha$ and $\beta$ in the following diagram:
\[
\xymatrix{
\Yh \ar[r]^-{\alpha} \ar[d]_{\pi} & f^{-1} \Yh' \ar[r]^-{\beta} \ar[d]^{\tilde{\pi}} & \Yh' \ar[d]^{\pi'} \\
\oeX \ar@{=}[r] & \oeX \ar[r]^f & \oeX' \; .
}
\]
Let $\gamma : G \to G'$ be the map corresponding to the morphism $\pi \to \tilde{\pi}$ in $\eT_{\oeX}$. Then $\alpha$ and hence $\psi$ is $G$-equivariant if $G$ acts via $\gamma$ on $f^{-1} \Yh'$ resp. $\Yh'$. Since $\pi$ and $\pi'$ were chosen in $\eT^{\good}$ we have $G = \Aut_{\oX} \Yh_{\oK}$ and $G' = \Aut_{\oX} \Yh'_{\oK}$. The choice of a point $\oy \in \Yh_{\oK} (\oK)$ over $\ox$ determines the homomorphisms $\varphi_{\oy}$ and $\varphi_{\oy'}$ where $\oy' = \psi (\oy)$ in the following diagram:
\begin{equation}
  \label{eq:40}
  \vcenter{\xymatrix{
\pi_1 (\oX , \ox) \ar[r]^{\varphi_{\oy}} \ar[d]_{f_*} & \Aut^{\op}_{\oX} (\Yh_{\oK}) \ar@{=}[r] & G^{\op} \ar[d]^{\gamma} \\
\pi_1 (\oX' , f(\ox)) \ar[r]^{\varphi_{\oy'}} & \Aut^{\op}_{\oX'} (\Yh'_{\oK}) \ar@{=}[r] & G'^{\op} \; .
}}
\end{equation}
We claim that this diagram commutes: Namely by $G$-equivariance of $\psi$ and the definition of $f_*$ we have:
\begin{eqnarray*}
  ((\gamma \verk \varphi_{\oy}) (\sigma)) (\oy') & = & \gamma (\varphi_{\oy} (\sigma)) (\psi (\oy)) = \psi (\varphi_{\oy} (\sigma) (\oy)) = (\psi \verk \varphi_{\oy} (\sigma)) (\oy) \\
 & = & \varphi_{\oy'} (f_* (\sigma)) (\oy') = ((\varphi_{\oy'} \verk f_*) (\sigma)) (\oy') \; .
\end{eqnarray*}
Hence $(\gamma \verk \varphi_{\oy}) (\sigma) = (\varphi_{\oy'} \verk f_*) (\sigma)$ for every $\sigma$ in $\pi_1 (\oX , \ox)$. 

Next, note that the commutative diagram
\[
\xymatrix{
\Yh_n \ar[dd]_{\psi_n} \\
 & \spec\eo_n \ar[ul]_{y_n} \ar[dl]^{y'_n} \\
\Yh'_n
}
\]
induces a commutative diagram of isomorphisms:
\begin{equation}
  \label{eq:41}
  \vcenter{\xymatrix{
\Gamma (\Yh_n , \psi^*_n (\pi^{'*}_n E'_n)) \ar@{=}[r] & \Gamma (\Yh_n , \pi^*_n (f^*_n E'_n)) \ar[d]^{\wr \, y^*_n} \\
\Gamma (\Yh'_n , \pi^{'*}_n E'_n) \ar[u]^{\psi^*_n} \ar[r]^{\overset{y^{'*}_n}{\sim}} & E'_{f(x_n)} \; .
}}
\end{equation}
Note here that $\pi$ and $\pi'$ were both chosen in $\eT^{\good}$. We thus get a diagram
\begin{equation}
  \label{eq:42}
  \vcenter{\xymatrix{
G^{\op} \ar[r] \ar[d]_{\gamma} & \Aut_{\eo_n} \Gamma (\Yh_n , \pi^*_n (f^*_n E'_n)) \ar[d]^{\mathrm{via}\,\psi^*_n} \ar[r]^-{\overset{\mathrm{via}\,y^*_n}{\sim}} & \Aut_{\eo_n} (f^* E')_{x_n} \ar@{=}[d] \\
G'^{\op} \ar[r] & \Aut_{\eo_n} \Gamma (\Yh'_n , \pi^{'*}_n E'_n) \ar[r]^-{\mathrm{via}\,y^{'*}_n} & \Aut_{\eo_n} (E'_{f (x_n)}) \; .
}}
\end{equation}
The right square commutes because diagram (\ref{eq:41}) commutes. As for the left square, note the following equalities for $\sigma \in G$:
\begin{eqnarray*}
  (\mathrm{via}\,\psi^*_n) (\sigma^*_n) & := & (\psi^*_n)^{-1} \verk \sigma^*_n \verk \psi^*_n = (\psi^*_n)^{-1} \verk (\psi_n \verk \sigma_n)^* \overset{(1)}{=} (\psi^*_n)^{-1} \verk (\gamma (\sigma)_n \verk \psi_n)^* \\
& = & (\psi^*_n)^{-1} \verk \psi^*_n \verk \gamma (\sigma)^*_n = \gamma (\sigma)^*_n \; .
\end{eqnarray*}
Here (1) comes from the $G$-equivariance $\psi \verk \sigma = \gamma (\sigma) \verk \psi$ of $\psi$. Hence the outer square of (\ref{eq:42}) commutes as well and combining it with (\ref{eq:40}) we get the commutative diagram: 
\[
\xymatrix{
\pi_1 (\oX, \ox) \ar[r]^-{\rho_{f^* E',n}} \ar[d]_{f_*} & \Aut_{\eo_n} (f^* E')_{x_n} \ar@{=}[d] \\
\pi_1 (\oX' , f (\ox)) \ar[r]^{\rho_{E',n}} & \Aut_{\eo_n} (E'_{f(x_n)}) \; .
}
\]
Passing to the limit over $n$'s we obtain equation (\ref{eq:37}) as desired. Similar arguments show that the diagrams corresponding to (\ref{eq:36}) for $g^*$ and $h^*$ commute as well.
\end{proof}

\begin{theorem}
  \label{t20}
1) Let $\Ah / \eo_K$ be an abelian scheme. Then the full subcategory $\eB_{\Ah_{\eo}}$ of $\Vec_{\Ah_{\eo}}$ is closed under extensions and contains all line bundles on $\Ah_{\eo}$ which lie in $\Pic^0$.\\
2) Let $\oX$ be a smooth projective curve over $\oK$ with a  smooth and proper model $\oeX$. Then every extension of a trivial vector bundle on $\eX_{\eo}$ by a trivial vector bundle is an object of $\eB_{\eX_{\eo}}$. The category $\eB_{\eX_{\eo}}$ contains all line bundles on $\eX_{\eo}$ of degree zero relative to $\Spec \eo$.
\end{theorem}

\begin{proof}
  1) Let $0 \to E' \to E \to E'' \to 0$ be an exact sequence of vector bundles on $\Ah_{\eo}$ with $E'$ and $E''$ in $\eB_{\Ah_{\eo}}$. Fix $n \ge 1$. Then there are integers $N' , N'' \ge 1$ such that $\pi^*_{N'} E'_n$ and $\pi^*_{N''} E''_n$ are trivial bundles on $\Ah_n$. Setting $N = N' N''$ we get an exact sequence
\[
0 \longrightarrow \Oh^{r'} \longrightarrow \pi^*_N \Oh (E_n) \longrightarrow \Oh^{r''} \longrightarrow 0
\]
on $\Ah_n$ where $r' = \rank E' ,
r'' = \rank E''$. We claim that $\pi^*_{p^n} (\pi^*_N E_n) = \pi^*_{p^n N} E_n$ is a trivial bundle. For this it suffices to show that $\pi^*_{p^n}$ induces the zero map on
\[
\Ext^1_{\Ah_n} (\Oh , \Oh) = H^1 (\Ah_n , \Oh) \; .
\]
Consider the following diagram
\begin{equation}
  \label{eq:43}
  \vcenter{\xymatrix{
H^1 (\Ah_n , \Oh) \ar[r]^{\pi^*_M} \ar@{=}[d] & H^1 (\Ah_n , \Oh) \ar@{=}[d]\\
\Lie \Pic^0_{\Ah_n / \eo_n} \ar[r]^{\Lie \pi^*_M} \ar@{=}[d] & \Lie \Pic^0_{\Ah_n / \eo_n} \ar@{=}[d] \\
\Lie \hAh_n \ar[r]^{M \cdot \id} & \Lie \hAh_n \; .
}}
\end{equation}
The top square is commutative by \cite{BLR}, 8.4, Theorem 1. The bottom square commutes since it is known that the diagram
\[
\xymatrix{
\Pic^0_{\Ah_n / \eo_n} \ar[r]^{\pi^*_M} \ar@{=}[d] & \Pic^0_{\Ah_n / \eo_n} \ar@{=}[d] \\
\hAh_n \ar[r]^{\pi_M} & \hAh_n
}
\]
commutes and since $\Lie \pi_M = M \cdot \id$. Choosing $M = p^n$ in diagram (\ref{eq:43}) and noting that $\Lie \hAh_n$ is an $\eo_n$-module it follows that $\pi^*_{p^n} = 0$ on $H^1 (\Ah_n , \Oh)$. \\
Let $L$ be a line bundle on $\Ah_{\eo}$  in $\Pic^0$. Its reduction defines a class in $\Pic^0_{\Ah / \eo_K} (\eo_n) = \Pic^0_{\Ah_n / \eo_n} (\eo_n) $.\\
As noted before this is a torsion group and hence $L^{\otimes N}_n$ is trivial for some $N \geq 1$. On the other hand we have $\pi^*_N L_n \cong L^{\otimes N}_n$. Thus $L$ is an object of $\eB_{\Ah_\eo}$. 

2) There exists a finite extension $K / \Q_p$ and a smooth projective curve $X / K$ with a smooth and proper model $\eX / \eo_K$ such that $\oX = X \otimes_K \oK$ and $\oeX = \eX \otimes_{\eo_K} \eo_{\oK}$. We may assume that $X$ has a rational point $x \in X (K) = \eX (\eo_K)$. Set $\Bh = \Pic^0_{\eX / \eo_K}$ and $\Ah = \hat{\Bh}$. Consider the Albanese map
\[
h : \eX \longrightarrow \Ah \quad \mbox{with} \; h (x) = 0 \; .
\]
It induces an isomorphism
\[
h^* : \Pic^0_{\Ah / \eo_K} \silo \Pic^0_{\eX / \eo_K}
\]
realizing the biduality $\hAh = \hat{\hat{\Bh}} = \Bh$. Over fields these facts are well known. They extend over $\spec \eo_K$ by the universal mapping property of N\'eron models and by noting that an abelian scheme over $\spec \eo_K$ is the N\'eron model of its generic fibre. The base extension
\[
h^* : \Pic^0_{\Ah_{\eo} / \eo} \silo \Pic^0_{\eX_{\eo} / \eo}
\]
is an isomorphism as well. According to theorem \ref{t19} we have $h^* (\ob \eB_{\Ah_{\eo}}) \subset \ob \eB_{\eX_{\eo}}$. Since all line bundles on  $\Ah_{\eo}$ which are algebraically equivalent to zero lie in $\eB_{\Ah_{\eo}}$, it follows that all line bundles of degree zero on  $\eX_{\eo}$ lie in $\eB_{\eX_{\eo}}$. Next, the commutative diagram
\[
\xymatrix{
\Lie \Pic^0_{\Ah_{\eo} / \eo} \ar[r]^{\overset{\Lie h^*}{\sim}} \ar@{=}[d] & \Lie \Pic^0_{\eX_{\eo} / \eo} \ar@{=}[d] \\
H^1 (\Ah_{\eo} , \Oh) \ar[r]^{h^*} & H^1 (\eX_{\eo} , \Oh)
}
\]
shows that the map $h^*$ is an isomorphism on cohomology as well. It follows that $h^*$ induces an isomorphism $\Ext^1_{\Ah_\eo}(\Oh^r, \Oh^s) \simeq \Ext^1_{\eX_\eo}(\Oh^r,\Oh^s)$ for all $r,s \geq 1$. In particular, every bundle $E$ on $\eX_{\eo}$ which is an extension of a trivial vector bundle by a trivial vector bundle is isomorphic to a bundle of the form $h^* E'$ for some bundle $E'$ of the same type on $\Ah_{\eo}$. According to part 1), we know that $E'$ lies in $\eB_{\Ah_{\eo}}$. Because of the inclusion $h^* (\ob \eB_{\Ah_{\eo}}) \subset \ob \eB_{\eX_{\eo}}$ the bundle $E \simeq h^* E' $ lies in $\eB_{\eX_{\eo}}$ as was to be shown.
\end{proof}

\begin{rem}
  We hope that $\eB_{\eX_{\eo}}$ is closed under arbitrary extensions. Theorem \ref{t13} {\bf b} for $\eT^{\good}_{ss}$ may be useful for proving this.
\end{rem}

We now define a category of vector bundles on $X_{\C_p}$ which is somewhat analogous to the category of flat bundles on a Riemann surface. 

Let $\oX$ be a smooth projective variety over $\oK$ and set $X_{\C_p} = \oX \otimes \C_p$. For every finitely presented flat and proper model $\oeX$ of $\oX$ over $\spec \eo_{\oK}$, we have the restriction functor
\[
j^* = j^*_{\eX_{\eo}} : \Vec_{\eX_{\eo}} \longrightarrow \Vec_{X_{\C_p}}
\]
where $j : X_{\C_p} \hookrightarrow \eX_{\eo}$ is the open immersion of the generic fibre into $\eX_\eo = \oeX \otimes \eo$.

For an additive category $\mathfrak{D}$ let $\mathfrak{D} \otimes \Q$ be the $\Q$-linear category with the same objects as $\mathfrak{D}$ and $\Hom_{\mathfrak{D} \otimes \Q} ( D_1,D_2) = \Hom_{\mathfrak{D}}(D_1,D_2) \otimes \Q$. 

\begin{prop}
  \label{t21}
The induced functor
\[
j^* : \Vec_{\eX_{\eo}} \otimes \Q \longrightarrow \Vec_{X_{\C_p}}
\]
is fully faithful.
\end{prop}

\begin{proof}
  For bundles $E_1$ and $E_2$ on $\eX_{\eo}$ set $F = \uHom (E_1 , E_2)$. Then the claim follows from the commutative diagram
\[
\xymatrix{
\Hom_{\eX_{\eo}} (E_1 , E_2) \otimes \Q \ar[r]^{j^*} \ar@{=}[d] & \Hom_{X_{\C_p}} (j^* E_1 , j^* E_2)  \ar@{=}[d] \\
H^0 (\eX_{\eo} , F) \otimes \Q \ar[r]^{j^*} & H^0 (X_{\C_p} , j^* F)
}
\]
and flat base change \cite{EGAIII} 1.4.15 if we note that
\[
H^0 (\eX_{\eo} , F) \otimes \Q = H^0 (\eX_{\eo} , F) \otimes_{\eo} \C_p \; .
\]
\end{proof}

By the proposition the induced functors
\[
j^*_{\eX_{\eo}} : \eB_{\eX_{\eo}} \otimes \Q \longrightarrow \Vec_{X_{\C_p}} \quad \mbox{for a curve} \; \oX / \oK 
\]
and
\[
j^*_{\Ah_{\eo}} : \eB_{\Ah_{\eo}} \otimes \Q \longrightarrow \Vec_{A_{\C_p}} \quad \mbox{for an abelian variety} \; A / K 
\]
are fully faithful as well. Here $\Ah / \eo_K$ is the abelian scheme with $A = \Ah \otimes K$.

\begin{defn}
  \label{t22}
1) Let $\oX$ be a smooth projective curve over $\oK$. The full subcategory $\eB_{X_{\C_p}}$ of $\Vec_{X_{\C_p}}$ is defined as follows: $\ob \eB_{X_{\C_p}}$ consists of all vector bundles on $X_{\C_p}$ which are isomorphic to a bundle of the form $j^*_{\eX_{\eo}} E$ for some $\oeX$ as above and some $E \in \eB_{\eX_{\eo}}$.\\
2) Let $A / K$ and $\Ah / \eo_K$ be as above. The category $\eB_{A_{\C_p}}$ is defined as the full subcategory consisting of bundles on $A_{\C_p}$ which are isomorphic to a bundle of the form $j^*_{\Ah_{\eo}} E$ for some $E$ in $\eB_{\Ah_{\eo}}$. In other words, $\eB_{A_{\C_p}}$ is the essential image of the fully faithful functor $j^*_{\Ah_{\eo}} : \eB_{\Ah_{\eo}} \otimes \Q \to \Vec_{A_{\C_p}}$. The category $\eB_{A_{\C_p}}$ depends only on $\oA = A \otimes_K \oK$ and is equivalent via $j^*_{\Ah_{\eo}}$ to $\eB_{\Ah_{\eo}} \otimes \Q$. 
\end{defn}

Using the categories $\eB^{\fin}_{\eX_{\eo}}$ one defines similarly the category $\eB^{\fin}_{X_{\C_p}}$. It is a full subcategory of $\eB_{X_{\C_p}}$. 
In order to show that the category $\eB_{X_{\C_p}}$ has reasonable properties we need the following geometric fact.

\begin{prop}
  \label{t23}
Let $X$ be a smooth projective curve over a finite extension $K$ of $\Q_p$. Given proper, flat models $\eX_1 , \eX_2$ of $X$ over $\spec \eo_K$ there exists a diagram
\[
\eX_1 \xleftarrow{p_1} \eX_3 \xrightarrow{p_2} \eX_2
\]
where $\eX_3$ is a proper flat model of $X$ over $\spec \eo_K$ which can be assumed to be integral and normal. Moreover $p_1$ and $p_2$ restrict to the identity on the generic fibre $\eX_{1K} = \eX_{3K} = \eX_{2K} = X$.
\end{prop}

\begin{proof}
 Let $\eX'_3$ be the closure of the image of the morphism
\[
X \xrightarrow{\Delta} X \times_{\spec K} X \longrightarrow \eX_1 \times_{\spec \eo_K} \eX_2
\]
with the reduced subscheme structure. Then $\eX'_3$ is integral and proper over $\spec \eo_K$. It is flat since its generic point is mapped to the one of $\spec \eo_K$. Morphisms $p'_1 : \eX'_3 \to \eX_1$ and $p'_2 : \eX'_3 \to \eX_2$ restricting to the identity on the generic fibres are obtained by composing the closed immersions $\eX'_3 \hookrightarrow \eX_1 \times_{\spec \eo_K} \eX_2$ with the projection maps to $\eX_1$ and $\eX_2$. Passing to the normalization of $\eX_3$ in its function field the proposition follows.
\end{proof}

\begin{cor}
  \label{t24}
Let $\oX$ be a smooth, projective curve over $\oK$ with finitely presented flat and proper models $\oeX_1$ and $\oeX_2$ over $\spec \eo_{\oK}$. Then there is another such model $\oeX_3$ together with morphisms 
\[
\oeX_1 \xleftarrow{p_1} \oeX_3 \xrightarrow{p_2} \oeX_2
\]
restricting to the identity on the generic fibres. We have a commutative diagram of fully faithful functors:
\begin{equation}
\label{eq:44}
\vcenter{\xymatrix{
\eB_{\eX_{1\eo}} \otimes \Q \ar[rd]^{p^*_1} \ar@/^/[rrd]^{j^*_{\eX_{1\eo}}} & & \\
 & \eB_{\eX_{3\eo}} \otimes \Q \ar[r]^{j^*_{\eX_{3\eo}}} & \eB_{X_{\C_p}} \; .\\
\eB_{\eX_{2\eo}} \otimes \Q \ar[ru]_{p^*_2} \ar@/_/[rru]_{j^*_{\eX_{2\eo}}} & &
}}
\end{equation}
\end{cor}

\begin{proof}
  As the notation suggests, we can find $K / \Q_p , X , \eX_1 , \eX_2$ as in proposition \ref{t23} inducing $\oX , \oeX_1 , \oeX_2$ by base extension. Using theorem \ref{t19} and propositions \ref{t21} and \ref{t23}, the corollary follows. 
\end{proof}

Consider the category of all finitely presented flat and proper models $\oeX$ of $\oX$ over $\spec \eo_{\oK}$. Morphisms are over $\spec \eo_{\oK}$ and are supposed to induce the identity on $\oeX_{\oK} = \oX$. 

It follows from diagram (\ref{eq:36}) that for a morphism $f : \oeX \to \oeX'$ in the category of models of $\oX$, after fixing a base point $\ox \in \oX$ we have a commutative diagram
\begin{equation}
  \label{eq:45}
  \vcenter{
\xymatrix{
\eB_{\eX'_{\eo}} \ar[rr]^{f^*} \ar[dr]_{\rho} & & \eB_{\eX_{\eo}} \ar[dl]^{\rho}\\
 & \Rep_{\pi_1 (\oX , \ox)} (\eo) \; .
}}
\end{equation}

Let $\Rep_{\pi_1 (\oX , \ox)} (\C_p)$ be the category of continuous representations of $\pi_1 (\oX , \ox)$ on finite dimensional $\C_p$-vector spaces. Let $x_{\C_p} : \spec \C_p \to X_{\C_p}$ be the point of $X_{\C_p}$ induced by $\ox \in \oX (\oQ_p)$. Taking the fibres in $x_{\C_p}$ of vector bundles and maps between them defines a functor
\[
x^*_{\C_p} : \eB_{X_{\C_p}} \longrightarrow \Vec_{\C_p} \; .
\]
We now define a functor
\begin{equation}
  \label{eq:46}
  \rho : \eB_{X_{\C_p}} \longrightarrow \Rep_{\pi_1 (\oX, \ox)} (\C_p)
\end{equation}
which on the underlying vector spaces is $x^*_{\C_p}$. For any object $F$ in $\eB_{X_{\C_p}}$ choose a model $\oeX$ of $\oX$ as above, an object $E$ in $\eB_{\eX_{\eo}}$ and an isomorphism $\psi : F \silo j^* E$ in $\Vec_{X_{\C_p}}$. Let
\[
\psi_{x_{\C_p}} = x^*_{\C_p} (\psi) : F_{x_{\C_p}} \silo (j^* E)_{x_{\C_p}} = E_{x_{\eo}} \otimes_{\eo} \C_p
\]
be the induced isomorphism on fibres. We get a $\pi_1 (\oX , \ox)$-action on $F_{x_{\C_p}}$ by transporting the one on $E_{x_{\eo}} = \rho_E \in \Rep_{\pi_1 (\oX , \ox)} (\eo)$ via $\psi_{x_{\C_p}}$ to $F_{x_{\C_p}}$. The resulting action is independent of all choices. This is a formal consequence of corollary \ref{t24}, the commutative diagram (\ref{eq:36}) and the fact that $\rho$ on $\eB_{\eX_{\eo}}$ is a functor.

Correspondingly we have a functor
\begin{equation}
  \label{eq:47}
  \rho : \eB_{A_{\C_p}} \longrightarrow \Rep_{\pi_1 (\oA , 0)} (\C_p) \; .
\end{equation}
Note that by construction the objects in the image of $\rho$ in (\ref{eq:46}) and (\ref{eq:47}) all lie in the essential image of the natural functor
\begin{equation}
  \label{eq:48}
  \Rep_{\pi_1 (\oX , \ox)} (\eo) \otimes \Q \longrightarrow \Rep_{\pi_1 (\oX , \ox)} (\C_p) \; .
\end{equation}
This gives no further information though because of the following fact:

\begin{prop}
  \label{t25}
The natural functor (\ref{eq:48}) is an equivalence of categories.
\end{prop}

\begin{proof}
  Set $\pi = \pi_1 (\oX , \ox)$ and take $\Gamma_1 , \Gamma_2$ in $\Rep_{\pi_1 (\oX , \ox)} (\eo)$. For full faithfulness one has to check that the natural map
\[
\Hom_{\pi} (\Gamma_1 , \Gamma_2) \otimes \Q \longrightarrow \Hom_{\pi} (\Gamma_1 \otimes \Q , \Gamma_2 \otimes \Q)
\]
is an isomorphism. This will follow if the $\pi$-equivariant map
\[
\Hom_{\eo} (\Gamma_1 , \Gamma_2) \otimes \Q \longrightarrow \Hom_{\C_p} (\Gamma_1 \otimes \Q , \Gamma_2 \otimes \Q)
\]
is an isomorphism. But this has to be noted for $\Gamma_1 = \Gamma_2 = \eo$ only where it is clear.

For essential surjectivity, the first part of the argument is standard: Let $V$ be in $\Rep_{\pi} (\C_p)$. We have to show that $V$ is isomorphic to an object of the form $\Gamma \otimes \Q = \Gamma \otimes_{\eo} \C_p$ with $\Gamma$ in $\Rep_{\pi} (\eo)$. Fix a basis $e_1 , \ldots , e_r$ of $V$ and set $\tilde{\Gamma} = \bigoplus^r_{i=1} \eo \cdot e_i$. Then $\tilde{\Gamma} \otimes_{\eo} \C_p = V$. The stabilizer $\pi_{\tilde{\Gamma}}$ of $\tilde{\Gamma}$ in $\pi$ is open since
\[
\pi_{\tilde{\Gamma}} = \{ g \in \pi \tei g (e_i) \in \tilde{\Gamma} \; \mbox{and} \; g^{-1} (e_i) \in \tilde{\Gamma} \}
\]
and since the operation of $\pi$ on $V$ is continuous. Hence $\pi / \pi_{\tilde{\Gamma}}$ is finite. Let $\Gamma$ be the $\eo$-submodule of $V$ which is generated by the finitely many translates $h \cdot \tilde{\Gamma}$ for $h \pi_{\tilde{\Gamma}}$ in $\pi / \pi_{\tilde{\Gamma}}$. Then $\Gamma$ is $\pi$-invariant and $\Gamma \otimes_{\eo} \C_p = V$.

Clearly, $\Gamma$ is a finitely generated torsion-free $\eo$-module. Since $\eo$ is a Bezout domain it follows that $\Gamma$ is a free $\eo$-module (whose rank is necessarily equal to $r$), cf. \cite{B} Lemma 3.9. It follows that $\Gamma$ lies in $\Rep_{\pi} (\eo)$. 
\end{proof}

We can now formulate the main properties of our categories $\eB_{X_{\C_p}}$ and the functors $\rho$.

\begin{theorem}
  \label{t26}
Let $\oX , \oX'$ be smooth and proper curves over $\oQ_p$ and let $\ox$ in $\oX (\oQ_p)$ be a point. Consider also abelian varieties $A , A'$ over $K$ with good reduction. Finally, let
\[
f : \oX \to \oX' \quad \mbox{resp.} \quad g : \oA \to \oA' \quad \mbox{resp.} \quad h : \oX \to \oA
\]
be morphisms over $\oK$ such that $g (0) = 0$ and $h (\ox) = 0$. The following assertions hold:

{\bf a} The categories $\eB_{X_{\C_p}}$ and $\eB_{A_{\C_p}}$ are full, additive subcategories of $\Vec_{X_{\C_p}}$ resp. $\Vec_{A_{\C_p}}$. They are closed under tensor products, duals, internal homs and exterior powers. 

{\bf b} All vector bundles in $\eB_{X_{\C_p}}$ have degree zero. For all vector bundles in $\eB_{A_{\C_p}}$ the determinant line bundle is algebraically equivalent to zero.

{\bf c} The category $\eB_{A_{\C_p}}$ contains all line bundles algebraically equivalent to zero on $A_{\C_p}$ and is closed under extensions in $\Vec_{A_{\C_p}}$.\\
If $\oX$ has a  smooth and proper model over $\spec \oZ_p$, then $\eB_{X_{\C_p}}$ contains all line bundles of degree zero on $X_{\C_p}$. Moreover, all bundles on $X_{\C_p}$ which are extensions of a trivial vector bundle by a trivial vector bundle then belong to $\eB_{X_{\C_p}}$.

{\bf d} Pullback of vector bundles induces exact additive functors
\[
f^* : \eB_{X'_{\C_p}} \to \eB_{X_{\C_p}} \quad \mbox{resp.} \quad g^* : \eB_{A'_{\C_p}} \to \eB_{A_{\C_p}} \quad \mbox{resp.} \quad h^* : \eB_{A_{\C_p}} \to \eB_{X_{\C_p}}
\]
where for $h^*$ we require that $\oX$ has a smooth model $\oeX$ over $\eo_{\oK}$. These functors commute with tensor products, duals, internal homs and exterior powers (up to canonical isomorphisms.)

{\bf e} The functors described in (\ref{eq:46}) and (\ref{eq:47})
\[
\rho : \eB_{X_{\C_p}} \longrightarrow \Rep_{\pi_1 (\oX, \ox)} (\C_p) \quad \mbox{and} \quad \rho : \eB_{A_{\C_p}} \longrightarrow \Rep_{\pi_1 (\oA , 0)} (\C_p)
\]
are additive, exact and commute with tensor products, duals, internal homs and exterior powers.

{\bf f} The following diagrams are commutative:
\[
\begin{array}{c}
\def\objectstyle{\scriptstyle}
\def\labelstyle{\scriptstyle}
\xymatrix{\eB_{X'_{\C_p}} \ar[r]^{f^*} \ar[d]_{\rho} & \eB_{X_{\C_p}} \ar[d]^{\rho} \\
\Rep_{\pi_1 (\oX' , f (\ox))} (\C_p) \ar[r]^{F} & \Rep_{\pi_1 (\oX , \ox)} (\C_p)
} \quad 
\xymatrix{\eB_{A'_{\C_p}} \ar[r]^{g^*} \ar[d]_{\rho} & \eB_{A_{\C_p}} \ar[d]^{\rho} \\
\Rep_{\pi_1 (\oA' , 0)} (\C_p) \ar[r]^{G} & \Rep_{\pi_1 (\oA, 0)} (\C_p)
} \\ 
\def\objectstyle{\scriptstyle}
\def\labelstyle{\scriptstyle}
\xymatrix{\eB_{A_{\C_p}} \ar[r]^{h^*} \ar[d]_{\rho} & \eB_{X_{\C_p}} \ar[d]^{\rho} \\
\Rep_{\pi_1 (\oA , 0)} (\C_p) \ar[r]^{H} & \Rep_{\pi_1 (\oX , \ox)} (\C_p)
}
\end{array}
\]
Here $F,G,H$ are the functors induced by composition with
\[
\scriptstyle
f_* : \pi_1 (\oX , \ox) \longrightarrow \pi_1 (\oX' , f(\ox)) \; , \; g_* : \pi_1 (\oA , 0) \longrightarrow \pi_1 (\oA' , 0) \quad \mbox{and} \quad h_* : \pi_1 (\oX , \ox) \longrightarrow \pi_1 (\oA , 0) \; .
\]
\end{theorem}

\begin{proof}
  The claims in {\bf a} and {\bf e} are formal consequences of Corollary \ref{t24} and the preceeding results on the categories $\eB_{\eX_{\eo}}$.
Let us now show {\bf b}. Since $\eB_{X_{\C_p}}$ and $\eB_{A_{\C_p}}$ are closed under exterior powers, it suffices to show that all line bundles from the categories $\eB$ lie in $\Pic^0_{\oX/\overline{\Q}_p} (\C_p)$ and in $\Pic^0_{A/K}({\C_p})$, respectively. 

Let us first consider the case of abelian varieties. We have to show that for any line bundle $L$ in $\eB_{\Ah_{\eo}}$ the generic fibre $L_{\C_p}$ lies in $\Pic^0_{A/K} (\C_p)$. By definition, there exists some $N \ge 1$ such that $N^* L_1$ is trivial on $\Ah_1$, where $L_1$ and $\Ah_1$ denote the reductions modulo $p$. If $k$ denotes the residue class field of $\eo$, this implies that $N^* L_k$ is trivial on $\Ah_k$. Since the N\'eron--Severi group of $\Ah_k$ is torsion free, $L_k$ lies in $\Pic^0_{\Ah_k/k}(k)$. Now $\Ah$ is projective as an abelian scheme over the normal base $\eo_K$ by \cite{Ray2}, Th\'eor\`eme XI 1.4, hence $\Pic^0_{\Ah / \eo_K}$ is an open subscheme of the scheme $\Pic_{\Ah / \eo_K}$ (see e.g. \cite{BLR}, section 8.4). Since the reduction of the point in $\Pic_{\Ah / \eo_K} (\eo)$ induced by $L$ is contained in $\Pic^0$, the generic fibre $L_{\C_p}$ is also contained in $\Pic^0$, whence our claim.

Now we consider the curve case. Let $\oeX$ be a finitely presented, flat and proper model of $\oX$ over $\eo_{\oK}$ and let $L$ be an object in $\eB_{\eX_{\eo}}$. By definition, there is a proper, finitely presented morphism $
\pi : \Yh \longrightarrow \oeX$
over $\eo_{\oK}$, such that its generic fibre $\pi_{\oK} : \Yh_{\oK} \to X_{\oK}$ is a finite \'etale covering, and such that for its special fibre $\pi_k : \Yh_k \to \oeX_k$ the pullback $\pi^*_k L_k$ is trivial. Since $\oeX$ is irreducible and $\pi (\Yh)$ is closed and contains the generic fibre, $\pi$ must be surjective. Hence $\pi_k$ is also surjective. Therefore every irreducible component $S$ of $\oeX_k$ is dominated by an irreducible component $T$ of $\Yh_k$. We endow $S$ and $T$ with the reduced structures.

The restriction of $\pi_k$ to $T$ is finite as a dominant, proper morphism of irreducible curves, see \cite{EGAII}, (7.4.4). Hence the restriction of $L_k$ to $S$ has degree zero. By  the degree formula in \cite{BLR}, Prop. 5, p. 239, the bundle $L_k$ has degree zero on $\oeX_k$. 

The Euler characteristic is locally constant in the fibres of $\eX$. To see this, use  \cite{EGAIV}, section 8 in order to reduce to a noetherian situation and apply \cite{EGAIII}, (7.9.4). Hence, by Riemann--Roch (see \cite{BLR}, Thm. 1, p. 238), the degree is also locally constant in the fibres of $\eX$, so that $L_{\C_p}$ has indeed degree zero. 

Let us now show {\bf c}. Let $\oeX$ be the smooth and proper model of $\oX$. We have $\Pic^0_{X_{\C_p} / \C_p} (\C_p) = \Pic^0_{\eX_{\eo}/ \eo} (\eo)$ since $\Pic^0_{\eX_{\eo} / \eo}$ is an abelian scheme. Together with Theorem \ref{t20} it follows that line bundles of degree zero lie in $\eB_{X_{\C_p}}$. Flat base change gives us
\[
\Ext^1_{\eX_{\eo}} (\Oh, \Oh) \otimes \Q = H^1 (\eX_{\eo} , \Oh) \otimes_{\eo} \C_p = H^1 (X_{\C_p} , \Oh) = \Ext^1_{X_{\C_p}} (\Oh, \Oh) \; .
\]
Thus, for a given extension
\[
0 \longrightarrow \Oh \longrightarrow \Oh(F) \longrightarrow \Oh \longrightarrow 0
\]
on $X_{\C_p}$, there exists an $N \ge 1$ such that $N \cdot [0 \to \Oh \to \Oh(F) \to \Oh \to 0]$ is the restriction of an extension class $[0 \to \Oh \to \Oh(E) \to \Oh \to 0]$ in $\Ext^1_{\eX_{\eo}} (\Oh, \Oh)$ for some vector bundle $E$. According to theorem \ref{t20}, we know that $E$ is an object of $\eB_{\eX_{\eo}}$. Consider the pushout diagram:
\[
\xymatrix{
0 \ar[r] & \Oh \ar[r] \ar[d]_{ N \, \wr} & \Oh(F) \ar[r] \ar[d]_{\wr} & \Oh \ar[r] \ar@{=}[d]& 0 \\
0 \ar[r] & \Oh \ar[r] & \Oh(F')  \ar[r] & \Oh  \ar[r] & 0 \; .
}
\]
On the one hand, we have
\[
\scriptstyle
[0 \to \Oh \to \Oh(F') \to \Oh] = N [0 \to \Oh \to \Oh(F) \to \Oh] = [0 \to \Oh \to j^* _{\eX_{\eo}} \Oh(E) \to \Oh \to 0] \; .
\]
Hence $F' \cong j^*_{\eX_{\eo}} E$. On the other hand we see that $F \cong F'$. By the definition of $\eB_{X_{\C_p}}$ it follows that $F$ lies in $\eB_{X_{\C_p}}$. Similarly for abelian varieties. 
\end{proof}

Parts {\bf d} and {\bf f} follow from corollary \ref{t24} and theorem \ref{t19} once we have shown the following two assertions:\\
{\bf i} For every finitely presented proper and flat model $\oeX'$ of $\oX'$ there exists a model $\oeX$ of $\oX$ with the same properties and an $\eo_{\oK}$-linear morphism $\tilde{f} : \oeX \to \oeX'$ such that the diagram
\[
\xymatrix{
\oeX \ar[d]_{\tilde{f}} & \oX \ar[l] \ar[d]^f \\
\oeX' & \oX' \ar[l]
}
\]
is commutative.\\
{\bf ii} Let $\oeX$ be a  smooth and proper model of $\oX$. Then there are commutative diagrams
\[
\xymatrix{
\oAh \ar[d]_{\tilde{g}} & \oA \ar[l] \ar[d]^g \\
\oAh' & \oA' \ar[l]
} \hspace*{2cm} 
\xymatrix{
\oeX \ar[d]_{\tilde{h}} & \oX \ar[l] \ar[d]^h \\
\oAh & \oA \ar[l]
}
\]
where $\tilde{g} (0) = 0$ and $\tilde{h} (\ox) = 0$ viewing $\ox$ as a section in $\oeX (\eo_{\oK}) = \oX (\oK)$. 

As for {\bf i} note that the image of $f$ is a closed and connected subset of $\oX'$. Hence $f$ is either surjective or its scheme theoretic image consists of a single reduced closed point $\ox'$ of $\oX'$. In the latter case we have $f = \ox' \verk \lambda$ where $\ox'$ is viewed as a $\oK$-valued point of $\oX'$ and $\lambda$ is the structural morphism of $\oX$ over $\oK$. It is now clear how to extend $f$ to arbitary finitely presented proper and flat models of $\oX$ resp. $\oX'$.

If $f$ is surjective it is quasi-finite by \cite{EGAII} 7.4.4 and hence finite by properness. Choose a finite extension $K / \Q_p$ such that $f$ is obtained by base extension from a finite morphism $f_0 : X \to X'$ of smooth projective curves over $K$. Moreover, we may assume that $\oeX'$ comes from a proper flat model $\eX'$ over $\eo_K$ of $X'$. Let $\eX$ be the normalization of the integral scheme $\eX'$ in the function field of $X$. Then according to \cite{EGAII} 6.3.9 we obtain a commutative diagram
\[
\xymatrix{
\eX \ar[d]_{\tilde{f}_0} & X \ar[l] \ar[d]^{f_0} \\
\eX' & X' \ar[l] \; .
}
\]
By construction of $\eX$ its fibre over $K$ is identified with $X$. Since $\eX'$ is excellent the morphism $\tilde{f}_0$ is finite by \cite{EGAIV}, 7.8.3 (vi). In particular, $\eX$ is proper over $\eo_K$. Being normal by construction and in particular reduced and irreducible it follows from \cite{Ha}, III, 9.7 that $\eX$ is also flat over $\eo_K$. Performing the base extension to $\eo_{\oK}$ we obtain the desired diagram in {\bf i}. 

As for {\bf ii} we first descend the situation to a finite extension $K$ of $\Q_p$ and then use the universal property of the N\'eron model. For example since $\eX$ is smooth over $\eo_K$, restriction induces a bijection
\[
\res : \Hom_{\eo_K} (\eX , \Ah) \silo \Hom_K (X , A) \; .
\]

\begin{cor}
  \label{t27}
If $A$ is an elliptic curve with good reduction then $\eB^{(\abb)}_{A_{\C_p}}$ (defined for abelian varieties in Def. 19.2) is a full subcategory of the category $\eB^{(c)}_{A_{\C_p}}$ of vector bundles defined for curves in Def. 19.1.

Besides, the two definitions for the functor $\rho$ on $\eB^{(\abb)}_{A_{\C_p}}$ in the curve case and in the case of abelian varieties coincide.

$\eB^{(\abb)}_{A_{\C_p}}$ contains all vector bundles $E$ of the form $E = E_1 \oplus \ldots \oplus E_r$, where all $E_i$ are indecomposable bundles of degree zero.
\end{cor}

\begin{proof}
  The first two assertions are immediate. For the last one, we use Atiyah's classification of vector bundles on elliptic curves. 

By \cite{At}, Theorem 5, p. 432 every indecomposable vector bundle $E$ of degree zero on $A_{\C_p}$ is of the form $L \otimes F_r$ for a line bundle $L$ of degree zero and a vector bundle $F_r$ sitting in an exact sequence
\[
0 \longrightarrow \Oh \longrightarrow F_r \longrightarrow F_{r{-1}} \longrightarrow 0 \; .
\]
Hence $F_r$ is a successive extension of trivial line bundles and therefore by Theorem \ref{t26} contained in $\eB_{A_{\C_p}}^{(\abb)}$. Since Theorem \ref{t26} also says that our category is closed under tensor products and direct sums and contains all line bundles of degree zero, our claim follows.
\end{proof}

We expect that in fact the categories $\eB^{(\abb)}_{A_{\C_p}}$ and $\eB^{(c)}_{A_{\C_p}}$ coincide and that they contain precisely the vector bundles whose indecomposable components have degree zero. It should be possible to calculate the representations $\rho_{F_r}$ explicitly at least if $A$ has good ordinary reduction.

\begin{theorem} \label{t28}
  Let $X / K$ be a smooth, projective curve over $K$ and $A / K$ an abelian variety with good reduction.\\
{\bf a} The natural actions of $G_K$ on $\Vec_{X_{\C_p}}$ and $\Vec_{A_{\C_p}}$ induce actions on $\eB_{X_{\C_p}}$ and $\eB_{A_{\C_p}}$. The functors
\[
\rho : \eB_{X_{\C_p}} \longrightarrow \Rep_{\pi_1 (\oX , \ox)} (\C_p) \quad \mbox{and} \quad \rho : \eB_{A_{\C_p}} \longrightarrow \Rep_{\pi_1 (\oA, 0)} (\C_p)
\]
commute with the $G_K$-actions on the categories. In particular we have
\[
\rho_{\,\!^{\sigma}E} = \, ^{\sigma} \rho_E \quad \mbox{for} \; E \; \mbox{in} \; \eB_{X_{\C_p}} \quad \mbox{resp.} \quad \eB_{A_{\C_p}} \; .
\]
{\bf b} Assume additionally that $X$ has good reduction.\\
Recall the homomorphism (\ref{eq:10}) resp. (\ref{eq:20}):
\[
\alpha : \Pic^0_{A/K} (\C_p) = \hA (\C_p) \longrightarrow \Hom_c (TA , \C^*_p) = \Hom_c (TA , \eo^*)
\]
and
\[
\alpha : \Pic^0_{X/K} (\C_p) \longrightarrow \Hom_c (\pi_1 (\oX , \ox) , \C^*_p) = \Hom_c (\pi_1 (\oX , \ox) , \eo^*) \; .
\]
Then for any line bundle $L$ of degree zero on $X_{\C_p}$ resp. $A_{\C_p}$ we have:
\[
\alpha (L) = \rho_L \; .
\]
\end{theorem}

\begin{proof}
  Part {\bf a} is a formal consequence of Corollary \ref{t24} and theorem \ref{t18}. As for {\bf b} let $h : X \to A$ be the map of $X$ into its Albanese variety $A = \Alb_{X / K}$ with $h (x) = 0$. By theorem \ref{t26} {\bf c} line bundles of degree zero lie in $\eB_{A_{\C_p}}$ resp. $\eB_{X_{\C_p}}$. The third diagram in theorem \ref{t26} {\bf f} therefore induces a commutative diagram of abelian groups
\[
\xymatrix{\Pic^0_{A/K} (\C_p) \ar[r]^{\overset{h^*}{\sim}} \ar[d]_{\rho} & \Pic^0_{X / K} (\C_p) \ar[d]^{\rho} \\
\Hom_c (\pi_1 (\oA , 0) , \C^*_p) \ar[r]^{\overset{H}{\sim}} & \Hom_c (\pi_1 (\oX , \ox) , \C^*_p)}
\]
where $H$ is induced by $h_* : \pi_1 (\oX , \ox) \to \pi_1 (\oA , 0)$. \\
On the other hand by the very definition (\ref{eq:20}) of $\alpha$ for curves via the map $\alpha$ in (\ref{eq:10}) for its Albanese variety $A$, the same diagram commutes if we replace the $\rho$'s by $\alpha$'s. Hence it suffices to prove assertion {\bf b} in the case of abelian varieties.

Let $\ha$ be a point in $\hA(\C_p) = \hAh(\eo)$ and denote by $L$ be the corresponding  line bundle on $\Ah_\eo$. If $N$ is big enough,  $\ha_n \in \hAh_N (\eo_n)$ is an $N$-torsion point corresponding to the reduction  $L_n$. Now $\alpha(\ha)$ is defined as the limit of the maps $TA \to \Ah_N (\eo_n) \xrightarrow{\varphi_n} \eo^{\times}_n$, where 
$\varphi_n$ is the image of $\ha_n$ under the Cartier duality homomorphism $\hAh_N (\eo_n) \to \Hom (\Ah_{n,N} , \Ge_{m, \eo_n})$.

Recall that we identified $\hAh$ with $\uExt^1 (\Ah , \Ge_m)$. Hence the $\Ge_m$-torsor $\tilde{L} = L \ohne \{ \mathrm{zero\,section} \}$ associated to $L$ can be endowed with the structure of an extension of $\Ah$ by $\Ge_m$. Besides, with this identification the inclusion $\uHom (\Ah_N , \Ge_m) \simeq \hAh_N \hookrightarrow \hAh$ is given by pushout with respect to the exact sequence $0 \to \Ah_N \to \Ah \xrightarrow{N} \Ah \to 0$. Since  $\varphi_n : \Ah_{n,N} \to \Ge_{m, \eo_n}$ corresponds to $\ha_n$, and hence to $\tilde{L}_n$, the extension $\tilde{L}_n$ is given by pushout with respect to $\varphi_n$:
\begin{equation}
  \label{eq:49}
\vcenter{
\xymatrix{
0 \ar[r] & \Ah_{n,N} \ar[r] \ar[d]_{\varphi_n} & \Ah_n \ar[r]^N \ar[d]_{f_n} & \Ah_n \ar[r] \ar@{=}[d] & 0 \\
0 \ar[r] & \Ge_{m, \eo_n} \ar[r]^j & \tilde{L}_n \ar[r] & \Ah_n \ar[r] & 0
}
}
\end{equation}
Then $s_n = (f_n , \id) : \Ah_n \to \tilde{L}_n \times_{\Ah_n} \Ah_n = N^* \tilde{L}_n$ is a section of the trivial extension $N^* \tilde{L}_n$, hence it gives a trivialization of $N^* \tilde{L}_n$. 
Let us denote by $e$ the unit section of $\Ah_n$. Then the map $j$ induces an isomorphism $e^* \tilde{L}_n \simeq \Ge_{m, \eo_n}$, which maps $e^* s_n$ to $1 \in \Ge_{m, \eo_n}(\eo_n) = \eo_n^*$. Now $\rho_{L}$ is defined as the limit of the maps $TA \longrightarrow \Ah_N (\eo_n) \xrightarrow{\psi_n} \eo_n^*$, where $\psi_n$ is obtained by transferring the natural action 
of $a_n \in \Ah_N(\eo_n)$ on $\Gamma(\Ah_n, N^* L_n)$ to an action on $\Gamma(\Spec \eo_n , e^* L_n) \simeq \eo_n$. Hence $\psi_n(a_n)$ can be calculated as follows: We take $1 \in \eo_n^*
{\simeq} \Gamma(\Spec \eo_n , e^* \tilde{L}_n)$. We have seen that its lift in $\Gamma(\Ah_n, N^* \tilde{L}_n)$ is equal to $s_n$. Applying the translation map $\tau_{a_n}$, we get the section
$s_n \verk \tau_{a_n}$. The pullback of this section via $e^*$ is equal to $f_n(a_n) \in 
\Gamma(\Spec \eo_n , e^* \tilde{L}_n)$, which corresponds to $\varphi_n(a_n) \in \eo_n^*$ under the isomorphism
$j$. 
Hence we have shown $\alpha (L) = \rho_L$.
\end{proof}

According to sections 3 and 4 the map $\Lie \alpha$ is one of the maps occuring in the Hodge--Tate decomposition on the first cohomology. The next theorem gives another relation of our constructions with this Hodge--Tate map.

\begin{theorem}
  \label{t29} 
Consider $X / K$ and $A / K$ with good reductions as above. We write $\Ext^1_{\eB_{X_{\C_p}}} (\Oh, \Oh)$ and $\Ext^1_{\eB_{A_{\C_p}}} (\Oh, \Oh)$ for the Yoneda groups of isomorphy classes of extensions $0 \to \Oh \to \Oh (E) \to \Oh \to 0$, where $E$ lies in $\eB_{X_{\C_p}}$, respectively in $\eB_{A_{\C_p}}$. Then there are commutative diagrams:
\[
\xymatrix{
\Ext^1_{\eB_{X_{\C_p}}} (\Oh, \Oh) \ar[r]^-{\rho_*} \ar@{=}[d] & \Ext^1_{\Rep_{\pi_1 (\oX, \ox)} (\C_p)} (\C_p , \C_p) \ar@{=}[d] \\
H^1 (X , \Oh) \otimes_K \C_p \ar[r]^{\rm Hodge-Tate} & H^1_{\et} (\oX , \Q_p) \otimes \C_p
}
\]
and
\[
\xymatrix{
\Ext^1_{\eB_{A_{\C_p}}} (\Oh , \Oh) \ar[r]^-{\rho_*} \ar@{=}[d] & \Ext^1_{\Rep_{\pi_1 (\oA, 0)}(\C_p)} (\C_p , \C_p) \ar@{=}[d] \\
H^1 (A , \Oh) \otimes_K \C_p \ar[r]^{\theta^*_A {\rm from\,(\ref{eq:15})}} & H^1_{\et} (\oA , \Q_p) \otimes \C_p \; .
}
\]
\end{theorem}

Before we give the proof let us make two remarks.

\begin{rems}
  {\bf 1} The second diagram in theorem \ref{t29} gives the following novel construction of the Hodge Tate map $\theta^*_A$. Consider a class $c$ in $H^1 (\Ah, \Oh) \otimes_{\eo_K} \eo$. It can be viewed as an extension 
\[
0 \to \Oh \to \Oh (E) \to \Oh \to 0
\]
of locally free sheaves on $\Ah_{\eo}$. The bundle $E$ lies in $\eB_{\Ah_{\eo}}$. Hence, for every $n \ge 1$ there is some $N \ge 1$ in fact, $N = p^n$ will do, such that $\pi^*_{N,n} E_n$ is the trivial rank two bundle on $\Ah_n$. Here $\pi_N : \Ah_{\eo} \to \Ah_{\eo}$ is the $N$-multiplication map, and $\pi_{N,n}$ its reduction mod $p^n$. The short exact sequence $0 \to \eo_n \to E_{n,0} \to \eo_n \to 0$ of fibres along the zero section of $\Ah_n$ becomes $TA$-equivariant if $TA$ acts trivially on the $\eo_n$'s and via the projection $TA \to A_N (\oK)$ and the isomorphism
\[
\Gamma (\Ah_n , \pi^*_{N,n} E_n) \overset{\res}{\silo} E_{n,0}
\]
on $E_{n,0}$. Passing to projective limits, we get a short exact sequence \\
$0 \to \eo \xrightarrow{i} E_0 \xrightarrow{q} \eo \to 0$ of $TA$-modules. Set $g_1 = i (1)$ and choose $g_2 \in E_0$ such that $q (g_2) = 1$. Then $E_0$ is a free $\eo$-module on $g_1$ and $g_2$, and the action of $\gamma \in TA$ on $E_0$ is given in terms of the basis $g_1, g_2$ by a matrix of the form $\left( 
  \begin{smallmatrix}
    1 & \beta (\gamma) \\ 0 & 1
  \end{smallmatrix} \right)$ where $ \beta : TA \to \eo$ is a continuous homomorphism. Note that $\beta$ does not depend on the choice of $e_1$. Viewing $\beta$ as an element of $H^1_{\et} (\oA , \Z_p) \otimes \eo$ we have $\theta^*_A (c) = \beta$.

{\bf 2} Another interpretation of theorem \ref{t29} and the preceeding remark makes the analogy clearer between these statements and the fact that the following diagram commutes
\[
\def\objectstyle{\scriptstyle}
\def\labelstyle{\scriptstyle}
\xymatrix{
\Pic^0_{A/K}(\C_p) \ar@{=}[r] & H^1 (A_{\C_p} , \Oh^*)^0 \ar[r]^{\log} \ar[d]_{\alpha} & H^1 (A_{\C_p} , \Oh) \ar@{=}[r] \ar[d]^{\theta^*_A} & \Lie \Pic^0_{A/K} (\C_p) \ar[d]^{\Lie \alpha} \\
 & \Hom_c (\pi_1 (\oX , \ox) , \C^*_p) \ar[r]^{\log_*} & \Hom_c (\pi_1 (\oX , \ox) , \C_p) \ar@{=}[r] & \Lie \Hom_c (\pi_1 (\oX , \ox) , \C^*_p) \; .
}
\]

Namely, consider the Lie group $U = \left\{ \left( 
    \begin{smallmatrix}
      1 & * \\ 0 & 1
    \end{smallmatrix} \right) \right\} \subset \GL_2$. It is isomorphic to $\Ge_a$, as is its Lie algebra $\Lie U$. Under these identifications the $p$-adic logarithm map becomes the identity and theorem \ref{t29} asserts that the following diagram is commutative:
\[
\def\objectstyle{\scriptstyle}
\def\labelstyle{\scriptstyle}
\xymatrix{
\Ext^1_{\eB_{A_{\C_p}}} (\Oh, \Oh) \ar@{=}[r] & H^1 (A_{\C_p} , U (\Oh)) \ar[r]^{\overset{\log = \id}{\sim}} \ar[d]_{\rho_*} & H^1 (A_{\C_p} , \Lie U (\Oh)) \ar@{=}[r] & H^1 (A_{\C_p} , \Oh) \ar[d]^{\theta^*_A} \\
 & \Hom_c (\pi_1 (\oA , 0) , U (\C_p)) \ar[r]^{\log_* = \id} & \Hom_c (\pi_1 (\oA , 0) , \Lie U (\C_p)) \ar@{=}[r] & H^1_{\et} (\oA , \Q_p) \otimes \C_p \; .
}
\]
Thus in the first diagram the underlying group is $\Ge_m$, whereas in the second it is $U \cong \Ge_a$. Ony may wonder whether there are generalizations to non-commutative groups $U$, where the logarithm maps are no longer homomorphisms.
\end{rems}

\begin{proofof}
  {\bf theorem \ref{t29}} By theorem \ref{t26} {\bf c} the categories $\eB_{X_{\C_p}}$ and $\eB_{A_{\C_p}}$ contain all vector bundles which are extensions of the trivial line bundle by itself. Hence we have
\[
\Ext^1_{\eB_{X_{\C_p}}} (\Oh, \Oh) = \Ext^1_{\Vec_{X_{\C_p}}} (\Oh , \Oh) = H^1 (X_{\C_p} , \Oh) = H^1 (X, \Oh) \otimes_K \C_p
\]
and similarly over $A_{\C_p}$. 

Set $\Ch = \Rep_{\pi_1 (\oX , \ox)} (\C_p)$. Consider the natural isomorphism
\begin{equation}
  \label{eq:50}
  \Ext^1_{\Ch} (\C_p , \C_p) \silo \Hom_c (\pi_1 (\oX , \ox) , \C_p)
\end{equation}
defined as follows. For an extension $0 \to \C_p \xrightarrow{i} V \xrightarrow{\varepsilon} \C_p \to 0$ define a continuous homomorphism $\psi : \pi_1 (\oX , \ox) \to \C_p$ by setting $\psi (\gamma) = i^{-1} (\gamma \cdot v - v)$ where $v \in V$ is such that $\varepsilon (v) =1 $. Composing with the isomorphisms
\[
\Hom_c (\pi_1 (\oX , \ox) , \C_p) = \Hom_c (\pi_1 (\oX , \ox) , \Q_p) \otimes_{\Q_p} \C_p = H^1_{\et} (\oX , \Q_p) \otimes_{\Q_p} \C_p
\]
we get the natural isomorphism
\[
\Ext^1_{\Ch} (\C_p , \C_p) \silo H^1_{\et} (\oX , \Q_p) \otimes \C_p
\]
used in the statement of the theorem. Similarly for abelian varieties. We now reduce the curve case to the case of abelian varieties. Let $h : X \to A$ be the map of $X$ into its Albanese variety with $h (x) = 0$. Since all functors involved are exact, the third diagram in theorem \ref{t26} induces the middle sqare in the commutative diagram:
\[
\xymatrix{
H^1 (A_{\C_p} , \Oh) \ar[r]^{\overset{h^*}{\sim}} \ar@{=}[d] & H^1 (X_{\C_p} , \Oh) \ar@{=}[d] \\
\Ext^1_{\eB_{A_{\C_p}}} (\Oh , \Oh) \ar[r]^{h^*} \ar[d]_{\rho_*} & \Ext^1_{\eB_{X_{\C_p}}} (\Oh, \Oh) \ar[d]^{\rho^*} \\
\Ext^1_{\Rep_{\pi_1 (\oA , 0)}(\C_p)} (\C_p , \C_p) \ar@{=}[d] \ar[r]^{H_*} & \Ext^1_{\Rep_{\pi_1 (\oX , \ox)}(\C_p)} (\C_p , \C_p) \ar@{=}[d] \\
H^1 (\oA , \Q_p) \otimes \C_p \ar[r]^{\overset{h^*}{\sim}} & H^1 (\oX , \Q_p) \otimes \C_p \; .
}
\]
The Hodge--Tate decomposition is functorial, moreover the maps $h^*$ on top and bottom of this diagram are known to be isomorphisms. It therefore suffices to verify the assertion of the theorem in the case of abelian varieties.

Since $H^1 (\Ah_{\eo} , \Oh) \otimes_{\eo} \C_p = H^1 (A_{\C_p}, \Oh)$, we find for every element in \\
$\Ext^1_{\eB_{A_{\C_p}}} (\Oh , \Oh)$ a $p$-power multiple lying in $\Ext^1_{\eB_{\Ah_{\eo}}} (\Oh , \Oh)$. Since $\rho_*$ is a homomorphism between Yoneda $\Ext$-groups, it suffices to show that
\[
\xymatrix{
\Ext^1_{\eB_{\Ah_{\eo}}} (\Oh , \Oh) \ar[r]^{\rho_*} \ar@{=}[d] & \Ext^1_{\Rep_{TA} (\eo)} (\eo , \eo) \ar[d]^{\wr} \\
H^1 (\Ah_{\eo} , \Oh) \ar[r]^{\theta^*_A} & \Hom_c(T_p A , \eo)
}
\]
commutes. Consider an extension $0 \to \Oh \stackrel{i}{\rightarrow} \Fh \to \Oh \to 0$ on $\Ah_\eo$, where $\Fh = \Oh(E)$ for some $E$  in $\eB_{\Ah_{\eo}}$. 

Let $\Ih = \uIso_{\Ext} (\Oh^2 , \Fh)$ be the (Zariski) sheaf on $\Ah_{\eo}$ which associates to an open subset $U \subseteq \Ah_{\eo}$ the set of isomorphisms $\varphi : \Oh^2_U \to \Fh_U$ of extensions, i.e. such that
\[
\xymatrix{
0 \ar[r] & \Oh_U \ar[r] \ar@{=}[d] & \Oh^2_U \ar[r]^p \ar[d]_{\varphi} & \Oh_U \ar[r] \ar@{=}[d] & 0 \\
0 \ar[r] & \Oh_U \ar[r]^i & \Fh_U \ar[r] & \Oh_U \ar[r] & 0
}
\]
commutes. Then $c \in \Ge_a(U) = \Gamma (U, \Oh)$ acts in a natural way on $\Ih (U)$ by mapping
\[
\varphi \longmapsto \varphi + i \verk f_c \verk p \; , 
\]
where $f_c : \Oh_U \to \Oh_U$ is multiplication by $c$.

The class $[\Ih]$ of $\Ih$ in $H^1 (\Ah_{\eo} , \Oh) = H^1_{\Zar} (\Ah_{\eo} , \Ge_a)$ coincides with the image of the extension class given by $\Fh$ under the isomorphism 
\[
\Ext^1 (\Oh , \Oh) \silo H^1 (\Ah_{\eo} , \Oh) \; .
\]
Note that the association
\begin{equation}
  \label{eq:51}
(T \xrightarrow{t} \Ah_{\eo}) \longmapsto \uIso_{\Ext} (\Oh^2_T , t^* \Fh)
\end{equation}
also defines a sheaf on the  flat site over $\Ah_{\eo}$, which is represented by a $\Ge_a$-torsor $Z \to \Ah_{\eo}$. 

We have $H^1_{\Zar} (\Ah_{\eo} , \Ge_a) = H^1_{\rm fppf} (\Ah_{\eo} , \Ge_a) = \Ext^1 (\Ah_{\eo} , \Ge_a)$, so that $Z$ can be endowed with a group structure sitting in an extension of $\Ah_{\eo}$ by $\Ge_a$:
\[
0 \longrightarrow \Ge_a \stackrel{j}{\longrightarrow} Z \longrightarrow \Ah_{\eo} \longrightarrow 0 \; .
\]
Hence there is a homomorphism $h : \omega_{\hAh_{\eo}} \to \Ge_a$ such that $Z$ is the pushout of the universal vectorial extension:
\[
\xymatrix{
0 \ar[r] & \omega_{\hAh_{\eo}} \ar[r] \ar[d]_h & E_{\eo} \ar[r] \ar[d] & \Ah_{\eo} \ar[r] \ar@{=}[d] & 0 \\
0 \ar[r] & \Ge_{a} \ar[r] & Z \ar[r] & \Ah_{\eo} \ar[r] & 0 \; .
}
\]
Recall that $\theta_A : T_p A \to \omega_{\hAh} (\eo)$ is defined as the limit of the maps
\[
\theta_{A,n} : \Ah_{p^n} (\eo_{\oK}) \longrightarrow \omega_{\hAh} (\eo) / p^n \omega_{\hAh} (\eo)
\]
associating to $a \in \Ah_{p^n} (\eo_{\oK})$ the class of $p^n b$ for an arbitrary preimage $b \in E (\eo)$ of $a$.

By \cite{Ma-Me}, chapter I, the natural isomorphism
\[
\Hom_{\eo}(\omega_{\hAh} (\eo) , \eo) \stackrel{\sim}{\longrightarrow} \mbox{Lie}(\hAh_{\eo})(\eo) \stackrel{\sim}{\longrightarrow}  H^1 (\Ah_{\eo} , \Oh) = \mbox{Ext}^1(\Ah_\eo, \Ge_{a})
\]
sends a map to the corresponding pushout of $E_{\eo}$.  Hence $\theta^*_A : H^1 (\Ah_{\eo} , \Oh) \to \Hom_c (TA , \eo)$ maps the extension class of $\Fh$ to the limit of the maps
\[
\Ah_{p^n} (\eo_{\oK}) \xrightarrow{\theta_{A,n}} \omega_{\hAh} (\eo) / p^n \omega_{\hAh} (\eo) \xrightarrow{h_n} \eo_n \; , 
\]
where $h_n$ is induced by $h$.

This map can also be described as follows: For $a \in \Ah_{p^n} (\eo_{\oK}) \subseteq \Ah_{p^n} (\eo)$ choose a preimage $z \in Z (\eo)$. Then
\[
h_n \verk \theta_{A,n} (a) = \mbox{class of} \; p^n z \; \mbox{in} \; \Ge_a(\eo) / p^n \Ge_a(\eo) = \eo_n \; .
\]
Set $Z_n = Z \otimes_{\eo} \eo_n$. Let $\pi_{p^n}$ denote multiplication by $p^n$ on $\Ah_n$. We have seen in the proof of Theorem \ref{t20} that $\pi_{p^n}^* \Fh_n$ is a trivial extension. Hence $\pi_{p^n}^* Z_n$ is trivial in $\Ext(\Ah_n, \Ge_{a, \eo_n})$, and there is a splitting   $r : \Ah_n \to \pi_{p^n}^* Z_n$ of the extension
\[
0 \longrightarrow \Ge_{a , \eo_n} \longrightarrow \pi_{p^n}^* Z_n \longrightarrow \Ah_n \longrightarrow 0
\]
over $\eo_n$. Let $F: \pi_{p^n}^* Z_n \rightarrow Z_n$ denote the projection, and denote by 
$a_n \in \Ah_n(\eo_n)$ the point induced by $a$. Then $F(r(a_n))$ projects to zero in $\Ah_n$, hence it is equal to $j(c)$ for some $c \in \Ge_a(\eo_n) = \eo_n$. Since for any $a' \in
 \Ah_n(\eo_n)$ with $p^n a' = a_n $ the point $F(r(a'))$ is a preimage of $a_n$, we have $h_n \verk \theta_{A,n} (a) = c \in \eo_n$. Besides,  $Z$ represents the functor (\ref{eq:51}), so that the map $r: \Ah_n \rightarrow \pi_{p^n}^* Z_n$  corresponds to a trivialization
\[
\xymatrix{
0 \ar[r] & \Oh_{\Ah_n} \ar[r] \ar@{=}[d] & \Oh^2_{\Ah_n} \ar[r] \ar[d]_{\varphi} & \Oh_{\Ah_n} \ar[r] \ar@{=}[d] & 0 \\
0 \ar[r] & \Oh_{\Ah_n} \ar[r] & \pi_{p^n}^* \Fh_n \ar[r] & \Oh_{\Ah_n} \ar[r] & 0 \; .
}
\]
The point $F(r(a_n))$ in the kernel of $Z_n(\eo_n) \rightarrow \Ah_n(\eo_n)$ corresponds to the
trivialization
\[
\alpha: \Oh^2_{\spec \eo_n} \stackrel{a_n^* \varphi}{\longrightarrow} a_n^* \pi_{p^n}^* \Fh  \stackrel{\sim}{\longrightarrow} 0^* \Fh,
\]
where $0$ is the zero element in $\Ah_n(\eo_n)$. 
Besides, the trivialization 
\[
\beta: \Oh^2_{\spec \eo_n} \stackrel{0^* \varphi}{\longrightarrow} 0^* \pi_{p^n}^* \Fh  \stackrel{\sim}{\longrightarrow} 0^* \Fh
\]
is given by the zero element in $Z_n$. By definition, $\alpha = \beta + i \circ f_c \circ p$.
If we denote the canonical basis of $\Gamma( \Ah_n, \Oh^2_{\Ah_n})$ by $e_1,e_2$, and the induced basis of $\Gamma(\Ah_n, \pi_{p^n}^* \Fh)$ by  $f_1 , f_2$, it follows that
$ a_n^* f_2 - 0^* f_2 = i(c)$. 

On the other hand, the image of $E$ under 
\[
\rho_* : \Ext^1_{\eB_{\Ah_{\eo}}} (\Oh , \Oh) \longrightarrow \Ext^1_{\Rep_{TA} (\eo)} (\eo , \eo) = \Hom_c(T_p A , \eo)\]
is the  homomorphism $\gamma: T_pA \rightarrow \eo$, such that $\gamma \mod p^n$ maps the $p^n$-torsion point
$a_n \in \Ah_n(\eo_n)$ to the element 
\[i^{-1}( 0^*(\tau_{a_n}^* f_2)- 0^* f_2) = i^{-1}(a_n^* f_2 - 0^*f_2),\] 
where $\tau_{a_n}$ is translation by $a_n$, and hence to $c = h_n \circ \theta_{A,n}(a)$.
This proves our claim.
\end{proofof}
%\newpage
%\input{sec7}
\section{Open questions}
In this section we discuss some further questions raised by the preceeding constructions. Let $X$ be a smooth projective curve over $K$ with good reduction and let $\ox \in X (\oK)$ be a geometric point.

{\bf a}
\rm  Is the functor $\rho : \eB_{X_{\C_p}} \to \Rep_{\pi_1 (\oX , \ox)} (\C_p)$ fully faithful? 

Since $\eB_{X_{\C_p}}$ is closed under internal homs, $\rho$ is faithful if and only if for all bundles $E$ in $\eB_{X_{\C_p}}$ the natural map $H^0 (X_{\C_p} ,E) \to E_{x_{\C_p}}$ is injective. Note that for line bundles of degree zero this is true. It is also known to be true for stable bundles of degree zero. Suppose that $X$ has a smooth and proper model $\eX$ over $\eo_K$. Then it can be proved that at least the restriction of $\rho$ to the essential image of $\eB_{\eX_{\eo}}$ in $\eB_{X_{\C_p}}$ is faithful. 

Similarly, $\rho$ is full if and only if the natural map $H^0 (X_{\C_p} , E) \to E^{\pi_1 (\oX , \ox)}_{x_{\C_p}}$ is surjective, where $\pi_1 (\oX , \ox)$ acts via $\rho_E$ on $E_{x_{\C_p}}$. Again this is true for line bundles $L$ of degree zero since $\rho_L = \alpha (L)$ is non-trivial whenever $L$ is non-trivial by theorem \ref{t4}. The proof used $p$-adic Hodge theory.\\

{\bf b}
\rm  What is the essential image of the functor $\rho$? In theorem \ref{t12} as a consequence of theorem \ref{t7} this was determined on the full subcategory of $\eB_{X_{\C_p}}$ of line bundles of degree zero. 

In particular $\rho$ is not essentially surjective. This leads to the next question.\\

{\bf c}
\rm  Is there an additional structure, a kind of $p$-adic Higgs field $\phi$ on a vector bundle, such that to certain pairs $(E , \phi)$ one can canonically attach a $p$-adic representation? For $\phi = 0$ the theory should reduce to the one of the present paper. For line bundles, a classical Higgs field is just a one-form. Via Hodge--Tate theory it gives rise to an element of $\Hom_c (\pi_1 (\oX , \ox) , \C_p (1))$. What one really needs is an element of $\Hom_c (\pi_1 (\oX , \ox) , \C^*_p)$ though, so that some modification is required.\\
Ideally one could hope that a suitable category of vector bundles plus field on $X_{\C_p}$ would be equivalent to $\Rep_{\pi_1 (\oX , \ox)} (\C_p)$. \\

The classical theory of vector bundles on compact Riemann surfaces suggests the following question:

{\bf d}
\rm  Is a vector bundle $E$ on $X_{\C_p}$ in $\eB_{X_{\C_p}}$ if and only if all its indecomposable components have degree zero? If so, then in particular all stable bundles of degree zero would lie in $\eB_{X_{\C_p}}$. Are the corresponding $p$-adic representations of $\pi_1 (\oX , \ox)$ irreducible?\\

Further questions come from Hodge--Tate theory:

{\bf e}
\rm  For $E$ in $\eB_{X_{\C_p}}$ we may view $\rho_E$ as a sheaf on $\oX$. Is there a Galois-equivariant Hodge--Tate decomposition:
\[
H^1_{\et} (\oX , \rho_E) = H^1_{\cont} (\pi_1 (\oX , \ox) , \rho_E) \cong H^1 (X_{\C_p} , \Oh (E)) \oplus H^0 (X_{\C_p} , \Omega^1 (E)) (-1)?
\]
The present theory is a $\C_p$-theory. Is there a refinement to a $\eB_{dR}$-theory?\\

{\bf f}
\rm  Let $\pi : \eX' \to \eX$ be a smooth and proper morphism of smooth and proper $\eo_K$-schemes. Assume that $X = \eX_K$ is a curve. The relative de Rham cohomology $E = R^i \pi_{K*} (\Omega^{\hullet}_{X'/X})$ is a vector bundle on $X$ with the integrable Gauss--Manin connection. Is $E$ in $\eB_{X_{\C_p}}$ and if so, is
\[
\rho_E \cong R^i \pi_{K*} (\Q_p)_{\ox} \otimes_{\Q_p} \C_p?
\]
Note that $R^i \pi_{K*} (\Q_p)$ is a smooth $\Q_p$-sheaf on $\oX_{\et}$. The group $\pi_1 (\oX , \ox)$ therefore acts on its stalk in $\ox$. 

%\newpage
%\input{lit}

\end{document}